\documentclass[12pt,letterpaper]{amsart}%
\usepackage{amsmath}
\usepackage{txfonts}
\usepackage{graphicx}
\usepackage{latexsym,amsmath,amsthm,amssymb,amsfonts,mathrsfs}
\usepackage[utf8]{inputenc}
\usepackage{xcolor}
\usepackage{lipsum}%
\usepackage{tikz}
\usepackage{xcolor}
\usepackage{comment}
\usepackage{newtx}
\usepackage{relsize}

\usepackage{geometry}
\providecommand{\U}[1]{\protect\rule{.1in}{.1in}}
\numberwithin{equation}{section}
\newtheorem{theorem}{Theorem}[section]
\newtheorem{lem}[theorem]{Lemma}
\newtheorem{thm}[theorem]{Theorem}
\newtheorem{pro}[theorem]{Proposition}
\newtheorem{cor}[theorem]{Corollary}
\newtheorem{defi}[theorem]{Definition}
\newtheorem{rem}[theorem]{Remark}

\def\S2{\mathbb{S}^2}
\def\s{\,\,\,\,}

\def\endproof{$\hfill\Box$}
\def\R{\mathbb{R}}
\def\C{\mathbb{C}}
\def\vol{\mbox{vol}}
\def\S{\mathcal{S}}
\def\M{\mathcal{M}}
\def\K{\mathbb{K}}
\def\co{\operatorname{cont}}
\def\vol{\mathrm{Area}\,}
\def\diam{\mathrm{diam}\,}

\def\avint{\mathop{\rlap{\ \,--}\int}\nolimits}

\allowdisplaybreaks

\begin{document}

\title[Uniform convergence of metrics]{Uniform Convergence of Metrics on Alexandrov Surfaces with Bounded Integral Curvature}
\author{Jingyi Chen and Yuxiang Li}
\address{Department of Mathematics, 
The University of British Columbia, Vancouver, BC, Canada}
\email{jychen@math.ubc.ca}
\address{Department of Mathematics, Tsinghua University, Beijing, China}
\email{liyuxiang@mail.tsinghua.edu.cn}

\subjclass[2020]{Primary 53C45; Secondary 53C43, 53C22.}

\thanks{Chen is partially supported by NSERC Discovery Grant, No. 22R80062 and Li is partially supported by  NSFC, Grant No. 12141103}

%\begin{document}

\begin{abstract}
We prove uniform convergence of metrics $g_k$ on a closed surface with bounded integral curvature (measure) in the sense of A.D. Alexandrov, under the assumption that the curvature measures $\K_{g_k}=\mu^1_k-\mu^2_k$, where $\mu^1_k,\mu^2_k$ are nonnegative Radon measures converging weakly to measures $\mu^1,\mu^2$ respectively, and $\mu^1$ is less than $2\pi$ at each point (no cusps). This is the global version of Yu. G. Reshetnyak's well-known result on uniform convergence of metrics on a domain in $\C$, and answers affirmatively the open question on the metric convergence on a closed surface. 
We also give an analytic proof of the fact that a (singular) metric $g=e^{2u}g_0$ with bounded integral curvature 
%in the sense of Alexandrov 
on a closed Riemannian surface $(\Sigma,g_0)$ can be approximated by smooth metrics in the fixed conformal class $[g_0]$. % in terms of distance functions, curvature measures and conformal factors.
 Results on a closed surface with varying conformal classes and on complete noncompact surfaces are obtained as well.
\end{abstract}

%\clearpage

\maketitle

\vspace{-1cm}
\tableofcontents

\section{Introduction}

In the 1960s, Reshetnyak developed an analytic approach centred around his theory of subharmonic metrics to study Alexandrov surfaces. He showed \cite{R2} that a Radon measure $\mu$ on a surface induces a distance function $d$ with  $d(x,y)<+\infty$ if $\mu(\{x\})<2\pi$ and $\mu(\{y\})<2\pi$, and proved a fundamental convergence result in a relatively compact domain in $\C$ \cite[Theorem 7.3.1]{R},  \cite[Theorem III]{R1}.

\begin{thm}[Reshetnyak]\label{Resh}
Let $\Omega\subset\C$ be a relatively compact domain with piecewise smooth boundary and let $\mu^1_n, \mu^{2}_n$ be nonnegative Radon measures supported in $\Omega$ weakly converging to Radon measures $\mu^1$ and $\mu^2$ respectively as $n\to\infty$.  Suppose $\mu_n=\mu_n^1-\mu_n^2,\mu=\mu^1-\mu^2$ and $g_n=e^{2u_n}|dz|^2,g=e^{2u}|dz|^2$ where 
$$
u_n(z)=-\frac{1}{\pi}\int_\C\log|z-\zeta|d\mu_{n}(\zeta),\s u(z) =-\frac{1}{\pi}\int_\C\log|z-\zeta|d\mu(\zeta).
$$
If $\mu^1(\{z\})<2\pi$ for all $z\in\Omega$, then $d_{g_n}$ converges to $d_g$ uniformly on every compact set of $\Omega$.
\end{thm}
This result plays an important role in Reshetnyak's proof of the local existence of generalized isothermal coordinates on Alexandrov surface with bounded curvature. Huber \cite{Huber} showed that the expression $ds^2=e^{2u(z)}|dz|^2$ (called ``line element'' in \cite{R}), where $u$ is representable as the difference of two subharmonic functions,  is invariant when moving from one isothermal chart $z$ to anther; consequently,  an orientable Alexandrov surface with bounded integral curvature is isometric to a Riemann surface equipped with a distance $d(x,y)=\inf_\gamma\int_\gamma e^{u(z)}|dz|$; the converse is also true:  On a (connected) Riemann surface there is a unique metric $d(x,y)$ for a conformally invariant line element making the surface an Alexandrov surface with bounded integral curvature.% (cf. \cite[p.102-p.103]{R}).

\vspace{.1cm}

Troyanov  \cite[Problem 9.1]{Troyanov1} asked the important global question: Is Reshetnyak's convergence theorem valid on a closed surface? In this paper, we give an affirmative answer to this open problem.

%n this paper, we study a basic question on uniform convergence of metrics on Alexandrov surfaces with bounded integral curvature. 

\vspace{.1cm}

  %{\color{blue}The signed measure $\mu$ is uniquely (Jordan) decomposed into a difference of two nonnegative Radon measures $\mu=\mu^+-\mu^-$. }
% Using the Green's function of $\Delta_{g_0}$, $G:\Sigma\times\Sigma\to\R\cup\{+\infty\}$ of $(\Sigma,g_0)$ 

%A metric space $(S,d)$ is an Alexandrov surface of bounded curvature  if $S$ is a closed topological surface and it satisfies two axioms: (A1) $d$ is intrinsic and (A2) any point has a neighborhood $G$ such that the total (upper) excess of any sequence of mutually disjoint simple triangles (i.e. of disk type and no two points on the edge can be joined by an arc outside of the triangle with shorter length) in $G$ is bounded above by a constant $C(G)$. 

%In the rich and long history of the development of Alexandrov surfaces of bounded curvature, there are notably three definitions \cite{R} while their equivalence is provided by deep theorems in the subject.  The axiomatic definition turns out to be equivalent to: 

A topological surface $S$ has a metric $d:S\times S\to\R$ with {\it bounded integral curvature} in the sense of Alexandrov (\cite{Troyanov1}, cf. \cite[6.1]{R}) if $d$ is continuous and 
\begin{enumerate}
\item[(i)] $d$ induces the manifold topology of $S$,
\item[(ii)] $d$ is intrinsic, i.e. for any $x,y\in S$ there exist curves $\gamma_n\in C^0([0,1],S)$ with $\gamma_n(0)=x,\gamma_n(1)=y$ so
 that the $d$-length of $\gamma_n$ converges to $d(x,y)$,
\item[(iii)]  there exist Riemannian metrics $g_n$ on $S$ with $\int_S|K_n|d\mu_{g_n}<C$, where $K_n$ is the Gauss curvature of 
 $g_n$ and $C$ is a constant, such that $d_{g_n}$ uniformly converge to $d$. 
\end{enumerate}

The above definition %of surfaces with bounded integral curvature in the sense of Alexandrov 
is equivalent to the original one given in \cite{AZ-1962} (see \cite[Remark (4)]{Troyanov1}). For a compact surface $S$ with an intrinsic metric $g$ of curvature bounded below by a constant $K\geq -1$ in the sense of Alexandrov (the triangle comparison) it is shown in \cite{Richard} (see also \cite{Slu})  that $g$ can be approximated by smooth metrics $g_n$ on $S$ (in the sense $d_{g_n}\to d_g$) with $K_{g_n}>-1$, hence $(S,g)$ has bounded integral curvature in the sense of Alexandrov. 
%Moreover,  in view of Corollary \ref{Riemann surface}, there is a conformal background smooth metric $h$ of constant curvature such that $d_g = d_{e^{2u}h}$ for some $u\in W^{1,q}(S)$ for any $q\in[1,2)$.

%The main purpose of this article is to establish uniform convergence of distance functions defined by metrics with bounded integral curvature on a closed surface (orientable or nonorientable) with a fixed conformal class $[g_0]$ as background,  directly corresponds to Theorem \ref{Resh}.   

\vspace{.2cm}

Let $\Sigma$ be a smooth surface with a Riemannian metric $g_0$.  Suppose that $u\in L^1_{\rm loc}(\Sigma)$ is an integrable function on $\Sigma$ with a well-defined Laplacian $\Delta_{g_0}u$ as a signed Radon measure $\mu$ satisfying 
\begin{equation}\label{def}
\int_\Sigma\varphi\,d\mu(g_u) = \int_\Sigma\left(\varphi \,K(g_0)-u\,\Delta_{g_0}\varphi\right)dV_{g_0},\s \mbox{for any $\varphi\in C^\infty_0(\Sigma)$,}
\end{equation}
 and set $\K_{g_u}:= \mu(g_u)$. The elliptic regularity implies : $u\in W^{1,q}_{\rm loc}(\Sigma)\footnote{In light of Weyl's lemma (%\cite{Weyl},
cf.  \cite[Theorem 2.3.1]{Morrey}),  if $u\in L^1_{\rm loc}(\Sigma,g_0)$ 
then $u-I_{\mu}\in C^\infty(\Sigma)$  where $I_\mu(x)\in W^{1,q}(\Sigma)$ (see Proposition \ref{B-M}); so $u\in W^{1,q}_{\rm loc}(\Sigma,g_0)$ for $q\in[1,2)$.  }$ and $e^u\in W^{1,q}_{\rm loc}(\Sigma)$ (Corollary \ref{MTI}). Let $\mathcal{M}(\Sigma,g_0)$ denote the set of $g=e^{2u}g_0$ for $u\in L^1_{\rm loc}(\Sigma)$ with \eqref{def}. Define 
\begin{equation}\label{distance-definition}
d_{g,\Sigma}(x,y)=\inf\left\{\int_\gamma e^u ds_{g_0}: \mbox{$\gamma$ is a piecewise smooth curve from $x$ to $y$ in $\Sigma$}\right\}
\end{equation}
where $\gamma$ is parametrized by its arclength parameter in $g_0$. The trace embedding theorem for Sobolev functions (cf. \cite[Theorem 18.1]{Leoni}, applied on each of the finitely many smooth pieces of $\gamma$) ensures integrability of $e^u$ along $\gamma$.

\vspace{.1cm}

The main goal of this paper is to prove the global version of Reshetnyak's convergence theorem: 
\begin{thm}\label{main-ChenLi}
Let $(\Sigma,g_0)$ be a closed surface and $g_k=e^{2u_k}g_0\in\M(\Sigma,g_0)$.  Assume that 
$\K_{g_k}=\mu_k^1-\mu_k^2$, where $\mu_k^1,\mu_k^2$ are nonnegative Radon measures. Assume $\mu_k^1,\mu_k^2$ converge to $\mu^1,\mu^2$ as measures, respectively, $\diam (\Sigma,g_k)=1$ and $\mu^1(\{x\})<2\pi$ for any $x$ in $\Sigma$. Then $u_k$
converges weakly to a function $u$ in $W^{1,q}$ for any $1\leq q<2$, $\K_g=\mu$ and $d_{g_k,\Sigma}$ converges to
$d_{g,\Sigma}$ uniformly where $g=e^{2u}g_0$.
%Moreover, for any domain $\Omega\subset\Sigma$, 
%$
%Area(\Omega,g_k)\rightarrow Area(\Omega,g).
%$
\end{thm}

The assumption $\mu^1(\{x\})<2\pi$ cannot be dropped, see the example in Appendix.
% where we glue Hulin-Troyanov's metric with the spherical metric.

 %A uniform bound on $\diam(\Omega,g_n)$ for Theorem \ref{Resh} follows from the assumption (see Remark \ref{diameter}). 
 
Theorem \ref{main-ChenLi} does not follow from Theorem \ref{Resh} by patching coordinate charts, the reason is that moving from one  chart $z$ to another $w$ yields a sequence of harmonic functions $u_k(z)-u_k(w)$ which is hard to control.  We develop a blow-up analysis around points where the curvature measure concentrates and prove the convergence by ruling out the trivial bubbles.

 \vspace{.1cm}
 
We now describe the idea in the proof of Theorem \ref{main-ChenLi}.  As $g_k$ is conformal to $g_0$ (by a possibly nonsmooth factor $e^{2u_k}$), it is convenient to use isothermal coordinates. Let $D$ be the unit disk in $\mathbb R^2$.  
%Suppose $(\Sigma,g)$ has bounded integral curvature in the sense of Alexandrov.  The background is $g=e^{2u}g_0$ for a Riemannian metric $g_0$ and $u\in L^1_{loc}(\Sigma,g_0)$; and $g_k=e^{2u_k}g_0$ when a sequence is relevant.  We highlight some key ingredients in our arguments: 

\vspace{.15cm}

1) When $|\K_g|(D)$ is small,  distance functions $d_{g_k,\Sigma}$ begin to converge (Proposition \ref{conv0}). This is achieved by proving that the distance function is comparable to the euclidean distance (Theorem \ref{estimate.d_g}) in this case. $\K_g(\{x\})<2\pi$ is crucial in establishing the so-called 3-circle type integral estimates which mimics the Fourier expansion of harmonic functions along a cylinder.

\vspace{.1cm}

2) When $\left|\K\right|_g(D)$ is finite, distance functions converge uniformly on compact sets away from (finitely many) curvature concentration points by 1) above. However,  the smallness required in 1) may not hold even on smaller disks, due to curvature concentration. We control the diameters in a scaling procedure and the key observation is that when $\K_g(\{x\})<2\pi$ only trivial bubbles develop (Proposition \ref{ghost}). Thus the blow-up analysis suffices for showing distance convergence.

\vspace{.1cm}

Theorem \ref{main-ChenLi} and the results established toward its proof enable us to approximate any metric in $\M(\Sigma,g_0)$ by Riemannian metrics with bounded total curvature. Consequently, any metric in $\M(\Sigma,g_0)$ has bounded integral curvature in the sense of Alexandrov, see Theorem \ref{2pi-main}. This important fact is known to Reshetnyak and Huber  (cf. \cite{Troyanov1}). 

\begin{thm}\label{2pi-main}
Let $(\Sigma,g_0)$ be a closed surface and $g=e^{2u}g_0\in \M(\Sigma,g_0)$ with $|\K_g|(\Sigma)<+\infty$. Assume $d_{g,\Sigma}$ is finite in $\Sigma\times\Sigma$. Then there exists smooth metric $g_k=e^{2u_k}g_0$, such that 
 \begin{enumerate} 
\item $u_k\to u$ in $W^{1,q}(\Sigma,g_0)$ for any $q\in[1,2)$, 
\item $\K_{g_k}\to \K_{g}$ in the sense of distributions,
\item $d_{g_k,\Sigma}\rightarrow d_{g,\Sigma}$ uniformly.   In particular, $g$ is a metric of bounded integral curvature in the sense of Alexandrov.
\end{enumerate}  
Moreover,  for any $R$ and $x\in\Sigma$, we have
\begin{equation}\label{vol.com}
\frac{\vol(B_R^g(x))}{\pi R^2}\leq 1+\frac{1}{2\pi}\K_g^-(\Sigma).
\end{equation}
\end{thm}

The essential idea in the proof of Theorem \ref{2pi-main} can be summarized as follows. When $\K_g(\{x\})<2\pi$ for any $x\in\Sigma$,  we first find a sequence $u_k\rightharpoonup u$ in $W^{1,q}$ with $\Delta u_k$ bounded in $L^1$ (Proposition \ref{app.glo}).  Then we scale $g_k=e^{2u_k}g_0$ to have fixed diameter, so we can apply Theorem \ref{main-ChenLi} to the normalized metrics.  When $\K_g(\{p_0\})\geq 2\pi$, as $d_{g,\Sigma}$ is finite it reduces to $\K_g(\{p_0\})=2\pi$ (Theorem \ref{2 option}).  On $D$, we can express $\mbox{$v=$ a harmonic function + the Poisson kernel for the signed measure } \K_g$ for $g=e^{2v}g_{\rm euc}$. Then using a cut-off function to mollify the Green's function $\log|x|$, we can construct metrics $g_k$ out of a sequence of functions approximating $v$ in $W^{1,q}$ such that $|\K_{g_k}|(D)\to|\K_g|(D)$ and $\K_{g_k}(\{p_0\})<2\pi$.  This discussion is not valid if $\K_g>2\pi$ somewhere; therefore we only treat the finite distance situation.

%This paper is self-contained. We feel that it will be convenient for the reader to have all details which are utilized to carry through our arguments in a unified single presentation. Our main source of reference on Alexandrov surfaces with bounded curvature is Reshetnyak's \cite{R}.%an excellent survey written by a major contributor to the field. 

%For the reader's convenience, we list the places where some of the key constants first appear in the paper: \noindent $\tau_0=\tau_0(1)$: Corollary \ref{MTI}; $\tau_1$: Lemma \ref{Sobolev.distance}; $\tau_2(\epsilon)$: Theorem \ref{estimate.d_g}; $\tau_3$: Corollary \ref{defi.dis2}; $\tau_0$, $\tau_0'$: Lemma \ref{3-circle}, \ref{3-circle'}.
		
%\begin{itemize}
%		\item $\tau_0=\tau_0(1)$: $e^u\in W^{1,1}$ on $\Sigma\setminus A_{\tau_0}$. cf. Corollary \ref{MTI}.
	%	\item $\tau_1$: $d_g\in W^{1,3}(\Sigma\backslash A_{\tau_1})$, cf. Lemma \ref{Sobolev.distance}
		%\item $\tau_2(\epsilon)$: cf. Theorem \ref{estimate.d_g}.
		%\item $\tau_3$: cf. Corollary \ref{defi.dis2}.
		%\item $\tau_0$, $\tau_0'$: for 3-circle lemma. cf. Lemma \ref{3-circle} and \ref{3-circle'}.
%\end{itemize}

For varying background conformal classes, we have 
\begin{thm}\label{main3}
Let $\Sigma$ be a closed surface of genus $\geq 1$. 
Assume that $h_k$ and $h_0$ are smooth metrics on $\Sigma$ with  $h_k\rightarrow h_0$ in the $C^2$-topology and $K_{h_k}=-1$ or $0$.
Let  $g_k=e^{2u_k}h_k\in\M(\Sigma,h_k)$.  Assume that $\K_{g_k}$ converges to a signed Radon measure $\mu$,   $\K_{g_k}^+$ converges to a Radon measure $\mu'$, and one of the following holds:
\begin{enumerate}
\item $\diam(\Sigma,g_k)=1$ and $\mu'(\{x\})<2\pi$ for any $x$ in $\Sigma$;
\item $d_{g_k,\Sigma}$ converges to a continuous distance function $d$ on $\Sigma$.
\end{enumerate}
Then, after passing to a subsequence, $u_k$
converges weakly to a function $u$ in $W^{1,q}$ for any $1\leq q<2$, $\K_g=\mu$ and $d_{g_k,\Sigma}$ converges to
$d_{g,\Sigma}$ uniformly, where $g=e^{2u}g_0$.
\end{thm}
\begin{comment}
{\color{red}This paper is self-contained. We feel that it will be convenient for the reader to have all details which are utilized to carry through our arguments in a unified single presentation. Our main source of reference on Alexandrov surfaces with bounded curvature is \cite{R}, \cite{Troyanov1}.}
\end{comment}

Global convergence has been addressed in \cite{Debin} under the assumption that the curvature measures at each point is less than $2\pi-\delta$ and the contractibility radius has a positive lower bound while conformal structures may vary. The contractibility radius is introduced in \cite{Debin}, replacing the role of injectivity radius of a smooth metric, to measure the longest loop which bounds a disk centered at a point. 
 %A positive lower bound prevents area collapsing in a limiting procedure. In fact, if imposing this condition then we can show that the conformal classes will stay in a bounded region hence subconverge {\color{blue}and the curvature measure will less than $2\pi$ at all points, why?}, see Corollary \ref{Debin-improve} and the last section of this paper.  
%Gromov-Hausdorff convergence of two-dimensional Alexandrov spaces with bounded integral curvature was studied 
On the other hand, it is shown in \cite{Shioya} that the space of Riemannian metrics on a closed surface with uniformly bounded total absolute curvature and diameter is precompact in the Gromov-Hausdorff distance and the limiting space may not be a topological surface.

It is tempting to show that the first nonzero eigenvalues of the Laplacians of the smooth metrics $e^{2u_k}g_0$ in Theorem \ref{2pi-main} converge to that of $e^{2u}g_0$. It is also curious to understand the completion of the moduli of conformal classes on a closed surface, with respect to the uniform distance convergence.

\vspace{.1cm}

{\bf Acknowledgements} We are grateful to Professor Zhichao Wang for his careful reading of the manuscript and useful suggestions. We would like to thank Professor Marc Troyanov for his interest in this work. We are grateful for the referees for their valuable suggestions. 

\section{Preliminaries on signed Radon measure $\Delta u$}

\subsection{Gauss curvature measures}

%Before state our results,  we introduce terminologies and notations that we will use. 

Let $\Sigma$ be a smooth surface without boundary (not necessarily compact) with a Riemannian metric $g_0$, the Gauss curvature $K(g_0)$ and the area element $dV_{g_0}$.  
For any $g_u=e^{2u}g_0\in\M(\Sigma,g_0)$ %on $\Sigma$. 
%Let $\M(\Sigma,g_0)$ be the set of $g_u$ so that there is a {\it signed} Radon measure $\mu(g_u)$ satisfying 
%\begin{equation}\label{def}
%\int_\Sigma\varphi \,d\mu(g_u)=\int_\Sigma \left(\varphi \,K(g_0) - u \Delta_{g_0}\varphi\right)dV_{g_0},\s \forall \varphi\in C^\infty_0(\Sigma).
%\end{equation}
  we introduce notations 
\begin{equation}\label{notation}
dV_{g_u}:= e^{2u} dV_{g_0} \ \ \  \mbox{and} \ \ \ 
\K_{g_u}:= \mu(g_u).
\end{equation}
%In these notations \eqref{def} reads
%\begin{equation*}
%\int_\Sigma\varphi \,\K_{g_u}=\int_\Sigma \left(\varphi K(g_0)+\nabla_{g_0} u \nabla_{g_0}\varphi\right)dV_{g_0}. 
%\end{equation*}
We call the signed Radon measure $\K_{g_u}$ the {\it Gauss curvature measure} for the measurable tensor $g_u$. For a positive constant $c$ it holds $\K_{cg_u}=\K_{g_u}$. When $u$ is smooth, 
$$
K(g_u) = e^{-2u}\left(K(g_0)-\Delta_{g_0}u\right).
$$

In an isothermal coordinate chart $(x,y)$ for $g_0$, we can write $g_0= e^{2u_0}g_{\rm euc}$ for a locally defined smooth function $u_0$. Any $g\in\M(\Sigma,g_0)$ is locally as $g=e^{2v}g_{\rm euc}$, where $v\in L^{1}_{\rm loc}(\Sigma)$
and
\begin{equation}
-\Delta v \ dxdy =\K_{g_v}
\end{equation}
as distributions and $\Delta=\frac{\partial^2}{\partial x^2}+\frac{\partial^2}{\partial y^2}$.

\begin{rem}
	The above definition of Gauss curvature measure is global and holds for nonorientable surfaces as well since the divergence theorem is valid there. Locally it 
	is essentially the one given by Reshetnyak via his subharmonic metrics in a domain of $\C$ with $g_0=|dz|^2$ where $u$ admits the representation in Theorem \ref{Resh} (the logarithmic potential of $g$). Huber showed that one can patch the subharmonic metrics from one isothermal chart to another on an orientable surface.
%Definition \ref{def} is global and holds for nonorientable surfaces as the divergence theorem holds for nonorientable surface. 
	\end{rem}

%\subsection{A Gauss-Bonnet formula}
%When $\Delta u$ is a signed Radon measure, t

We now state a Gauss-Bonnet formula in the nonsmooth settiing. 
Let $u\in W^{1,1}(D_{R_2}\backslash D_{R_1})$. From the trace embedding theorem for Sobolev functions
%, \cite[Section 5.6]{Attouch-Buttazzo-Michaille}), 
we know that $u(re^{i\theta})$ is defined in $L^1(\mathbb S^1)$ for any $r\in[R_1,R_2]$. Then 
$$
u^*(r)=\frac{1}{2\pi}\int_0^{2\pi}u(re^{i\theta})\,d\theta
$$
is well-defined  on $[R_1,R_2]$. %When $-\Delta u=\K_g$ measures the change of the radial derivative of $u^*$ for a.e. $r$. 
The following Gauss-Bonnet formula is essentially a Green's formula and we will present a proof in Appendix.

\begin{lem}\label{Gauss-Bonnet-Green}
Then for almost every $s,t\in(R_1,R_2)$,  $s<t$,  we have 
\begin{equation}\label{GB1}
t\frac{du^*}{dr}(t)- s\frac{du^*}{dr}(s)=-\frac{1}{2\pi}\,\K_g(D_t\backslash {D}_s).
\end{equation}
There exists $E\subset[0,R]$ with ${\mathscr L}^1(E)=0$, such that for any $\{t_k\}\subset[0,R]\backslash E$ with $t_k\rightarrow 0$ it holds
\begin{equation}\label{GB2}
\K_g(\{0\}) =-2\pi\lim_{t_k\to 0}t_k\frac{du^*}{dr}(t_k).
\end{equation}
\end{lem}

When $u$ is smooth, \eqref{GB1} is the classical Gauss-Bonnet on the annulus (both sides of \eqref{GB1} equal 0) and the limiting case \eqref{GB2} for a disk captures $\K_g$ at 0 as a measure ($K_g dxdy$ at 0 is zero for smooth $u$).

\subsection{Uniform estimates}

Let $\mu$ be a signed Radon measure on a domain $\Omega$. It is well-known that there is
a Radon measure $|\mu|$ and a $|\mu|$-measurable function $\nu$ with $|\nu|=1$, such that 
\begin{equation*}\label{signed Radon}
	\int_\Omega\varphi d\mu=\int_{\Omega}\varphi\nu d|\mu|,\s \forall \varphi
	\in C_0(\Omega).
\end{equation*}
Setting $\mu^{\pm}=|\mu| \llcorner\nu^\pm$ then 
$\mu=\mu^+-\mu^-$ and $ |\mu|=\mu^++\mu^-$.  This decomposition is unique according to the Jordan decomposition theorem. It implies that at each $x\in\Omega$ at least one of $\mu^+(\{x\})$ and $\mu^-(\{x\})$ is 0. In particular, we will write
\begin{equation}\label{signed Radon}
\K_g = \K^+_g -\K^-_{g} \s \mbox{and}\s |\K_g| = \K^+_g+\K^-_g.
\end{equation}
A weak solution on $\Omega\subset\mathbb R^2$ to 
\begin{equation}\label{the equation}
\Delta u = -\mu
\end{equation}
means 
$$
\int_\Omega u\Delta\varphi \,dxdy = -\int_\Omega \varphi v d|\mu|, \s\s \forall\varphi\in C^\infty_0(\Omega). 
$$
%where $v$ is the $|\mu|$-measurable function  in \eqref{signed Radon}.

\vspace{.1cm}

With slight modification, the proof in \cite{B-M} can be adapted to the current setting. 

\begin{pro}[Brezis-Merle]\label{B-M}  
Given a signed Radon measure $\mu$ supported in $D\subset \R^2$ with $0<|\mu|(\R^2)<+\infty$, let 
$$
I_\mu(x)=-\frac{1}{2\pi}\int_{\R^2}\log|x-y|d\mu(y).
$$
Then $I_\mu \in W^{1,q}_{\rm loc}(\R^2)$ for any $q\in[1,2)$ and weakly solves the 
 equation: 
\begin{equation}\label{equation.I}
-\Delta I_\mu=\mu.
\end{equation}
Moreover, we have
\begin{equation}\label{I}
\int_{D_R}\left|I_\mu(x)\right|^qdx\leq C(q,R)|\mu|(\mathbb R^2)^q,
\end{equation}
\begin{equation}\label{nabla.q}
r^{q-2}\int_{D_r(x)}|\nabla I_\mu|^qdx\leq C(q)|\mu|(\R^2)^q,\s \forall x, r,
\end{equation}
and
\begin{equation}\label{whynot}
\int_{D_{R}} e^{\frac{(4\pi-\epsilon)|I_\mu|}{|\mu|(\R^2)|}}dx\leq CR^\frac{\epsilon}{2\pi},\s \forall R>0\s and\s \epsilon\in(0,4\pi)
\end{equation}
where $C$ is a constant independent of $\epsilon, R, \mu$. 
\end{pro}

\begin{rem}\label{diameter}
In Theorem \ref{Resh}, the assumption implies that $|\mu_n|(\R^2)$ is uniformly bounded.  In light of Proposition \ref{B-M} and compactness of $\overline\Omega$,  it follows that $u_n$ is uniformly bounded in $W^{1,q}(\R^2)$.  The trace embedding theorem for Sobolev functions then 
asserts uniform boundedness of $\diam(\Omega,g_k)$.  
\end{rem}

 For a positive measure $\mu$ supported in $D$ with $\mu(D)<\frac{2\pi}{q}$, item 1) below is observed in 
 \cite[Corollary 4.3]{Troyanov3} (cf. \cite[TH. 3.1]{R2}).
 % that $\|e^{2u}\|_{L^q(D)}$ depends only on $q$ and $\mu(D)$ for the logarithmic potential 
 %$u=-\frac{1}{2\pi}\int_{\C}\log|z-\zeta|d\mu(\zeta)$. This is essentially 1) below.
 % the second conclusion there may not be found in \cite{Troyanov3} as explicit as $\|e^u\|_{W^{1,q}(D_{1/2})}\leq\beta'$. 

\begin{cor}\label{MTI} 
Let $\mu$ be a signed Radon measure on $D$ with $|\mu|(D)<\tau$. Suppose that  $u$ 
solving \eqref{the equation} weakly and $\|u\|_{L^1(D)}<\gamma$.    
Then
\begin{enumerate}
\item[1)] $u\in W^{1,q}(D_{1/2})$ for any $q\in[1,2)$. Moreover, 
$$
\|\nabla u\|_{L^q(D_{1/2})}<C(q)
(\|u\|_{L^1(D)}+|\mu|(D)).
$$
\item[2)] for any $p<\frac{4\pi}{\tau}$ there exists  $\beta=\beta(\tau,p,\gamma)$ such that
$$
\int_{D_{1/2}}e^{p|u|}dx\leq \beta.
$$
\end{enumerate}
Moreover, for any $1\leq q<2$ and $\frac{2q}{2-q}<\frac{4\pi}{\tau}$,
$$
\|e^u\|_{W^{1,q}(D_{1/2})}\leq \beta',
$$
where $\beta'$ only depends on $q$, $\tau$ and $\gamma$.
\end{cor}
\proof Extend $\mu$ to a measure on $\R^2$ by $\mu(A) = \mu(A\cap D)$ for any $A\subseteq \R^2$, and denote the extension by $\mu$ for simplicity. If $\mu=0$ then $u$ is a smooth harmonic function as it is a distributional solution \cite[Theorem 2.3.1]{Morrey},  so the corollary holds, see argument below for $u^{\rm har}$.  Now assume $\mu$ is not the zero measure.
According to Proposition \ref{B-M}, $v:=I_\mu$ solves \eqref{the equation} in $\R^2$ weakly. 
Let $u^{\rm har}=u-v$. Then $\int_D u^{\rm har}\Delta\varphi=0$
for any $\varphi\in C_0^\infty(D)$. By Weyl's Lemma,  $u^{\rm har}$ is a smooth harmonic function on $D$.  
Then by \eqref{I} with $q=1, R=1$ there for $v$  
$$
\|u^{\rm har}\|_{L^1(D)}\leq\|u\|_{L^1(D)}+\|v\|_{L^1(D)}< \|u\|_{L^1(D)}+C|\mu|(D)\leq \gamma+ C\tau.
$$ 
By the mean value theorem for harmonic functions,
\begin{equation}\label{C1}
\|u^{\rm har}\|_{C^0(D_{3/4})}\leq\sup_{D_{3/4}}\frac{1}{D_{1/4}(x)}\left|\int_{D_{1/4}(x)}u^{\rm har}dx\right|\leq
C(\|u\|_{L^1(D)}+|\mu|(D)).
\end{equation}
Then, it follows from the elliptic estimates that
$$
\|u^{\rm har}\|_{C^1(D_{1/2})}<C(\|u\|_{L^1(D)}+|\mu|(D)).
$$
Together with \eqref{nabla.q}, we complete the proof of 1).

%{\color{red}
%By \eqref{C1}
%\begin{equation}\label{mean}
%\|u^{har}-\bar{u}^{har}\|_{L^\infty(D_{1/2})}<C
%\end{equation}
%where $\bar{u}^{har}$ is the mean value of $u^{har}$ over $D$. }

Taking $\epsilon =4\pi\left(1-\frac{|\mu|(D)}{\tau}\right)>0$ we have 
$p<\frac{4\pi}{\tau}=\frac{4\pi-\epsilon}{|\mu|(\R^2)}$. By Proposition \ref{B-M}, 
%{\color{red}	$$
%\|v\|_{L^1(D)}<C\tau\s \mbox{and} \s\int_De^{p|v |}<C.$$}

	\begin{align}\label{useonce}
		\int_{D_{1/2}}e^{p|u|}= \int_{D_{1/2}}e^{p|u^{\rm har}+v|}
		\leq e^{pC(\|u\|_{L^1(D)}+|\mu|(D))}\int_{D_{1/2}}e^{p|v|}\leq C(\tau,p,\gamma).
	\end{align}
	
%{\color{red}
%\begin{align*}\label{useonce}
%\int_{D_{{1/2}}}e^{p|u|}&= \int_{D_{1/2}}e^{p|u^{har}+v|}\nonumber\\
%&\leq\int_{D_{1/2}}e^{p|u^{har}-\bar{u}^{har}|+p|v-\bar{v}|+p|\bar{u}|} \s\s\s \mbox{(as $\bar{u} = \bar{u}^{har}+\bar{v}$)}\nonumber\\
%&\leq C\int_{D_{1/2}}e^{p|v|+p|\bar{v}|+p|\bar{u}|} \s\s\s\s\s\s\s\s \mbox{(by \eqref{mean})}\\
%&\leq Ce^{p|\bar{v}|+p|\bar{u}|}\int_{D_{1/2}}e^{p|v|} \nonumber\\
%&\leq C(\tau,p,\gamma). \nonumber
%\end{align*}
%no need to average}

When  $2q/(2-q)<{4\pi}/{\tau}$, we have $q<{4\pi}/{\tau}.$ Let $1/q'+{\tau}/{4\pi}=1$. Also by $2q/(2-q)<{4\pi}/{\tau}$, we have $1/qq'=(1-{\tau}/{4\pi})/q>1/q^2+1/(2q)$. As $1\leq q<2$, we obtain $qq'<2$. Using (generalized) H\"older's inequality and \eqref{useonce}, 
\begin{align*}
\|\nabla e^{u}\|_{L^q(D_{1/2})}&\leq\|\nabla u\|_{L^{q'q}(D_{1/2})}\|e^{u}\|_{L^\frac{4\pi}{\tau}(D_{1/2})} \leq C\|\nabla u\|_{L^{qq'}(D_{1/2})}\|e^{u}\|_{L^\frac{4\pi}{\tau}(D_{1/2})} \\
&\leq C \left(\|\nabla u^{\rm har}\|_{L^{qq'}(D_{1/2})}+\|\nabla v\|_{L^{qq'}(D_{1/2})}\right).%{\color{red}e^{\bar{u}}?}.
\end{align*}
By \eqref{nabla.q} we see $\|\nabla v\|_{L^{qq'}(D_{{1/2}})}<C$. Estimates for harmonic functions lead to 
%Gradient estimates for harmonic functions (cf. \cite[Theorem 2.10]{GT}) and the mean value theorem assert 
$$
\|\nabla u^{\rm har}\|_{L^{qq'}(D_{1/2})}\leq \sup_{D_{1/2}} |\nabla u^{\rm har}|\leq C\sup_{D_{3/4}}|u^{\rm har}| \leq C \|u^{\rm har}\|_{L^1(D)}<C. 
$$
Now the desired result follows. \endproof

%\subsection{$W^{1,q}$ estimate} 

\begin{lem}\label{rq}
Let  $u\in W^{1,1}(D)$
solve \eqref{the equation}. If $\|\nabla u\|_{L^1(D)}<A$, then for $q\in [1,2), r<1$ 
\begin{equation}\label{II}
\|\nabla u\|_{L^q(D_{r})}\leq CAr^{\frac{2}{q}}+C(q)r^{\frac{2-q}{q}}|\mu|(D).
\end{equation}
\end{lem}
\proof
%As $\mathcal{D}$ is dense in $W^{1,q}_0\cap C^0(D)$, for any $\varphi\in W^{1,0}_0\cap C^0(D)$, 
%\begin{equation}\label{weak.equation}
%\int_D\nabla u\nabla\varphi dx=\int_D\varphi d\mu.
%\end{equation}
Extend $\mu$ to a signed Radon measure on $\R^2$
by setting $\mu(D^c)=0$ and write
$u=u^{\rm har}+I_\mu$, 
where $u^{\rm har}\in C^\infty(D)$ is harmonic. 
%Then, by Proposition \ref{B-M},  $u\in \cap_{q\in [1,2)}W^{1,q}_{loc}(D).$
Using properties of harmonic functions as above %in the proof of Corollary \ref{MTI} %(we only demonstrate the case $D_r\subset D_{1/2}$, namely centered at $0$)
\begin{align*}\label{II}
\|\nabla u\|_{L^q(D_{r})} &\leq \|\nabla u^{\rm har}\|_{L^q(D_r)}+ \|\nabla I_\mu\|_{L^q(D_r)}\nonumber\leq Cr^{\frac{2}{q}}\sup_{D_{1/2}} |\nabla u^{\rm har}| + C(q) r^{\frac{2-q}{q}}|\mu|(D)  \\&\leq Cr^{\frac{2}{q}} A + C(q) r^{\frac{2-q}{q}}|\mu|(D).
\end{align*}
\endproof

\begin{rem}
	When $u$ is smooth, under an area growth condition for $g=e^{2u}g_{\rm euc}$, a priori estimates for $\|\nabla u\|_{L^q},q\in(1,2)$ were obtained in \cite[Theorem 1.3]{LST}.

%\eqref{nabla.q} is not general true without the assumption $\|\nabla u\|_{L^1(D)}<A$. For example, $u_k=kx^1-c_k$ is harmonic, but $\|\nabla u_k\|_{L^q(D_r(x))}\rightarrow+\infty$ for any $D_r(x)\subset D$.  In \cite{LST}, it was proved that $\|nabla u\|_{L^1(D)}<A$ is equivalent to the euclidean volume growth when $u$ is smooth. 
%In a forthcoming paper, we will show such a result is also valid when $u\in W^{1,1}$.
\end{rem}

The following global gradient estimate is known (cf. \cite[Proposition 2.19]{Troyanov1}). It will be used both for a fixed conformal structure and for varying conformal structures on a closed surface in the proof of Theorem \ref{main-ChenLi} and Theorem \ref{main3}, respectively. 
%The result below will be used when we consider varying conformal structures on a closed orientable surface and the conformal structures stay in a bounded region in the moduli space of Riemann surfaces with fixed genus. 

\begin{lem}\label{global.gradient.estimate}
Let $\mu$ be a signed Radon measure defined on a closed Riemannian surface $(\Sigma,g)$ and $u\in L^1(\Sigma)$ solves 
$
-\Delta_g u=\mu.
$
We assume $\|g-g_0\|_{C^{2,\alpha}}<a$.
Then, for any $r>0$ and $q\in [1,2)$ there exists $C=C(q)$ such that
$$
r^{q-2}\int_{B_r(x)}|\nabla_{g}u|^q\leq C \, |\mu|(\Sigma)^q.
$$
where $B_r(x)$ is the geodesic ball in $g$. 
\end{lem}

\subsection{Approximation}

\begin{lem}\label{local}
Let $u\in L^1(D)$ with $-\Delta u=\mu$, where $\mu$ is a signed Radon measure compactly supported in $D$. Then there exist $u_k,f_k\in C^\infty_0(\R^2)$ with $-\Delta u_k=f_k$ and 
\begin{enumerate}
\item[{(1)}] $u_k$ converges to $u$ in $W^{1,q}_{\rm loc}(D)$ for any $q\in[1,2)$; 
\item[{(2)}]there are smooth functions $f_k^1,f_k^2\geq 0$ so that $f_k=f_k^1-f_k^2$ with $\|f_k\|_{L^1}\leq |\mu|(\R^2)+\frac{1}{k}$ and $f_k^1dx\rightharpoonup\mu^+$,  $f_k^2dx\rightharpoonup\mu^-$, 
$f_kdx\rightharpoonup \mu$ as measures. 
\end{enumerate}
\end{lem}

\proof Let $0\leq \eta\in C^\infty_0(D)$ with 
$\int_{\R^2}\eta=1$ and $\eta_k(x)=\eta(\frac{x}{\epsilon_k})/\epsilon^2_k$ where $\epsilon_k\rightarrow 0$. Denote 
$$
u_k(x)=\int_{\R^2}u(y)\eta_k(x-y)\,dy
$$
and
$$ 
f_k^1(x)=\int_{\R^2}\eta_k(x-y)\,d\mu^+(y),\s f_k^2(x)=\int_{\R^2}\eta_k(x-y)\,d\mu^-(y).
$$
By the dominated convergence theorem, $u_k,f_k^1, f_k^2\in C^\infty_0(\R^2)$. By Corollary \ref{MTI} 1), $u\in W^{1,q}_{\rm loc}(D)$, so $u_k\to u$ in $W^{1,q}_{\rm loc}(D)$ (cf. \cite[Theorem 4.1]{E-G}). Further, $f_k^1 dx\rightharpoonup\mu^+$ and $f_k^2 dx\rightharpoonup\mu^-$:  $\forall\phi\in C^\infty_0(\R^2)$, it holds 
 \begin{align*}
 \lim_{k\to\infty}\int_{\R^2} \phi(x) \int_{\R^2} \eta_k(x-y)d\mu^\pm(y) dx &= \lim_{k\to\infty}\int_{\R^2}\int_{\R^2} \phi(x)\eta_k(x-y)dx d\mu^\pm(y)\\
 &=\int_{\R^2}\phi(y)d\mu^\pm(y).
\end{align*}
 Let $f_k=f_k^1-f_k^2$ and $\mbox{supp}\,(f_k)\subset\subset D_2$ for large $k$.
By \cite[Theorem 1.40 (iii)]{E-G}, 
$$
\int_{D_2}\left(f_k^1+f_k^2\right)dx\rightarrow |\mu|(D_2).
$$
Without loss of generality, we assume $\|f_k\|_{L^1(\R^2)}\leq|\mu|(\R^2)+\frac{1}{k}$.  Moreover, 

\begin{align*}
\int_{\R^2}\nabla u_k\nabla \varphi \,dx&=\int_{\R^2}\int_{\R^2}u(y)\nabla_x\eta_k(x-y)\nabla_x \varphi(x) \,dydx\\
&=\int_{\R^2}\left(\int_{\R^2}\nabla_y u(y) \eta_k(x-y) \right)dy\nabla_x\varphi(x) \,dx\\
%&=\int_{\R^2}\left(\int_{\R^2}- \varphi(x)\nabla_x\eta_k(x-y) dx\right)\nabla_y u(y)\,dy =\int_{\R^2}\left(\int_{\R^2}\varphi(x)\nabla_y\eta_k(x-y)  dx\right)\nabla_y u(y)\,dy \\
&=\int_{\R^2}\nabla_y\left(\int_{\R^2}\eta_k(x-y) \varphi(x) dx\right)\nabla_y u(y)\,dy \\
&=\int_{\R^2}\left(\int_{\R^2}\eta_k(x-y) \varphi(x) \,dx\right)d\mu(y)\\
&=\int_{\R^2}\left(\int_{\R^2}\eta_k(x-y) d\mu(y)\right)\varphi(x) \,dx=\int_{\R^2}f_k(x)\varphi(x) \,dx.
\end{align*}
\endproof

%A partition of unity leads to a global statement: 
\begin{pro}\label{app.glo}
Let $\Sigma$ be a surface with a Riemmanian metric $g$. Let $u\in L^{1}(\Sigma,g)$ such that $-\Delta u$ is a signed Radon measure $\mu$. Then
there exists $u_k\in C^\infty(\Sigma)$ so that 
\begin{enumerate}
\item $u_k$ converges to $u$ in $W^{1,q}$;
\item $\|\Delta u_k\|_{L^1(\Sigma, g)}<C$ and $-\Delta u_k$ converges to $\mu$ in the sense of distributions;
\item There are smooth functions $F_k^1,F_k^2\geq 0$, such that $-\Delta u_k=F_k^1-F_k^2$ and $F_k^1$, $F_k^2$ converge to
$\mu^+$, $\mu^-$ in the sense of distributions, respectively.
\end{enumerate}
\end{pro}
\proof Let $\{h_\alpha\}$ be a partition of unity subordinate to an open covering of $\Sigma$ by coordinate disks. Let $\mu_\alpha=\Delta(h_\alpha u)$ be the signed Radon measure.  So $\sum \mu_\alpha = \mu$. 
From Lemma \ref{local} for each $\alpha$, there is a sequence $u^\alpha_k\to h_\alpha u$ in $W^{1,1}$. Define $u_k = \sum h_\alpha u^\alpha_k$. As $k\to\infty$, $u_k$ tends to $\sum h_\alpha u = u$. 
\endproof

\begin{lem}\label{W1q.convergence}
Let $\mu_k$ be a signed Radon measure on $D$ and $u_k\in L^1(D)$
solve 
$-\Delta u_k=\mu_k$
weakly for each $k$. Assume that  $\mu_k$ converges to a Radon measure $\mu$ weakly and
$u_k\rightarrow u$ in $L^1(D)$. Then
\begin{enumerate}
\item[1)] $-\Delta u=\mu$ holds weakly and $u_k$ converges to $u$ weakly in $W^{1,q}_{\rm loc}(D)$, $\forall q\in[1,2)$;
\item[2)] If $|\mu_k|(D)\rightarrow 0$ %in addition,  
then $u_k\rightarrow u$
in $W^{1,q}_{\rm loc}(D),\forall q\in[1,2)$ and $u$ is a harmonic function on $D$.
\end{enumerate}
\end{lem}
\proof 
Given $\varphi\in C_0^\infty(D)$, we have
$$
-\int_D u_k\Delta\varphi=\int_D\varphi d\mu_k.
$$
Letting $k\to\infty$ 
$$
-\int_Du\Delta\varphi dx=\int_D\varphi d\mu.
$$ 
Hence $u$ solves $-\Delta u=\mu$ weakly.
Next, we let $\varphi\in C^\infty(\R^2)$. Then for any $r\in(0,1)$,
$$
\int_{D_r}\nabla u_k\nabla\varphi dx=-\int_{D_r}u_k\Delta\varphi+\int_{\partial D_r}u_k\frac{\partial\varphi}{\partial r}.
$$
We claim 
$$
\int_{\partial D_r}u_k\frac{\partial\varphi}{\partial r}\rightarrow\int_{\partial D_r}u\frac{\partial\varphi}{\partial r}.
$$
It suffices to prove the claim for any convergent subsequence. By Corollary \ref{MTI} 1) and that $u_k\to u$ in $L^1(D)$, by passing to a subsequence, we assume $u_k\rightharpoonup u$ in $W^{1,q}(D_r)$. Now the claim follows from 
the trace embedding theorem.  
Then 
$$
\int_{D_r}\nabla u_k\nabla\varphi dx\rightarrow-\int_{D_r}u\Delta\varphi+\int_{\partial D_r}u\frac{\partial\varphi}{\partial r}=\int_{D_r}\nabla u\nabla\varphi,
$$
so $u_k\rightharpoonup u$ in $W^{1,q}_{\rm loc}(D)$. As $|\mu_k|(D)\rightarrow 0$ by assumption, \eqref{nabla.q} and \eqref{I} imply 
$
\|I_{\mu_k}\|_{W^{1,q}(D)}\rightarrow 0.
$
 Then 
$$
\|u_k^{\rm har}- u\|_{L^1(D)}\leq \|I_{\mu_k}\|_{L^1(D)}+\|u_k-u\|_{L^1(D)}\rightarrow 0.
$$
Now, for any $\phi\in C^\infty_0(D)$
$$
\int_D u \, \Delta \phi  = \int_D (u-u^{\rm har}_k) \, \Delta \phi \rightarrow 0.
$$ 
Hence $u$ is harmonic on $D$. Then for any compact region $D'$ in $D$
$$
\|u_k-u\|_{W^{1,q}(D')} \leq \|I_{\mu_k}\|_{W^{1,q}(D')} + \| u^{\rm har}_k -u\|_{W^{1,q}(D')} \rightarrow 0,
$$
as $u_k^{\rm har}-u$ is harmonic and tends to 0 in $L^1$. 
\endproof

\section{Distance function}

\subsection{Distance of singular metric as Sobolev function}

Let $(\Sigma,g_0)$ be a Riemannian surface without boundary. For $g=e^{2u}g_0\in\M(\Sigma,g_0)$ recall 
$
\K_g = \left(K(g_0)-\Delta_{g_0} u \right) dV_{g_0}.
$
The $\K_g$-measure of a point may not be 0 for nonsmooth $u$. To investigate curvature concentration, set
\begin{equation}\label{defi.A}
A_\epsilon=\left\{x\in\Sigma:|\K_g|(\{x\})\geq\epsilon\right\}, \ \ \ \epsilon>0.
\end{equation}
As the Radon measure $|\K_g|$ is locally finite, the set $A_\epsilon$ is discrete with no accumulation points.

Cover $\Sigma\backslash A_\epsilon$ by open sets $U_\alpha\subset\Sigma\backslash A_\epsilon$ so that each $U_\alpha$ is conformal to $D$. On $U_\alpha$, $g_0=e^{2u_0}g_{\rm euc}$ for some smooth $u_0$. Then $g=e^{2(u-u_0)}g_{\rm euc}$,  $\K_g =-\Delta (u-u_0) dV_{g_{\rm euc}}$ and $|\K_g|(D)\leq 2\epsilon$
by choosing $U_\alpha$ small. 
Applying Corollary \ref{MTI} to $u-u_0$ yields $e^{u}\in W^{1,1}_{\rm loc}(\Sigma\backslash A_\epsilon)$ when $2\epsilon\leq \frac{4}{3}\pi$. 

Let $x,y\in\Sigma$ and $\gamma$ be a piecewise smooth curve from $x$ to $y$ in $\Sigma$.  For any $t\notin \gamma^{-1}(A_\epsilon)$, we can find an interval $(t-\delta,t+\delta)$ not intersecting $\gamma^{-1}(A_\epsilon)$. By the trace embedding theorem,
%(cf. \cite[Theorem 3.45]{HT}) %or \cite[section 1.3]{MS}
 $e^{u}$ is measurable on $\gamma(t-\delta,t+\delta)$. Since $\gamma^{-1}(A_\epsilon)$ is at most countable, $\ell_g(\gamma):=\int_\gamma e^{u}ds_{g_0}$ is well defined.
Define $d_{g,\Sigma}:\Sigma\times\Sigma\rightarrow [0,+\infty]$ by \eqref{distance-definition}. 
%$$
%d_{g,\Sigma}(x,y)=\inf\left\{\ell_{g}(\gamma)=\int_{\gamma}e^{u}ds_{g_0}:\gamma \,\,\hbox{is a piecewise smooth curve from $x$ to $y$ in $\Sigma$}\right\}.
%$$
Note that $d_{g,\Sigma}(x,y)<+\infty$
for any $x,y\in \Sigma\backslash A_\epsilon$. This is because we can take a piecewise smooth curve $\gamma$ from $x$ to $y$ in $\Sigma
\backslash A_\epsilon$ with $\int_{\gamma}e^u<+\infty$ by the trace embedding theorem.

\begin{lem}\label{Sigma12} Let %$(\Sigma,g_0)$ be a Riemannian surface and 
$g\in\M(\Sigma,g_0)$ and let $\Omega_1,\Omega_2$ be relatively compact domains in $\Sigma$ with piecewise smooth boundary. If $\overline{\Omega}_1\subset\Omega_2$ then 
\begin{enumerate}
\item[(i)]  $d_{g,\Omega_2}(x,\partial\Omega_1)\!=\!\inf\left\{\ell_{g}(c)\left|\right. c:[0,1]\rightarrow\Sigma, c(0)=x, c(1)\in\partial\Omega_1, c((0,1))\subset\Omega_1\right\}, \forall x \!\in\Omega_1;$
\item[(ii)] $d_{g,\Sigma}(\partial \Omega_1,\partial\Omega_2)=\inf\left\{\ell_{g}(c)\left|\right.c:[0,1]\rightarrow\Sigma, c(0)\in\partial\Omega_1, c(1)\in\partial\Omega_2, 
 c((0,1))\subset\Omega_2\backslash\overline{\Omega}_1\right\};$
\item[(iii)]  $d_{g,\Omega_2}(x,y)\leq d_{g,\Omega_1}(x,y),\forall x,y\in \Omega_1;$
\item[(iv)]  if $x,y\in \Omega_1$ and 
$d_{g,\Omega_2}(x,y)\neq d_{g,\Sigma}(x,y)$, then
$
d_{g,\Sigma}(x,y)\geq d_{g,\Sigma}(\partial \Omega_1,\partial\Omega_2).
$
\end{enumerate}
\end{lem}
\proof 

(i) Let $c$ be a curve in $\Omega_2$ from $x$ to a point $y\in\partial \Omega_1$. If $c$ leaves $\Omega_1$ and $x$ is interior in $\Omega_1$ then $c$ must hit $\partial\Omega_1$ first before departing $\Omega_1$ (if $c$ departed $\Omega_1$ from an interior point of $\Omega_1$ then $\Sigma$ would not be a manifold there). Any such $c$ does not affect the infimum in the definition of $d_{g,\Omega_1}$.

\vspace{.1cm}

(ii) Let  $c_k$ be a curve from $x_k^1\in\partial\Omega_1$ to $x_k^2\in\partial\Omega_2$, and $\lim_{k\to+\infty}\ell_g(c_k)=d_{g,\Sigma}(\partial\Omega_1,\partial\Omega_2)$. If $c_k((0,1))\subset\Omega_2\backslash\overline{\Omega}_1$, we let $c_k'=c_k$. Otherwise, we let $t_1$ be the greatest $t$ for $\gamma(t)\in\partial\Omega_1$ and $t_2$  the least $t$ for $\gamma(t)\in\partial\Omega_2$, and
define
$
c_k'(t)=c_k(t_1+(t_2-t_1)t).
$
Since $\ell_g(c_k')\leq \ell_g(c_k)$, we get (ii).

\vspace{.1cm}

(iii) The conclusion follows from the fact that any curve connecting $x,y$ in $\Omega_1$ is also in $\Omega_2$.

\vspace{.1cm}

(iv)
Let $c_k$ be a curve in $\Sigma$ from $x\in\Omega_1$ to  $y\in\Omega_1$ such that $\ell_g(c_k)\rightarrow d_{g,\Sigma}(x,y)$. If $c_k$ leaves $\Omega_2$ then $c_k$ must meet $\partial\Omega_1$ and $\partial\Omega_2$, so $\ell_g(c_k)\geq d_g(\partial\Omega_1,\partial\Omega_2)$, in turn, $d_{g,\Sigma}(x,y)\geq d_{g,\Sigma}(\partial\Omega_1,\partial\Omega_2)$. Otherwise, $c_k\subset\Omega_2$, then $\ell_g(c_k)\rightarrow d_{g,\Omega_2}(x,y)$, but this contradicts $d_{g,\Omega_2}(x,y)>d_{g,\Sigma}(x,y)$. 
 \endproof

\vspace{.1cm}

%\begin{rem}
It follows from (i) and (ii): for any $x\in\Omega_1\subset\Omega_2\subset\Omega_3\subset\Omega_4$, there holds
$$
d_{g,\Omega_3}(\partial \Omega_1,\partial\Omega_2)=d_{g,\Omega_4}(\partial \Omega_1,\partial\Omega_2)\s \mbox{and}\s d_{g,\Omega_3}(x,\partial\Omega_1)=
d_{g,\Omega_2}(x,\partial\Omega_1).
$$
In other words, $d_{g,\Omega_3}(\partial \Omega_1,\partial\Omega_2)$ and $d_{g,\Omega_2}(x,\partial\Omega_1)$ only depend on  $\Omega_1$ and $\Omega_2$. We will denote them by $d_g(\partial \Omega_1,\partial\Omega_2)$ and $d_g(x,\partial\Omega_1)$ respectively. 
%\end{rem}

%Let $(\Sigma,g_0)$ be a Riemannian surface  and $g\in\mathcal{M}(\Sigma,g_0)$. 
Suppose that $\Omega\subset\Sigma$ is a  bounded domain  with piecewise smooth boundary.  
%By the trace embedding theorem, $\partial\Omega$ is rectifiable in $g$ and $\ell_g(\partial\Omega)<\infty$.
Define 
$$
\diam(\Omega,g)=\sup_{x,y\in\overline\Omega}d_{g,\overline{\Omega}}(x,y).
$$
Assume $x,y\in\overline{\Omega}$. If $x,y\in\partial\Omega$, we have
$
d_{g,\overline{\Omega}}(x,y)\leq\ell_g(\partial\Omega).
$
If $x\in\Omega$ and $y\in\partial\Omega$, there is a curve $\gamma:[0,1]\rightarrow\Sigma$, $\gamma(0)=x$, $\gamma(1)\in\partial\Omega$,  $\gamma((0,1))\subset\Omega$, such that
$
\ell_g(\gamma)\leq d_{g}(x,\partial\Omega)+\epsilon.
$
Then 
$$
d_{g,\overline{\Omega}}(x,y)\leq d_{g,\overline{\Omega}}(x,\gamma(1))+d_{g,\overline{\Omega}}(\gamma(1),y)\leq \ell_g(\gamma)+\ell_g(\partial\Omega)\leq d_g(x,\partial\Omega)+\ell_g(\partial\Omega)+\epsilon.
$$ 
Letting $\epsilon\rightarrow 0$ leads to 
$$
d_{g,\overline{\Omega}}(x,y)\leq d_g(x,\partial\Omega)+\ell_g(\partial\Omega). 
$$
In a similar way, when $x,y\in\Omega$, we have
$$
d_{\overline{\Omega}}(x,y)\leq d_g(x,\partial\Omega)+d_g(y,\partial\Omega)+\ell_g(\partial\Omega). 
$$
Hence,
\begin{equation}\label{diam.estimate}
\diam(\Omega,g)\leq 2\sup_{x\in\Omega}d_g(x,\partial\Omega)+\ell_g(\partial\Omega).
\end{equation}

%%Let $(\Sigma,g_0)$ be a Riemann surface,  $g\in\mathcal{M}_{loc}(\Sigma,g_0)$ and $\Omega$ a relatively compact domain in $\Sigma$.  A geodesic neighborhood of $\Omega$ is a neighborhood $U$ of $\Omega$ relatively compact in $\Sigma$ with 
%$$
%d_{g,U}|_{\Omega\times\Omega}=d_{g,\Sigma}|_{\Omega\times\Omega},
%$$
%i.e., for any $x$,$y\in\Omega$, there exists a piecewise smooth curve $\gamma_k$ from $x$ to $y$ in $U$ and 
%$$
%\ell(\gamma_k)\rightarrow d_{g,\Sigma}(x,y).
%$$ 
%We say $U$ is a quasiconvex neighborhood  of $\Omega$, if 
%$$d_g(\partial U,\partial\Omega)>diam(\Omega,d_{g,U}).$$
%\end{defi}

\begin{defi}
Let $(\Sigma,g_0)$ be a Riemannian surface,  $g\in\mathcal{M}(\Sigma,g_0)$ and $\Omega$ a connected relatively compact domain in $\Sigma$. A {\bf quasi-geodesic convex neighbourhood} of $\Omega$ is a neighbourhood $U$ of $\Omega$ relatively compact in $\Sigma$, which satisfies:
for any $x$, $y\in\Omega$, and a curve $\gamma$ from $x$ to $y$ in $\Sigma$,
there exists another curve $\gamma'$ from $x$ to $y$, such that 
$\gamma'\subset U$ and $\ell_g(\gamma')\leq\ell_g(\gamma\cap U).$
\end{defi}
Obviously, when $U$ is a quasi-geodesic convex neighborhood of $\Omega$,  there holds
$$
d_{g,U}|_{\Omega\times\Omega}=d_{g,\Sigma}|_{\Omega\times\Omega},
$$
and this relation is not true if $d_{g,U}$ is replaced by $d_{g,\Omega}$, e.g. $\Omega$ is a nonconvex domain in $\R^2$. A similar idea is contained in \cite[Lemma 2.2.1]{R}. 

The metric surface $(\Sigma,g)$ is {\it complete} if any Cauchy sequence with respect to $d_{g,\Sigma}$ converges. 
Observed basic facts include (cf. \cite{AZ-1962}, \cite[2.2-2.3]{R}): for the intrinsic metric $g$, the closed ball  $\overline{B}_r(x),x\in\Sigma$ is compact and completeness implies that any two points can be joined by a shortest curve in $\Sigma$.

\begin{lem}\label{suff.geodesic}
Let $(\Sigma,g_0)$ be a Riemannian surface and $g\in\M(\Sigma,g_0)$.  Let $\Omega$ be a relatively compact domain in $\Sigma$ with piecewise smooth boundary and $U$ a neighborhood of $\Omega$ relatively compact in $\Sigma$.  
If $d_{g}(\partial\Omega,\partial U)>\diam(\Omega,d_{g,\Sigma})$ then $U$ is a quasi-geodesic convex neighborhood of $\Omega$.
\end{lem}
\proof
If there existed $x,y\in\Omega$ so that $d_{g,U}(x,y)\not=d_{g,\Sigma}(x,y)$,  by (iv) in Lemma \ref{Sigma12} we would have $d_{g,\Sigma}(x,y)\geq d_{g}(\partial \Omega,\partial U)>\diam(\Omega,d_{g,\Sigma})$, contradicting the definition of diameter. 
\endproof
%1) $F_{U,\Omega}(d_{g,\Sigma})>0$ if and only if $F_{U,\Omega}(d_{g,U})>0$.
%2) if $F_{\Omega,U}(d_{g,U})>0$ (or $F_{U,\Omega}(d_{g,\Sigma})>0$),
%then $U$ is a geodesic neighborhood of $\Omega$ in $(\Sigma,g)$.

\vspace{.1cm}

For $g=e^{2u}g_0\in\M(\Sigma,g_0)$, from the discussion above, $d_{g,\Sigma}(x,y)$ is finite for any  $x,y\in\Sigma \backslash A_{\frac{4}{3}\pi}$,  where $A_{\frac{4}{3}\pi}$ is the discrete set where $|\K_g|$ concentrates as in \eqref{defi.A}. By Corollary \ref{defi.dis2}, $d_{g,\Sigma}$ can be realized in a disk $D_{r_0}$ provided $|\K_g|(D)$ is small.% < \tau_1$ with $r_0,\tau_1$ depending on $A$ for $\|\nabla u\|_{L^1(D)}<A$. 

%Next, we observe that $e^{u+\varphi}\in L^q$ by Corollary \ref{MTI} hence $e^u\in L^q$ as $\varphi$ is smooth. 
%Finally, since the Lebesgue points exist $g_{euc}$ a.e. we conclude $d_{g,\Sigma}(p,\cdot)\in W_{loc}^{1,q}(\Sigma\backslash\cap_{\epsilon>0}A_\epsilon)$. \endproof

%\vspace{.1cm}

\subsection{A distance comparison theorem for small total curvature measure}
%We focus on controlling deviation of the distance function from the euclidean distance if the total curvature on a disk is small.  For the purpose of proving our main theorems 
%We are interested in precise estimates in terms of $u$ in the conformal factor $e^{2u}$.  
Let $u_{x,r}$ denote the average of $u$ over $D_r(x)\subset D$. 

%Denote $$u_{x,r}=\frac{1}{\pi r^2}\int_{D_r(x)}u(y)dy.$$

%We begin with a few preparatory results and then present the proof of the theorem in due course.

\begin{lem}\label{lower.curve}
Assume $g=e^{2u}g_{\rm euc}\in\M(D)$ 
with $\|\nabla u\|_{L^1(D)}<A$. For $q\in (1,2)$ and any $\epsilon>0$, there is $c=c(\epsilon,q)>0$ so that for any piecewise smooth curve $\gamma$ from $0$ to $x$ in $D$ it holds
$$ 
\int_\gamma e^u\geq e^{u_{0,|x|}}|x|e^{-c|x|^{1-2/q}\|\nabla u\|_{L^q(D_{2|x|})}-\epsilon}.
$$
\end{lem}

%\proof Since $1<q<2$, we can select $s$ so that $2-q<s<1$. Taking $\epsilon=1$ in Lemma \ref{diff.of.d.1}, we see $\mathcal H^s(S) \leq \Lambda'$, where $S=\{y\in D:|u(y)-u_{0,|x|}|>\|\nabla u\|_{L^q(D)}\}$. Now \begin{align*} \int_{\gamma}e^u &= e^{u_{0,|x|}} \int_\gamma e^{u-u_{0,|x|}} \\&\geq e^{u_{0,|x|}}\int_{\gamma}e^{-|u-u_{0,|x|}|} \\&\geq e^{u_{0,|x|}}\int_{\gamma\backslash S}e^{-|u-u_{0,|x|}|} \\&\geq e^{u_{0,|x|}} \int_{\gamma\backslash S} e^{-\|\nabla u\|_{L^q(D)}}\\&\geq e^{u_{0,|x|}} e^{-\|\nabla u\|_{L^q(D)}} |x|\end{align*}

\proof 
Fix $x\in D$ and set $r=|x|$. Let $t_1=\min\left\{t:|\gamma(t)|=r\right\}$, $\gamma_1=\gamma|_{[0,t_1]}$ and $x_1=\gamma(t_1)$. Since
$$
\int_\gamma e^u\geq \int_{\gamma_1}e^{u},
$$
it suffices to prove that there exists $c=c(\epsilon,q)$ such that 
$$ 
\int_{\gamma_1} e^u\geq e^{u_{0,r}}re^{-c r^{1-2/q}\|\nabla u\|_{L^q(D_{2r})}-\epsilon}.
$$

Applying Lemma \ref{diff.of.d.1} to
$u'=u(\frac{x}{r})$, we see that for any $\epsilon_1>0$ there is $\lambda$ so that
\begin{align*}
\mathcal{H}_\infty^1&\left(\left\{y\in D_r: |u(y)-u_{0,{r}}| >\lambda r^{1-\frac{2}{q}}\|\nabla u\|_{L^q(D_{2r})}\right\}\right)\\
&=r\,\mathcal{H}_\infty^1\left(\left\{y\in D: |u'(y)-{u'}_{0,{1}}| >\lambda \,\|\nabla u'\|_{L^q(D_{2})}\right\}\right)\leq r\epsilon_1.
\end{align*}
Denote 
$$
S_\lambda=\left\{y\in D_r:|u(y)-u_{0,r}|\leq\lambda r^{1-\frac{2}{q}}\|\nabla u\|_{L^q(D_{2r})}\right\},\s T_\lambda=D_r\backslash S_\lambda.
$$
We have
\begin{align*}
\int_{\gamma_1}e^{u}&=e^{u_{0,r}}\int_{\gamma_1} e^{u-u_{0,r}}\geq e^{u_{0,r}}\int_{\gamma_1} e^{-|u-u_{0,r}|}\\
&\geq e^{u_{0,r}}\int_{\gamma_1\cap S_\lambda}e^{-\lambda r^{1-\frac{2}{q}}\|\nabla u\|_{L^q(D_{2r})}}
\geq  e^{u_{0,r}}e^{- \lambda r^{1-\frac{2}{q}}\|\nabla u\|_{L^q(D_{2r})}}\mathcal{H}^1(\gamma_1\cap S_\lambda).
\end{align*}

To estimate the 1-dimensional Hausdorff measure of $\gamma_1\cap S_\lambda$,  let $\pi$ be the orthogonal projection from $\R^2$ to the straight line passing through $0$ and $x_1$. By \cite[Proposition 3.5]{Maggi},  we have
\begin{eqnarray*}
\mathcal{H}^1(\gamma_1\cap S_\lambda)\geq\mathcal{H}^1(\pi(\gamma_1\cap S_\lambda)),\s \mathcal{H}^1(\gamma_1\cap T_\lambda)\geq\mathcal{H}^1(\pi(\gamma_1\cap T_\lambda))
\end{eqnarray*}
and
\begin{eqnarray*}
\mathcal{H}_\infty^1(\gamma_1\cap T_\lambda)\geq\mathcal{H}_\infty^1(\pi(\gamma_1\cap T_\lambda)). 
\end{eqnarray*}
Noting that
$$
\pi(\gamma_1\cap S_\lambda)
\cup \pi(\gamma_1\cap T_\lambda)\,{\supset\overline{0x_1}},
$$
we have
$$
\mathcal{H}^1(\gamma_1\cap S_\lambda)\geq r-\mathcal{H}^1(\pi(\gamma_1\cap T_\lambda))
=r-\mathcal{H}^1_\infty(\pi(\gamma_1\cap T_\lambda))
\geq r-\mathcal{H}^1_\infty(\gamma_1\cap T_\lambda)
\geq r(1-\epsilon_1), 
$$
here we used \cite[Proposition 3.5 and Theorem 3.10]{Maggi}. Now, 
$$
\int_\gamma e^u \geq e^{u_{0,r}}e^{- \lambda r^{1-\frac{2}{q}}\|\nabla u\|_{L^q(D_{2r})}}r{(1-\epsilon_1)}.
$$
Choose $\epsilon_1$ so that $e^{-\epsilon}=(1-\epsilon_1)$, and
take $c=\lambda$. \endproof

\begin{thm}\label{estimate.d_g} 
Assume $g=e^{2u}g_{\rm euc}\in\M(D)$ 
with $\|\nabla u\|_{L^1(D)}<A$. Then for any $\epsilon>0$ there are constants 
$r(\epsilon),\tau(\epsilon)>0$ depending on $\epsilon,A$, such that if $|\K_g|(D)<\tau(\epsilon)$ then 
$$ 
e^{u_{0,|x|}-2\epsilon}\leq\frac{d_{g,D}(0,x)}{|x|}\leq e^{u_{0,|x|}+2\epsilon},\s \forall x\in D_{r(\epsilon)}.
$$
\end{thm}

%\noindent{\it Proof of Theorem \ref{estimate.d_g}.} 
\proof We write $d_g$ for $d_{g,D}$. 
By Lemma \ref{rq}, we can choose $\tau(\epsilon),r(\epsilon)$ such that 
$c|x|^{1-2/q}\|\nabla u\|_{L^q(D_{2|x|})}<\epsilon$  for $x\in D_{r(\epsilon)}$. By Lemma \ref{lower.curve},  
$$
\frac{d_{g,D}(0,x)}{|x|}\geq e^{-2\epsilon+u_{0,|x|}}, \ \ \ |x|<r(\epsilon). 
$$

Next, we prove the other inequality by contradiction.
Suppose there exist $\epsilon>0$, $u_k$ and $x_k\in\partial D_{r_k}$ with $r_k=|x_k|\rightarrow 0$, such that 
$$
|\K_{g_k}|(D)\rightarrow  0\s\mbox{ and }\s
\frac{d_{g_k,D}(0,x_k)}{r_k}>e^{\epsilon+(u_k)_{0,r_k}}.
$$
By Lemma \ref{rq}
$$
r_k^{1-2/q}\|\nabla u_k\|_{L^q(D_{2r_k}(x))}\leq
C(Ar_k+|\K_{g_k}|(D))\rightarrow 0, \forall x\in D_\frac{1}{4}.
$$
Let $u_k'(x)=u_k(r_kx)-(u_k)_{0,r_k}$. The above inequality and the Poincar\'e inequality imply $u_k'\to
0$ in $W^{1,q}_{\rm loc}(\R^2)$. We fix a $q\in(1,2)$. By Corollary \ref{MTI},  $e^{|u_k'|}$ is bounded in $L^{q'}(D_2)$ where $q'=\frac{q}{q-1}$. By the mean value theorem, $|e^y -1|\leq e^{|y|}|y-0|$ and let $y=u'_k(x)$. 
%$|e^{u_k'(x)}-1|=|e^{u_k'(x)}-e^0|\leq e^{|u_k'(x)|}|u_k'(x)-{0}|$, which implies that
Hence
$$
\|e^{u_k'}-1\|_{L^1(D_2)}\leq\|e^{|u_k'|}\|_{L^{q'}(D_2)}
\|u_k'\|_{L^q(D_2)}.
$$
 Then $e^{u_k'}$ converges to $1$ in $L^1$. Moreover, 
$$
\int_{D_2}|\nabla e^{u_k'}|=\int_{D_2}e^{u_k'}|\nabla u_k'|
\leq\|e^{u_k'}\|_{L^{q'}(D_2)}\|\nabla u_k'\|_{L^q(D_2)}
\rightarrow 0.
$$
Then $e^{u_k'}$converges to $1$ in $W^{1,1}(D_2)$. 
Applying the fact that the trace operator is compact (cf. \cite[Corollary 18.4]{Leoni}) to $e^{u_k'}-1$, we have
$$
\frac{d_{g_k,D}(0,x_k)}{e^{(u_k)_{0,r_k}}r_k}\leq \frac{1}{e^{(u_k)_{0,r_k}}r_k}\int_{\overline{ox_k}}e^{u_k}
= \int_0^{\frac{x_k}{r_k}}e^{u_k'}\rightarrow 1.
$$
Then $$\frac{d_{g_k,D}(0,x_k)}{e^{(u_k)_{0,r_k}}r_k}<e^\epsilon$$ for large $k$, therefore 
$$
\frac{d_{g_k,D}(0,x_k)}{r_k}<e^{\epsilon+(u_k)_{0,r_k}}
$$
but this contradicts the choice of $\epsilon, u_k, x_k$. 
\endproof

A consequence of Theorem \ref{estimate.d_g} is that {\it length minimizing} is realized locally when the total Gauss curvature measure is small.   We set
$$
C_P=\inf_{u\in W^{1,3/2}(D),\int_{D_{1/2}}u=0}\frac{\|\nabla u\|_{L^{3/2}(D)}}{\|u\|_{L^1(D)}}.
$$
By the Poincar\'e inequality (cf. \cite[Theorem 5.4.3]{Attouch-Buttazzo-Michaille}), $C_P>0$.

\begin{cor}\label{defi.dis2} 
Assume $g=e^{2u}g_{\rm euc}\in\M(D)$ 
with $\|\nabla u\|_{L^1(D)}<A$. Let $c_0=\min\left\{\frac{4}{3}\pi,\frac{\log 2}{C(\frac{3}{2})C_P}\right\}$, where $C(\frac{3}{2})$ is as in \eqref{II} with $q=\frac{3}{2}$. There exists $r_0$   depending on $A$, such that if $|\K_g(D)|<c_0$ then
for any $x\in D_{r_0}(y)\subset D_{2r_0}(y)\subset D_{1/2}$ 
\begin{equation}\label{dist(0,x).small.domain}
d_{g,D}(y,x)=d_{g,D_{2r_0}(y)}(y,x).%=\inf\left\{\int_\gamma e^u: \gamma\subset D_{2r_0},\gamma(0)=0,\gamma(1)=x\right\}=d_{g,D_{2r_0}}(0,x).
\end{equation}
Moreover, for any $\Omega\subset D$, if $D_{2r_0}\subset\Omega$ then
\begin{equation}\label{dist.small.domain}
d_{g,\Omega}|_{D_{r_0/4}}=d_{g,D_{2r_0}}|_{D_{r_0/4}}.
\end{equation}
\end{cor}

\proof We establish \eqref{dist(0,x).small.domain} first. 
It is obvious that $d_{g,D}(y,x)\leq d_{g,D_{2r_0}(y)}(y,x)$ as long as $D_{2r_0}(y)\subset D$. We now argue the other direction. 
Choose $\epsilon$ with $2e^{-2\epsilon}>1$ and let $r(\epsilon),\tau(\epsilon)$ be as in Theorem \ref{estimate.d_g}. We begin with $2r_0<r(\epsilon)$. 
Suppose that $\gamma\subset D$ is a curve from $\gamma(0)=y$ to $\gamma(1)=x\in D_{r_0}(y)$ but $\gamma\not\subset D_{2r_0}(y)$. Assume $t_0$ be the first $t$ such that $|\gamma(t_0)|=2r_0$. By Theorem \ref{estimate.d_g},
$$
2r_0e^{u_{y,2r_0}-\epsilon}\leq d_{g,D}(y,\gamma(t_0))\leq\int_\gamma e^u
$$ 
and 
$$
d_{g,D}(y,x)\leq  e^{u_{y,r_0}+\epsilon}r_0.
$$
By the Poincar\'e inequality (cf. \eqref{repeat}) and Lemma \ref{rq}, we have
$$
|u_{y,2r_0}-u_{y,r_0}|\leq c_pr_0^{-1/3}\|\nabla u\|_{L^{3/2}(D_{2r_0}(y))}
<C_P\left(c_0C(\frac{3}{2})+CAr_0\right).
$$
Then
$$
2e^{-2\epsilon}d_{g,D}(y,x) \leq 2r_0 e^{u_{y,r_0}-u_{y,2r_0}+u_{y,2r_0}-\epsilon}\leq e^{{C_P(c_0C(\frac{3}{2})+CAr_0)}}\int_\gamma e^u.
$$
Choosing  $r_0<\frac{r(\epsilon)}{2}$ small with
$
2e^{-2\epsilon} e^{-{C_P(c_0C(\frac{3}{2})+CAr_0)}} >1.
$
%then 
%$$
%d_{g,D}(0,x)< \int_\gamma e^u.
%$$
So $d_{g,D}(y,x)$ can only be realized by curves in $D_{2r_0}(y)$. 

Next, we see that \eqref{dist.small.domain} follows from \eqref{dist(0,x).small.domain} and 
$$
d_{g,D}(x_0,y_0)\leq d_{g,\Omega}(x_0,y_0)\leq d_{g,D_{2r_0}}(x_0,y_0)\leq d_{g,D_{r_0}(x_0)}(x_0,y_0)
$$ 
for any $x_0,y_0\in D_{r_0/4}$. \endproof

%Next, we prove \eqref{dist.small.domain}. By a rescaling argument, for any $x_0,y_0\in D_{r_0/4}$, and curve $\gamma$ from $x_0$ to $y_0$, if $\gamma$ is not in $D_{r_0}(x_0)$, then $\int_\gamma e^u>2e^{-\epsilon}d_{g,D_\frac{1}{2}(x_0)}(x_0,y_0)$. Note that $d_{g,D_\frac{1}{2}(x_0)}\geq d_{g,D}(x_0,y_0)$. Then $$d_{g,D_{r_0}(x_0)}(x_0,y_0)=d_{g,D}(x_0,y_0).$$ Since $$d_{g,D}(x_0,y_0)\leq d_{g,\Omega}(x_0,y_0)\leq d_{g,D_{2r_0}}(x_0,y_0)\leq d_{g,D_{r_0}(x_0)}(x_0,y_0), $$we get $d_{g,\Omega}(x_0,y_0)=d_{g,D_{2r_0}}(x_0,y_0)$. \endproof 

\vspace{.1cm}

Recall that when $g$ is smooth, the distance function $d_g(p, x)$ is Lipschitz and $|\nabla^{g} d_g(p, x)| = 1$ almost everywhere for $x$; hence if $g=e^{2u}g_0$, where $u$ and $g_0$ are smooth,  then $|\nabla^{g_0} d_g(p, x)| = e^{u}$ almost everywhere. 
For nonsmooth $g$, it is known that the components of $g$ belongs to some Sobolev space if the curvature of $g$ is bounded below in the sense of Alexandrov \cite[Proposition 2.8]{Ambro-Bert}.  We now show that $d_{g,\Sigma}$ is a Sobolev function on $\Sigma$ with finitely many points removed with the same estimate.  

\begin{lem}\label{Sobolev.distance}
Assume $(\Sigma,g_0)$ is a smooth surface and  $g=e^{2u}g_{\rm euc}\in \M(\Sigma,g_0)$. For any $\tau< 4\pi$,  it holds $d_{g,\Sigma}(p,\cdot)\in W^{1,q}_{\rm loc}(\Sigma\backslash A_\tau)$ for any $q\in[1,\frac{4\pi}{\tau})$ where $A_\tau=\left\{x: |\K_g|(\{x\})\geq\tau\right\}$. Moreover,
$$
|\nabla^{g_0} d_{g,\Sigma}(p,x)|\leq e^{u(x)},
$$
for a.e.  $x$ measured in $g_0$. In particular, $d_{g,\Sigma}$  is continuous on $\Sigma\setminus A_{2\pi}$.
\end{lem}  

%{\color{red}
%For simplicity, we set $\tau_1=\min\{\tau_0,\pi\}$. By Corollary \ref{MTI}, if $\|u\|_{L^1(\Sigma)}<\gamma$, then for any $\Omega\subset\subset\Sigma\backslash A_{\tau_1}$,
%$$
%\|d_{g,\Sigma}(p,x)\|_{W^{1,3}(\Omega)}<C(\gamma,\Omega).
%$$
%}

%\noindent{\it Proof of Lemma \ref{Sobolev.distance}.} 
\proof Take an isothermal coordinate system $(D,x)$ on $\Sigma\backslash A_{\frac{4}{3}\pi}$ around a point $p$. Assume 
$$
g_0=e^{2\varphi}g_{\rm euc}\s \mbox{and} \s w=u+\varphi.
$$
Suppose $0$ is a Lebesgue point of $u$ w.r.t. $g_{\rm euc}$ away from $\cup_{\epsilon>0}A_\epsilon$. Let $\epsilon_k\to 0$ and $\tau(\epsilon_k), r(\epsilon_k)$ be as in Theorem \ref{estimate.d_g}. There exist $r_k$ such that $|\mu_g|(D_{r_k})<\tau(\epsilon_k)/2$ and we assume $r_k\leq r(\epsilon_k)$. 
	Set $g_k = r_k^{-2} g$ and $w_k(x) =w(r_k{x})+\log r_k$. Then $|\mu_{g_k}|(D)<\tau(\epsilon_k)<{A}/{2}$. By \eqref{II} (cf. Lemma \ref{rq}) 
	$$
	\|\nabla w_k\|_{L^1(D)}= r_k^{-1}\|\nabla w\|_{L^1(D_{r_k})}\leq r^{-1}_kC(Ar^2_k +r_k\mu_{g_k}(D))
	\leq Cr_kA + \tau(\epsilon_k)<A
	$$ 
	by further shrinking $r_k$ if necessary. By Theorem \ref{estimate.d_g}
	$$
	e^{(w_k)_{0,|x|}-\epsilon_k}\leq\frac{d_{g_k,D}(0,x)}{|x|}\leq e^{(w_k)_{0,|x|}+\epsilon_k},\s \forall x\in D_{r_k}.
	$$
	Moreover, any curve $\gamma$ from $0$ to $x$ corresponds to a curve $\gamma_k$ from $0$ to ${r_k}x$ and vice versa, so 
	$$
	d_{g_k,D}(0,x) = \inf_\gamma \int_\gamma e^{w_k(y)}dy= \frac{1}{r_k}\inf_{\gamma_k} \int_{\gamma_k} e^{w({r_k}y)}d\left({r_k}y\right)= \frac{1}{r_k} d_{g,\Sigma}\left(0,{r_k}x\right)
	$$
	for  $x\in D_{r_k}\subset D_{r_0}$, 
	and 
	$$
	(w_k)_{0,|x|} = \avint_{D_{|x|}}w_k(y)dy =  \avint_{D_{|x|}} w\left({r_k}y\right)dy =\avint_{D_{{r_k}|x|}}w(z)dz= w_{0,{r_k}|x|}.
	$$
	%$$(w_k)_{0,|x|} = \frac{1}{\pi |x|^2}\int_{D_{|x|}}w_k(y)dy = \frac{1}{\pi |x|^2} \int_{D_{|x|}} w\left({r_k}y\right)dy =\frac{1}{\pi |x|^2r^2_k}\int_{D_{{r_k}x}}w(z)dz= w_{0,{r_k}|x|}.$$
	Hence
	$$
	e^{w_{0,r_k|x|}-\epsilon_k}\leq\frac{d_{g,\Sigma}(0,r_kx)}{r_k|x|}\leq e^{w_{0,{r_k|x|}}+\epsilon_k},\s \forall x\in D_{r_k}\subset D_{r_0}.
	$$
	This can be rewritten as 
	$$
	e^{w_{0,|x|}-\epsilon_k} \leq \frac{d_{g,\Sigma}(0,x)}{|x|}\leq e^{w_{0,|x|}+\epsilon_k}, \ \ \ \forall x\in D_{r_k^{2}}.
	$$
	Thus, since $0$ is a Lebesgue point, 
	$$
	e^{w(0)-\epsilon_k}\leq \liminf_{x\to 0}\frac{d_{g,\Sigma}(0,x)}{|x|} \ \mbox{and}\ \limsup_{x\to 0}\frac{d_{g,\Sigma}(0,x)}{|x|}\leq e^{w(0)+\epsilon_k}.
	$$
	Letting $\epsilon_k\to 0$ we see 
	$$
	\lim_{x\to 0} \frac{d_{g,\Sigma}(0,x)}{|x|} = e^{w(0)}.
	$$
	It follows that at any Lebesgue point $x$ we have 
	$$
	\lim_{h\to 0}\frac{\left|d_{g,\Sigma}(p,x+h)-d_{g,\Sigma}(p,x)\right|}{|h|}\leq \lim_{h\to 0}\frac{d_{g,\Sigma}(x+h,x)}{|h|}=e^{w(x)}.
	$$
	This shows that {if $d_{g,\Sigma}$ is differentiable at a Lebesgue point $x$ of $u$, then}
	$$
	\left|\nabla^{g_0}_x d_{g,\Sigma}(p,x)\right| = e^{-\varphi}\left|\nabla_x d_{g,\Sigma}(p,x)\right|\leq e^{w(x)-\varphi(x)}= e^{{u(x)}}. 
	$$

Now we show that $d_{g,\Sigma}$ has weak derivative in $L^q$. {Recall that a two dimensional $W^{1,p}$ function is differentiable almost everywhere for $p>2$ (cf. \cite[Theorem 6.5]{E-G})}. 
Let $h\in\R^n$ with $|h|=r<\frac{1}{2}$. Then for any $x\in D_{{1}/{2}}$,
$$
\left|d_{g,\Sigma}(p,x+h)-d_{g,\Sigma}(p,x)\right|\leq \int_{[x,x+h]}e^{w}=r
\int_{[0,h/|h|]}e^{w(x+ry)}dy.
$$
By Lemma \ref{rq},
$$
\|\nabla w(x+ry)\|_{L^1(D_2)}=r^{-1}\|\nabla w\|_{L^1(D_{2r}(x))}<C.
$$
By Corollary \ref{MTI}, $e^{w}\in L^q$. Applying the trace embedding theorem and the Poincar\'e inequality to $w(x+ry)$, 
$$
\int_{[0,h/r]}e^{w(x+ry)}\leq Ce^{\overline{w(x+ry)}}=Ce^{w_{x,|h|}}.
$$
Then 
$$
\left| d_{g,\Sigma}(p,x+h)-d_{g,\Sigma}(p,x) \right|\leq C e^{w_{x,|h|}}|h|.
$$
We have
\begin{align*}
\int_{D_{1/2}}e^{qw_{x,r}}dx&=
\int_{D_{1/2}}e^{\frac{q}{\pi r^2}\int_{D_r(x)}w(y)dy}dx=\int_{D_{1/2}}e^{\frac{q}{\pi }\int_{D}w(x+ry)dy}dx\\
&\leq\frac{1}{\pi}\int_{D_{1/2}}\int_{D}e^{qw(x+ry)}dydx\ \ \ \ \s(\mbox{Jensen's inequality})\\
&=\frac{1}{\pi}\int_{D}\int_{D_{1/2}}e^{qw(x+ry)}dxdy\leq\frac{1}{\pi}\int_{D}\int_{D}e^{qw(x)}dxdy=\int_{D}e^{qw(x)}dx.
\end{align*}
Then by  \cite[Proposition 9.3]{Brezis}, we conclude $d_{g,\Sigma}(p,\cdot)\in W^{1,q}(D_{1/2})$. When $\tau>2\pi$ we know $q>2$ and then $d_{g,\Sigma}\in C^0$ from the Sobolev embedding theorem. 
\endproof

\section{Convergence of distance functions with fixed conformal class}

We first introduce a set of finite ordered lists of points in a metric space such that any pair of adjacent points in a list are separated by at least a positive distance $a$ but not by $2a$.  Different lists may have different number of points. The points from a list (called $a$-string) will be used as endpoints of a polygonal curve.

\begin{defi}
In a metric space $(X,d)$, a finite collection of points is called an {\bf $a$-string} if the distance between any two  adjacent points in the collection lies in $[a,2a]$, $a>0$.  The set of all $a$-strings is denoted 
$$
\Gamma_a(X,d)=\left\{(x_0,x_1,\cdots,x_m):  x_i\in X,\s a\leq d(x_i,x_{i-1})\leq 2a,\s m\in \mathbb{Z}^+\right\}.
$$
For $\alpha=(x_0,\cdots,x_m)\in\Gamma_a(X,d)$, we define
$
\alpha^-=x_0, \alpha^+=x_m, \ell(\alpha)=m,
$
and for an arbitrary distance $d'$ on $X$ the $d'$-length of the $a$-string $\alpha$ by 
$$
\mathcal{L}_{d'}(\alpha)=\sum_{i=1}^m d'(x_i,x_{i-1}).
$$
\end{defi}

We will use $\ell_g(\gamma)$ to denote the length of a curve $\gamma$ in a metric $g$. 

\begin{lem}\label{disg}
Let $\gamma$ be a curve on a Riemannian surface $(\Sigma,g_0)$ parametrized by $t\in[0,1]$. 
If $\ell_{g_0}(\gamma)\in(a,+\infty)$,
then there exists $\alpha\in\Gamma_a(\Sigma,d_{g_0})$ such that
$\alpha\subset\gamma$ and $\alpha^-=\gamma(0)$, $\alpha^+=\gamma(1)$.
\end{lem}

\proof
Let $t_1$ be the smallest $t$ such that $d_{g_0}(\gamma(t_1),\gamma(0))=a$, and $t_2$ be the smallest $t>t_1$ with 
$d_{g_0}(\gamma(t_2),\gamma(t_1))=a$. Repeat this whenever possible to get $t_1, \cdots, t_{m}$ with $d_{g_0}(\gamma(t_i),\gamma(t_{i-1}))=a$ for $i=1,\cdots,m-1$ 
and $d_{g_0}(\gamma(t_m),\gamma(1))\leq a$. % or $\gamma(t_m)=\gamma(1)$.  
Clearly $m$ depends on $a$. Set $x_0=\gamma(0)$, $x_m=\gamma(1)$, $x_i=\gamma(t_i)$ for $i=1,\cdots,m-1$. Then 
$\alpha=(x_0,\cdots,x_m)$ fulfills the requirement. 
\endproof

\subsection{Singular metrics with small total curvature measure}
The main result in the subsection is: 
\begin{pro}\label{conv0}
Let $g_k=e^{2u_k}g_{\rm euc}\in\M(D)$ and $g=e^{2u}g_{\rm euc}\in\M(D)$. Assume that $\|\nabla u\|_{L^1(D)}<A$,
 $|\K_{g_k}(D)|<c_0$  and $u_k$ converges  to $u$ in $L^1_{\rm loc}(D)$. 
Then $d_{g_k,D}$ converges to $d_{g,D}$ in $C^0(D_{\frac{r_0}{4}}\times D_{\frac{r_0}{4}})$, where $c_0$ and $r_0$ are as in Corollary \ref{defi.dis2}. 
\end{pro}

\proof 
We have $\frac{4\pi}{c_0}>2$ . By Corollary \ref{MTI} and  Lemma \ref{Sobolev.distance}, for any fixed $r\in(0,1)$,
$d_{g_k,D}$ is bounded in $W^{1,q}(D_r\times D_r)$ for some $q>2$.
Then $d_{g_k,D}$ converges in $C^{0,\sigma}_{\rm loc}$ to a function $d$ by the Sobolev embedding theorem. 

\vspace{.1cm}
 
First, we show $d$ is a distance function. The triangle inequality and the symmetry follow from that $d$ is the limit of distance functions $d_{g_k,D}$. 
So it suffices to prove $d(x_0,y_0)>0$ for any $x_0\neq y_0$. Without losing generality we assume $y_0=0$. Let $\gamma$ be a curve from $0$ to $x_0$ in $D$. 
By  Lemma \ref{W1q.convergence} $u_k \rightharpoonup u$ in $W_{\rm loc}^{1,q}(D)$. Then by Lemma \ref{rq}
$$
|x_0|^{1-\frac{2}{q}}\|\nabla u_k\|_{L^q(D_{2|x_0|})}<C(|x_0|\|\nabla u\|_{L^1(D)}+C(q)c_0).
$$
Then, from Theorem \ref{estimate.d_g} and $u_k\rightarrow u$ in $L^1_{\rm loc}(D)$, we deduce 
\begin{equation}\label{mid step}
d_{g_k,D}(0,x_0) \geq e^{(u_k)_{0,|x_0|}}|x_0| e^{ -c|x_0|^{1-2/q}\|\nabla u_k\|_{L^q(D_{2|x_0|})}-\epsilon}\geq \delta_0(|x_0|,\|u\|_{W^{1,1}(D_{|x_0|})},\epsilon).
\end{equation}  
Therefore, $d(0,x_0)>0$.

Next, we prove $d(x,y)=d_{g,D}(x,y)$ on $D_{r_0/4}\times D_{r_0/4}$.
It suffices to show $d(0,x)=d_{g,D}(0,x)$ for any $x\in D_{r_0/4}$. The trace operator of Sobolev functions is compact (cf. \cite[Corollary 18.4]{Leoni}) and $\|e^{u_k}\|_{W^{1,q}(D_r)}$ is uniformly bounded by Corollary \ref{MTI}, thus
\begin{equation}\label{d<}
d(0,x) = \lim_{k\to\infty} d_{g_k,D}(0,x) \leq \inf_\gamma\lim_{k\to\infty}\int_\gamma e^{u_k}
= \inf_\gamma\int_\gamma e^u
= d_{g,D}(0,x).
\end{equation}

We now need to show
$$
d(0,x)\geq d_{g,D}(0,x).
$$
{\bf Step 1.} For a fixed $x$ and any $a<d(0,x)$ we claim
$$
d(0,x)=\inf\left\{\mathcal{L}_d(\alpha):\alpha\in\Gamma_a(\overline{D_{r_0/2}},d_{g_{\rm euc}}), \,\alpha^-=0,\,\alpha^+=x\right\}.
$$

It suffices to show that for any $\epsilon>0$ there exists $\alpha$ such that $\mathcal{L}_d(\alpha)\leq d(0,x)+\epsilon$.
By Corollary \ref{defi.dis2}, we can choose a curve $\gamma_k\subset \overline{D_{r_0/2}}$ with $\gamma_k(0)=0$, $\gamma_k(1)=x$, such that
$$
d_{g_k,D}(0,x)\leq \int_{\gamma_k}e^{u_k}\leq d_{g_k,D}(0,x)+\epsilon.
$$
Then $a<d_{g_k,D}(0,x)$ for all large $k$.
By Lemma \ref{disg}, we can find an $a$-string $\alpha_k=(x_0^k,\cdots,x_{m_k}^k)\in\Gamma_a(\overline{D_{r_0/2}})$, such that $\alpha_k\subset\gamma_k$ with $\alpha_k^-=0$ and
$\alpha_k^+=x$. Note 
$$
\mathcal{L}_{d_{g_k,D}}(\alpha_k)= 
\sum d_{g_k,D}(x_i^k,x_{i-1}^k)\leq\int_{\gamma_k}e^{u_k}\leq d_{g_k,D}(0,x)+\epsilon.
$$
Since $a$ is fixed  here and $x_i^k\in \overline{D_{r_0/2}}$, for $t\in[a,2a]$ we have 
$$
	|(u_k)_{x_i^k,t}|\leq C(a)\|u_k\|_{L^1(D)}<C
,$$ 
and then replacing $|x_0|$ in \eqref{mid step} by $t$ we see
$$
d_{g_k,D}(x_i^k,x_{i-1}^k)>\delta_0(a,\|u\|_{W^{1,1}(D_{r_0})}).
$$
So $m_k$ is bounded from above by a number independent of $\gamma$. Without loss of generality, we assume $m_k$ is fixed and
$\alpha_k\to\alpha_\infty$ as $k\rightarrow+\infty$, i.e.
$x_i^k\rightarrow x_i^\infty$ for each $i$. Since $d_{g_k,D}\to d$ in $C^{0,\sigma}_{\rm loc}(D\times D)$, $$d(x^k_i,x^\infty_i)\leq C|x^k_i-x^\infty_i|^\sigma\to 0\s\mbox{ and}\s 
d_{g_k,D}(x^k_i,x^k_{i+1})\rightarrow d(x^\infty_i,x^\infty_{i+1})$$ for each $i$. 
Then 
\begin{equation}\label{L_d}
\mathcal{L}_d(\alpha_\infty)=\lim_{k\rightarrow+\infty}\mathcal{L}_{d_{g_k,D}}(\alpha_k)\leq \lim_{k\rightarrow+\infty}d_{g_k,D}(0,x)+\epsilon=d(0,x)+\epsilon.
\end{equation}

\noindent{\bf Step 2.} Now we start to prove $d(0,x)\geq d_{g,D}(0,x)$. For any $\epsilon>0$, let $\tau(\epsilon),r(\epsilon)$ be as in Theorem \ref{estimate.d_g} (Note that we cannot take $\tau(\epsilon)$ independent of $\epsilon$ as we will let $\epsilon\to 0$). 
Define a finite set 
$$
\S=\left\{y\in \overline{D_{r_0/2}}:|\K_g|(\{y\})>\frac{\tau(\epsilon)}{4}\right\}.
$$
Fix a constant $\delta_0<\frac{1}{100}\cdot\mbox{(distance of points in $\S$)}$ and set
$
E_{\delta_0}=\bigcup_{y\in\S}\overline{D_{\delta_0}(y)}.
$
Select $\delta<\delta_0$ so that $|\K_g|(D_{2\delta}(y))<\frac{\tau(\epsilon)}{2}$ for any $y\in \overline{D_{2r_0}}\backslash E_{\delta_0}$ and choose $a<\min\left\{r(\epsilon)\delta,\frac{r_0}{4}\right\}$. 
Let $\alpha=(x_0,x_1,\cdots,x_m)\in\Gamma_a(\overline{D_{r_0/2}},d_{g_{\rm euc}})$ with $\alpha^{-}=0$,
$\alpha^+=x$ and
$$
d(0,x)\leq \mathcal{L}_d(\alpha)<d(0,x)+\epsilon.
$$
%Further, we may assume 
%$$
%|u_{k,x,t}-u_{x,t}|<\epsilon,\s \forall x\in D_\frac{1}{2},\s t\in[a,2a].
%$$

There are two cases. 

\vspace{.1cm}

\noindent Case 1: $E_{\delta_0}\cap\alpha=\emptyset$. By Lemma \ref{rq}, $\delta^{-1}\|\nabla u\|_{L^1(D_{\delta}(x_i))}<C$. Using Theorem \ref{estimate.d_g} for $u(\frac{x-x_i}{\delta})$, 
\begin{equation}\label{d-d_g.case1-1}
d_{g,D}(0,x)\leq \mathcal{L}_{d_{g,D}}(\alpha)\leq
\sum_{i=0}^{m-1} e^{\epsilon+u_{x_i, a_i}}a_i,
\end{equation}
where $a_i=|x_i-x_{i-1}|<\delta$. As $|\K_g|(D_{2\delta}(x_i))\leq\frac{\tau(\epsilon)}{2}$, we see
$
|\K_{g_k}|(D_{\delta}(x_i))<\tau(\epsilon)
$
when $k$ is sufficiently large. Then applying Theorem \ref{estimate.d_g} to $u_k(\frac{x-x_i}{\delta})$ we have 
\begin{equation}\label{d-d_g.case1-2}
\sum_{i=0}^{m-1}e^{\epsilon}e^{u_{x_i, a_i}}a_i
=e^{\epsilon}\lim_{k\rightarrow+\infty}\sum_{i=0}^{m-1}
e^{(u_k)_{x_i, a_i}}a_i
\leq e^{2\epsilon}\lim_{k\rightarrow+\infty}\sum_{i=0}^{m-1}d_{g_k,D}(x_i,x_{i+1})
=e^{2\epsilon}\mathcal{L}_d(\alpha).
\end{equation}
Then by \eqref{d-d_g.case1-1} and \eqref{L_d}
$$
d_{g,D}(0,x)\leq e^{2\epsilon}(d(0,x)+\epsilon).
$$

\noindent Case 2: $E_{\delta_0}\cap\alpha\neq\emptyset$. 
We define $\alpha'\subset\alpha$ via a ``shortening" procedure by removing substrings of $\alpha$ that depart $E_{\delta_0}$ and then return.
Given $D_{\delta_0}(y)$ for some $y\in\S$, we delete
$x_i$ if $x_i\in D_{\delta_0}(y)$ or if there exists
$x_{i-j}$, $x_{i-j+1}$, $\cdots$, $x_{i+j'}\subset \alpha$, such that
$x_{i-j}$ and $x_{i+j'}\in D_{\delta_0}(y)$, but 
$x_{i-j+1}$, $\cdots$, $x_{i+j-1}\notin D_{\delta_0}(y)$, see Figure 1.
Then $\alpha'$ can be divided into $\alpha_1$, $\cdots$, 
$\alpha_{m_0}$, such that 
\begin{enumerate}
\item[{(1)}] the euclidean distance from the last point in $\alpha_{i-1}$ to the first point in $\alpha_{i}$ is at most $2\delta_0+2a$,
\item[{(2)}] $E_{\delta_0}\cap\alpha_i=\emptyset$,
\item[{(3)}] $m_0-1$ is not bigger than the cardinality of $\S$.
\end{enumerate} 
Case 1 and (2) assert %The arguments for establishing \eqref{d-d_g.case1-1} and \eqref{d-d_g.case1-2} leads to 
$$
d_{g,D}(\alpha_i^-,\alpha_i^+)\leq e^{2\epsilon}\mathcal{L}_{d}(\alpha_i).
$$
Then 
$$
\mathcal{L}_d(\alpha)\geq \sum^{m_0}_{i=1} \mathcal{L}_d(\alpha_i)\geq \sum^{m_0}_{i=1} e^{-2\epsilon}d_{g,D}(\alpha_i^-,\alpha_i^+)\geq e^{-2\epsilon}\left(d_{g,D}(0,x)-
\sum^{m_0}_{i=1} d_{g,D}(\alpha_{i}^+,\alpha_{i+1}^-)\right).
$$
Recall that $d_{g,D}$ is in $C^{0,\sigma}$. So there is a constant $C$ independent of $\delta_0,a$ such that 
$$
d_{g,D}(\alpha_{i}^+,\alpha_{i+1}^-)\leq C(2\delta_0+2a)^\sigma.
$$
Therefore, by \eqref{L_d}
$$
d(0,x)+\epsilon\geq \mathcal{L}_d(\alpha)\geq e^{-2\epsilon}d_{g,D}(0,x) - m_0C(2\delta_0+2a)^\sigma.
$$
Let $\epsilon,a\to 0$ then $\delta_0\to 0$. We see $d(0,x)\geq d_{g,D}(0,x)$. 
\endproof

\scalebox{0.5}{
\begin{tikzpicture}
\path(0pt,0pt);
\draw (6,-5) circle (1.5);
\draw[line width=1pt]
(94pt, -282pt)
 -- (98pt, -257pt)
 -- (117pt, -248pt)
 -- (117pt, -222pt)
 -- (134pt, -213pt)
 -- (141pt, -188pt)
 -- (129pt, -166pt)
 -- (133pt, -138pt)
 -- (152pt, -133pt)
 -- (166pt, -115pt)
 -- (156pt, -85pt)
 -- (156pt, -51pt)
 -- (183pt, -51pt)
 -- (203pt, -65pt)
 -- (238pt, -74pt)
 -- (217pt, -91pt)
 -- (218pt, -111pt)
 -- (197pt, -130pt)
 -- (202pt, -149pt)
 -- (224pt, -159pt)
 -- (238pt, -187pt)
 -- (263pt, -199pt)
 -- (245pt, -229pt)
 -- (266pt, -248pt)
 -- (301pt, -244pt)
;
\draw (94pt, -286pt) node{$x_1$};
\draw (307pt, -244pt) node{$x_m$};
\draw (124pt, -166pt) node{$x_i$};
\draw (229pt, -159pt) node{$x_j$};
\draw (145pt, -140pt) node{$x_{i+1}$};
\draw (200pt, -154pt) node{$x_{j-1}$};
\end{tikzpicture}}
\hspace{2ex}
\scalebox{0.5}{
\begin{tikzpicture}%[overlay]
\path(0pt,0pt);
\filldraw (6,-5) circle (1pt);
\filldraw (133pt, -138pt) circle (1pt);
\filldraw (202pt, -149pt) circle (1pt);
\draw (6,-5) circle (1.5);
\draw[line width=1pt]
(94pt, -282pt)
 -- (98pt, -257pt)
 -- (117pt, -248pt)
 -- (117pt, -222pt)
 -- (134pt, -213pt)
 -- (141pt, -188pt)
 -- (129pt, -166pt);%
\draw[line width=1pt]
(224pt, -159pt)
 -- (238pt, -187pt)
 -- (263pt, -199pt)
 -- (245pt, -229pt)
 -- (266pt, -248pt)
 -- (301pt, -244pt)
;
\draw (94pt, -286pt) node{$x_1$};
\draw (307pt, -244pt) node{$x_m$};
\draw (124pt, -166pt) node{$x_i$};
\draw (229pt, -159pt) node{$x_j$};
\draw (141pt, -134pt) node{$x_{i+1}$};
\draw (202pt, -143pt) node{$x_{j-1}$};
\end{tikzpicture}}
\begin{center}
Figure 1 \\ $x_{i+1}$, $\cdots$, $x_{j-1}$ are deleted

$d(x_i,x_j)\leq d(x_i,x_{i+1})+d(x_{i+1},x_{j-1})+d(x_{j-1}+x_j)\leq 2\delta_0+2a$\\

\end{center}

%\begin{rem} For smooth $u_k$ converging weakly to $u$ in $W^{1,p},p>1$, Proposition \ref{conv0} was established in \cite[Theorem 1.6]{LST}where the method is different.  For nonsmooth metrics,  the condition $|\K_g|(\{x\})\})<2\pi,|\K_g|(\{y\})\})<2\pi$ at all $x,y$, in Theorem \ref{Resh} follows from smallness of total curvatture,  see Theorem \ref{2 option}.\end{rem}

\vspace{.1cm}

The distance convergent result on $D$ in Proposition \ref{conv0} can be used to deduce:

\begin{cor}\label{conv.length}
Let $g_k=e^{2u_k}g_0\in \M(\Sigma,g_0)$.  Assume $u_k\rightarrow u$ in $L^1_{\rm loc}(\Sigma)$ and $|\K_{g_k}|$ converges to a measure $\nu$ weakly. 
Then 
\begin{enumerate}
\item[(i)] Assume $K\subset\Sigma$ is compact and $\nu(\{x\})<c_0$ in $K$ where $g=e^{2u}g_{0}$ on $\Sigma$.
 For any  $\gamma_k\subset K$, if $\gamma_k(0)\rightarrow x_0$, 
$\gamma_k(1)\rightarrow y_0$ as $k\to\infty$ and $\ell_{g_k}(\gamma_k)<C$ for some constant $C$,  then
$$
\liminf_{k\rightarrow+\infty}\ell_{g_k}(\gamma_k)\geq d_{g,\Sigma}(x_0,y_0).
$$
\item[(ii)] Let $U,V$ be compact domains in $\Sigma$ and $U\subset V$. If $\nu(\{x\})<c_0$ in  $\overline{V\backslash U}$,  then 
$$
d_{g_k}(\partial U,\partial V)\rightarrow d_g(\partial U,\partial V).
$$
\end{enumerate}
\end{cor}

\proof (i)
We can use the $a$-string to localize to disks where Proposition \ref{conv0} is applicable. 
From \eqref{dist.small.domain} and Proposition \ref{conv0},
for any $x\in K$,  by selecting a suitable conformal chart w.r.t. $g_0$, we can find  $\delta_x$ so that $d_{g_k,\Sigma}(y,z)\rightarrow d_{g,\Sigma}(y,z)$
when $y$, $z\in B_{\delta_x}(x,g_0)$. Choose $x_1$, $\cdots$, $x_N$ s.t. 
$$
K\subset\bigcup_{j=1}^NB_{\delta_{x_j}/4}(x_j,g_0).
$$
Set $\delta=\min_j\delta_{x_j}$. For any  $d_{g_0}(x,y)<\frac{\delta}{2}$, we can 
find $x_{j}$ such that $x$, $y\in B_{\delta}(x_j,g_0)$. This implies 
\begin{equation}\label{d=d_g.loc.K}
d_{g,\Sigma}(x,y)\rightarrow d(x,y), \s\mbox{when $x,y\in K$ and $d_{g_0}(x,y)<\frac{\delta}{2}$.}
\end{equation}

Let $0<a<\frac{\delta}{4}.$ We claim

\begin{equation}\label{inf.osc.d}
\rho=\inf_k\inf_{x\in K}d_{g_k}(x,\partial B_{a}(x,g_0))>0.
\end{equation} 
If not, we could find $x_k\in K$ and $y_k\in \partial B_{a}(x_k,g_0)$ such that $d_{g_k,\Sigma}(x_k,y_k)
\rightarrow 0$. As $K$ is compact, we assume $x_k\rightarrow x_\infty$ and $y_k\rightarrow y_\infty$.  Then $d_{g_0}(x_\infty,y_\infty)=a$. We see a contradiction from \eqref{d=d_g.loc.K}.

\vspace{.1cm}

%Without loss of generality,  
We assume $x_0\neq y_0$ and $\liminf_{k\to+\infty}\ell_{g_k}(\gamma_k)=l_0\in(0,+\infty)$. By Lemma \ref{disg}, we can choose $\alpha_k=(x_0^k,\cdots,x_{m_k}^k)\in \Gamma_a(K,g_0)$ with $x_0^k=\gamma_k(0)$, $x_m^k=\gamma_k(1)$ and $\alpha_k\subset\gamma_k$. Obviously,
\begin{equation}\label{obvious}
\mathcal{L}_{d_{g_k,\Sigma}}(\alpha_k)=\sum_{i=1}^{m_k}d_{g_k,\Sigma}(x_k^{i-1},x_k^i)\leq\ell_{g_k}(\gamma_k).
\end{equation}
By \eqref{inf.osc.d} we see
$
m_k\leq \frac{\ell_{g_k}(\gamma_k)}{\rho} < \frac{C}{\rho}.
$
%Without loss of generality, 
Assume $m_k$ is a fixed $m$, and $\alpha_k\rightarrow \alpha_\infty$ since $\alpha_k\subset\gamma_k\subset K$ and $K$ is compact. The cardinality of $\alpha_\infty$ is still $m$ as $a\leq d_{g_0}(x^{i-1}_k,x^{i}_m)\leq 2a$. By \eqref{d=d_g.loc.K}, 
\begin{align*}
\lim_{k\rightarrow+\infty}\mathcal{L}_{d_{g_k,\Sigma}}(\alpha_k)&=\lim_{k\rightarrow+\infty}\sum_{i=1}^md_{g_k,\Sigma}(x_k^{i-1},x_k^i)
=\sum_{i=1}^md(x_\infty^{i-1},x_\infty^i)\\
&=\sum_{i=1}^md_{g,\Sigma}(x_\infty^{i-1},x_\infty^i)
\geq d_{g,\Sigma}(x_0,y_0).
\end{align*}
This together with \eqref{obvious} yield (i). 

\vspace{.1cm}

(ii) Without loss of generality, we assume $\partial U\cap\partial V=\emptyset$.

\vspace{.1cm}

Let ${\mathscr C}(U,V)$ be the set of curves $\gamma:[0,1]\rightarrow \overline{V\backslash U}$ with $\gamma(0)\in\partial U$, $\gamma_k(1)\in\partial V$, $\gamma_k((0,1))\subset V\backslash U$.
 By Lemma \ref{Sigma12} (ii),
$$
d_{g_k}(\partial U,\partial V)=\inf\left\{\ell_{g_k}(\gamma):\gamma\in {\mathscr C}(U,V)\right\}.
$$
For any $k$, we may choose $\gamma_k\in {\mathscr C}(U,V)$ with 
$$
\ell_{g_k}(\gamma_k)\leq d_{g_k}(\partial U,\partial V)+\frac{1}{k}.
$$
Let $\gamma_k(0)\rightarrow x_\infty$ and $\gamma_k(1)\rightarrow y_\infty$.  %{\color{red}Uniform boundedness of $d_{g_k}(\partial U,\partial V)$ can be seen from taking a finite covering of a fixed curve in ${\mathscr C}(U,V)$ by disks in $\Sigma$ and applying Proposition \ref{conv0} on each disk, together with the compactness of $U$ and $V$ which in particular implies $d_g(\partial U,\partial V)$ is finite. } 
We can find $\gamma\in {\mathscr C}(U,V)$ such that %$\gamma(0)=x_\infty,\gamma(1)=y_\infty$ and 
$$
\ell_{g}(\gamma)\leq d_{g}(\partial U,\partial V)+\epsilon.
$$
Then, by the trace embedding theorem, 
	$$
	\lim_{k\to\infty}d_{g_k}(\partial U,\partial V)\leq\lim_{k\to\infty}\ell_{g_k}(\gamma)\to\ell_g(\gamma)\leq d_g(\partial U,\partial V)+\epsilon.
	$$
Then, by (i) we have 
$$
d_g(\partial U,\partial V)\leq d_g(x_\infty,y_\infty)\leq\lim_{k\rightarrow+\infty}d_{g_k}(\partial U,\partial V).
$$
We complete the proof by letting $\epsilon\rightarrow 0$.
\endproof

%\begin{cor}\label{conv.domain}
%Let $g_k,g$ be as in Corollary \ref{conv.length}. Let $U,V$ be compact domains in $\Sigma$ and $U\subset V$. Assume $|\K_g|(\{x\})<\tau_1$ in a neighborhood of $\overline{V\backslash U}$.Then 
%$$d_{g_k}(\partial U,\partial V)\rightarrow d_g(\partial U,\partial V)$$ where $d_{g}(\partial U,\partial V)=\inf_{x\in\partial U,y\in\partial V}d_g(x,y).$\end{cor}

\vspace{.1cm}

We can deduce continuity of $d_{g,\Sigma}$ via Proposition \ref{conv0}, even when the point $|\K_g|$-mass is not smaller than $2\pi$, provided $d_{g,\Sigma}$ is finite on $\Sigma$. 
%(cf. Reshetnyak's condition $|\K_g|(\{x\})<2\pi,|\K_g|(\{y\})<2\pi$ at all $x,y$, see Theorem \ref{2 option}). We have 
\begin{cor}\label{continuity.distance}
Let $g=e^{2u}g_0\in\mathcal{M}(\Sigma,g_0)$.
Assume that $|\K_g|<+\infty$ and $d_{g,\Sigma}(x,y)<+\infty$
for any $x,y\in\Sigma$. Then $d_{g,\Sigma}$ is continuous in
$\Sigma\times\Sigma$. Moreover, if $S\subset\Sigma$ is a finite set, then
$
d_{g,\Sigma\setminus S}=d_{g,\Sigma}|_{\Sigma\setminus S}.
$
\end{cor}
\proof
Let $x_k\rightarrow x_0$ and $y_k\rightarrow y_0$. It suffices to prove $d_{g,\Sigma}(x_k,x_0)\rightarrow 0$, 
since
$$
|d_{g,\Sigma}(x_k,y_k)-d_{g,\Sigma}(x_0,y_0)|\leq d_{g,\Sigma}(x_k,x_0)+d_{g,\Sigma}(y_k,y_0).
$$

Since $|\K_g|<+\infty$, there are at most a finite set $A_{\frac{4}{3}\pi}$ of points on $\Sigma$ where $|\K_g|(\{x\})>\frac{4}{3}\pi$. Select an isothermal coordinate system centered around $x_0$ such that $x_0$ is the only possible one from $A_\frac{4}{3}\pi$, and assume $g=e^{2u'}g_{\rm euc}$ there.
By taking $r_0$ small we assume 
$$
|\K_g|(D_{r_0}\backslash\{0\})<\frac{4}{3}\pi.
$$
Define functions
$$
l_1(r)=d_{g,\Sigma}(\partial D_r,\partial D_{{2r}}),\s l_2(r)=\ell_g(\partial D_r).
$$
We claim that $l_1(r),l_2(r)\rightarrow 0$ as $r\to 0$.

First, since $d_{g,\Sigma}(0,y)$ is finite for a fixed $y\in\partial D_{r_0}$, we can select a curve $\gamma$ with $\gamma(0)=0$, $\gamma(1)=y$ and $\ell_g(\gamma)=\int_\gamma e^{u'}<+\infty$. Let $t_r$ be the smallest $t$ such that $\gamma(t_r)\in\partial D_{r}$. Then we have
$$
l_1(r)\leq\int_{\gamma|_{[0,t_r]}}e^{u'}.
$$
Obviously, $t_r\rightarrow 0$ as $r\rightarrow 0$, which yields that $l_1(r)\rightarrow 0$ as $r\rightarrow 0$.

Next, let $r_k\rightarrow 0$.    
Define 
$$
{u_k(x)=u'(r_kx)+\log r_k-c_k}
$$
and $g_k=e^{2u_k}g_{\rm euc}$, where $c_k$ is chosen so that the mean value of $u_k$ on $D_4\backslash D_{1/2}$ vanishes. By Lemma \ref{rq}, 
$$
\|\nabla u_k\|_{L^1(D_4\backslash D_\frac{1}{4})}
=\frac{1}{r_k}\|\nabla u'\|_{L^1(D_{4r_k}\backslash D_\frac{r_k}{4})}<C.
$$
By the Poincar\'e inequality, we may assume $u_k\to u_\infty$ weakly in $W^{1,q}(D_4\backslash D_{1/2})$ and strongly in $L^1(D_4\backslash D_{1/2})$. Set $g_\infty=e^{2u_\infty}g_{\rm euc}$.
By Corollary \ref{conv.length} (ii),
$$
\lim_{k\rightarrow+\infty}e^{-c_k}l_1(r_k)=
\lim_{k\rightarrow+\infty}d_{g_k,D_4}(\partial D_1,\partial D_2)=d_{g_\infty,D_4}(\partial D_1,\partial D_2)>0.
$$
So $c_k\rightarrow-\infty$. 
By Corollary \ref{MTI} and the trace embedding theorem,
$$
\int_{\partial D_1}e^{u_k}\to\int_{\partial D_1}e^{u_\infty}<\infty.
$$
 Then
\begin{equation}\label{length.of.circle}
l_2(r_k)=\int_{\partial D_{r_k}}e^{u'}=e^{c_k}\int_{\partial D_1}e^{u_k}\rightarrow 0.
\end{equation} 

Finally,
$$
d_{g,\Sigma}(0,x_k)\leq d_{g,\Sigma}(0,\gamma(t_{|x_k|}))+
d_{g,\Sigma}(\gamma(t_{|x_k|}),x_k)\leq \ell_{g}(\gamma|_{[0,t_{|x_k|}]})+l_2(|x_k|)\rightarrow 0,
$$
where we recall that $t_{|x_k|}$ is the first time when $\gamma$ hits $\partial D_{|x_k|}$ emanating from $0$. 

\vspace{.15cm}

Now, we prove $d_{g,\Sigma\setminus S}=d_{g,\Sigma}|_{\Sigma\setminus S}$. 
By Lemma \ref{Sigma12}, for any $x,y\in\Sigma\backslash A$
$$
d_{g,\Sigma\setminus S}(x,y)\geq d_{g,\Sigma}(x,y).
$$
On the other hand, for any $y,y'\in\Sigma\backslash S$, we select $\gamma\subset\Sigma$, such that $\gamma(0)=y$ and $\gamma(1)=y'$ and 
$\int_{\gamma}e^{u}\leq d_{g,\Sigma}(y,y')+\epsilon$. Let $S=\{p_0,...,p_m\}$. Take a disk $D_\delta$ centered at $p_0$ in the open surface $\Sigma$ (similar for other $p_i$). Assume $y,y'\not\in \overline{D}_\delta$ by choosing $\delta$ small. If $\gamma\cap\overline{D}_\delta\not=\emptyset$, 
 let $t_1$ be the first time that $\gamma$ intersects $\partial D_{\delta}$ and $t_2$ be the last time that $\gamma$ leaves $\partial D_{\delta}$. Let $\Gamma$ be a circular arc from 
$\gamma(t_1)$ to $\gamma(t_1)$ on $\partial D_\delta$, and set $\Gamma=\emptyset$ if $\gamma\cap\overline{D}_\delta=\emptyset$. 
 Then 
 \begin{align*}
 \int_\gamma e^u &\geq \big(\int_{\gamma|_{[0,t_1]}}+\int_{\gamma|_{[t_1,1]}}\big)e^u 
\geq\big(\int_{\gamma|_{[0,t_1]}}+\int_{\gamma|_{[t_1,1]}}+
\int_{\Gamma}\big)e^u-\int_{\partial D_\delta(x_0)}e^{u}\\
&\geq d_{g,\Sigma\setminus  S}(y,y')-l_2(\delta).
 \end{align*}
 Letting $\delta\to 0$, we get 
 $$
 d_{g,\Sigma\backslash S}(y,y')\leq d_{g,\Sigma}(y,y').
 $$
This completes the proof. 
\endproof

\vspace{.1cm}

The solution $u_k$ of $-\Delta u_k=\mu_k$ will only convergence after substract $c_k$ (average of $u_k$) and geometric conclusion  comes from knowing $c_k\rightarrow-\infty$ or $c_k\rightarrow c$. This is essentially H\'elein's convergence theorem \cite[Theorem 5.1.1]{H}.

\subsection{A  three-circle type lemma along cylinders}

In the development of the general theory of harmonic maps from a surface,  Parker's bubble tree construction \cite{Parker} is an important contribution after Sacks-Uhlenbeck's seminal work \cite{SU1}; more information, especially on compactness while controlling topological data such as homotopy class, can be obtained from refined analysis on the ``neck'' region connecting the regular region and the bubble region,  an integral version of the Hadamard three-circle theorem plays a key role there (cf. \cite{CT}, \cite{QT}).  Bubbling analysis of almost harmonic maps have developed by Ding-Tian \cite{DT}, Topping \cite{Topping} and others. We now demonstrate that the same idea leads to geometric estimates in our setting as the potential $u$ is related to harmonic functions.

We introduce some notations for $i\in\mathbb{Z}$:
\begin{eqnarray*}
S_i&=&S^1\times\{iL\},\\
Q_i&=&S^1\times[(i-1)L,iL],\\
L_i&=&\{0\}\times[(i-1)L,iL].
\end{eqnarray*}
Suppose $g=e^{2u}(dt^2+d\theta^2)$ is defined on a cylinder $Q=S^1\times[-L,4L]$ for $L>0$ where $-L$ is used purely for convenience of applying interior elliptic estimates and it can always be achieved by shifting a constant amount along the $t$-direction.  

%As before, we denote the length of a curve $\gamma$ in the metric $g$ by $\ell_g(\gamma)$.  
The {diameter} of a compact domain $\Omega$ for a metric $g$ is  
$$\diam(\Omega,g)=\sup_{x,y\in\overline\Omega}d_{g,\Omega}(x,y).$$  
In this section, $\K_g$ is the Gauss curvature measure for $g=e^{2u}(dt^2+d\theta^2)$.
\begin{lem}\label{3-circle}
Let $L$, $\Lambda$ and $\kappa$ be positive constants. Suppose that for a.e. $t\in [-L,3L]$, $L>1$ 
$$
\|\nabla u\|_{L^1(S^1\times[t,t+1])}<\Lambda \s \mbox{and}\s  \int_{S^1\times\{t\}}\frac{\partial u}{\partial t}<-2\pi \kappa<0.
$$
Assume $L>\frac{16\Lambda}{\kappa}$. Then there is $\tau_0=\tau_0(\kappa,\Lambda)<c_0$ so that if 
$$
|\K_g|(S^1\times[-L,4L])<\tau_0
$$ 
then 
\begin{eqnarray}\label{decay.SS}
d_{g,Q}(S_2,S_1)&<&e^{-\frac{\kappa}{2}L} d_{g,Q}(S_1,S_0),\\
\label{decay.diam/SS}
\frac{\diam(Q_1,g)}{d_g(S_0,S_1)}&<& 2\frac{e^{8\Lambda}
(1+8\Lambda)}{1-e^{-16\Lambda}},\\
\label{decay.S}
\ell_g(S_2)&<&e^{-\frac{\kappa}{2}L} \ell_g(S_1),\\
\label{decay.L}
\ell_g(L_2)&<&e^{-\frac{\kappa}{2}L}\ell_g(L_1).%\\\label{decay.A}
%Area(Q_2,g)&<&e^{-\kappa L}Area(Q_1,g).
\end{eqnarray}
\end{lem}
\proof %We only prove {\color{red} the first and the last inequality.??}
Assume \eqref{decay.SS} is not true. Then we can find $g_k=e^{2u_k}(dt^2+d\theta^2)$ with $|\K_{g_k}|(Q)\rightarrow 0$
and 
$$
d_{g_k,Q}(S_2,S_1)\geq e^{-\frac{\kappa}{2}L} d_{g_k,Q}(S_1,S_0).
$$ 
Then, by the Poincar\'e inequality and the Sobolev embedding theorem, we can choose $c_k$ such that $u_k-c_k$ converges in $L^1(Q)$ to some $v$. From Lemma \ref{W1q.convergence} {2), the convergence is   in $W^{1,q}_{\rm loc}(Q)$} and $v$ is harmonic. As $\nabla v$ is also harmonic in $Q$, by the mean value theorem there is $t$ such that 
$$
\|\nabla v\|_{C^0(S^1\times[0,3L])}\leq \frac{1}{|D_{\frac{1}{2}}|}\|\nabla v\|_{L^1(S^1\times [t,t+1])}\leq\frac{4}{\pi}\lim_{k\to+\infty}\|\nabla (u_k-c_k)\|_{L^1(S^1\times[t,t+1])}\leq \frac{4}{\pi}\Lambda.
$$ 
For $g_\infty = e^{2v} g_{\rm euc}$, by Corollary \ref{conv.length} (ii)
\begin{equation}\label{contr}
d_{g_\infty,Q}(S_2,S_1)\geq e^{-\frac{\kappa}{2}L} d_{g_\infty,Q}(S_1,S_0).
\end{equation}
On the cylinder, the harmonic function $v$ can be expanded as 
$$
v=a_0+at+\sum_{k=1}^\infty (a_k(t)\cos(k\theta)+b_k(t)\sin(k\theta)):=a_0+at+v'.
$$
%where $a_k(t)=c_ke^{kt}+\tilde{c}_ke^{-kt}$ and  $b_k(t)=c_k'e^{kt}+\tilde{c_k}'e^{-kt}$. 
Then 
$$
a=\frac{1}{2\pi}\int_{S^1\times\{t\}}\frac{\partial v}{\partial t}\s \mbox{and}\s a_0+at=\frac{1}{2\pi}\int_{S^1\times\{t\}}v.
$$
It follows 
$$
\left|v'(t,\theta)\right|=\left| v(t,\theta) - \frac{1}{2\pi}\int_{S^1\times\{t\}}v \right|\leq \pi\max_\theta\left| \frac{\partial v}{\partial\theta}(t,\theta)\right|\leq 4\Lambda, 
$$
and
$
|a|\leq \frac{4}{\pi}\Lambda.  
$
In particular,  this yields a lower bound on $a$. 
Note that
$$
a=\frac{1}{2\pi}\int_{S^1\times[0,1]}\frac{\partial v}{\partial t}d\theta dt=\lim_{k\rightarrow+\infty}
\frac{1}{2\pi}\int_{S^1\times[0,1]}\frac{\partial (u_k-c_k)}{\partial t}d\theta dt\leq-\kappa.
$$
Write $g_\infty=e^{2v}(dt^2+d\theta^2)$ and let $g_\infty'=
e^{2a_0+2at}(dt^2+d\theta^2)$, so $g_\infty=e^{2v'}g_\infty'$. Then
$$
e^{-8\Lambda}g_\infty'\leq g_\infty<e^{8\Lambda} g_\infty'.
$$
A direct computation shows 
\begin{equation}\label{direct computation}
d_{g_\infty'}(S_i,S_{i-1})=\frac{e^{aL}-1}{a}e^{a_0}e^{(i-1)aL}.
\end{equation}
In fact, 
$$
d_{g_\infty',Q}(S_i,S_{i-1})=
\inf\left\{\int_{\gamma}e^{v'}:\gamma:[0,1]\rightarrow Q_i, \gamma[0]\in S_{i-1}, 
\gamma(1)\in S_i\right\}.
$$
For any curve $\gamma(s)=(t(s),\theta(s)):[0,1]\rightarrow Q_i$, 
with $\gamma(0)\in S_{i-1}$ and $\gamma(1)\in S_i$, 
\begin{align*}
\int_\gamma e^{v'}&=\int_{(i-1)L}^{iL}e^{a_0+at}\left(\left|\frac{dt}{ds}\right|^2+\left|\frac{d\theta}{ds}\right|^2\right)^{1/2}ds\geq\int_{\{0\}\times[(i-1)L,iL]}e^{a_0+at}\left|\frac{dt}{ds}\right|ds\\
&\geq\int_{(i-1)L}^{iL}e^{a_0+at}dt=\frac{e^{aL}-1}{a}e^{a_0}e^{(i-1)aL},
\end{align*}
and the equality can be attained by the curve $\theta(s)=0, t(s)=Ls$. Then we get from \eqref{contr} 
$$
e^{-\frac{\kappa}{2}L}\leq\frac{d_{g_\infty}(S_2,S_1)}{d_{g_\infty}(S_1,S_0)}
\leq e^{8\Lambda}\frac{d_{g_\infty'}(S_2,S_1)}{d_{g_\infty'}(S_1,S_0)}= e^{8\Lambda+aL},
$$
but this contradicts the assumption $L>\frac{16\Lambda}{\kappa}$ as $a\leq -\kappa$.  Hence we have \eqref{decay.SS}.

\vspace{.1cm}

We now prove \eqref{decay.diam/SS} by contradiction. Assume there were $g_k=e^{2u_k}(dt^2+d\theta^2)$ with $|\K_{g_k}|(Q)\rightarrow 0$
and 
\begin{equation}\label{con.diam/SS}
\frac{\diam(Q_1,g_k)}{d_{g_k}(S_0,S_1)}\geq 2\frac{e^{8\Lambda}
(1+8\Lambda)}{1-e^{-16\Lambda}}.
\end{equation} 
As in the proof of \eqref{decay.SS}, 
we assume $u_k-c_k$ {converges
to a harmonic function $v$  in $W^{1,q}$}. %, where $v$ is as in the proof of \eqref{decay.SS}. 

%By the definition of diameter, t
There exist points $x_k=(\theta_k,t_k)$ and $y_k=(\theta_k',t_k')$ on $Q_1$ such that 
$$
d_{g_k,Q_1}(x_k,y_k)\geq {\diam}(Q_1,g_k)-\frac{e^{c_k}}{k}.
$$ 
Assume $\theta_k\rightarrow\theta_\infty$ and $\theta_k'\rightarrow\theta_\infty'$.
As $x_k, y_k$ are on the loop  $
\Gamma_{k}=\partial ([\theta_k,\theta_k']\times[0,L])$ (include the case $\theta_k=\theta_k'$ where $[0,L]$ is counted twice with opposite direction), then  
\begin{align*}
\lim_{k\to+\infty}e^{-c_k}&d_{g_k,Q_1}(x_k,y_k)\leq\lim_{k\to+\infty}\frac{1}{2}e^{-c_k}\ell_{g_k}(\Gamma_k)\\
&\leq\frac{1}{2}\lim_{k\to+\infty}\left(\int_{S_0\cup S_1}e^{u_k-c_k}+
\int_{\{\theta_k,\theta_k'\}\times[0,L]}e^{u_k-c_k}\right)\\
&=\frac{1}{2}\left(\int_{S_0\cup S_1}e^{v}+
\int_{\{\theta_\infty,\theta_\infty'\}\times[0,L]}e^{v}\right)\\
&\leq\frac{1}{2}e^{4\Lambda+a_0}\left(\int_{S_0\cup S_1}e^{at}+
\int_{\{\theta_\infty,\theta_\infty'\}\times[0,L]}e^{at}\right)\\
&=\frac{1}{2}e^{4\Lambda+a_0}\left((e^{aL}+1)2\pi+2\frac{e^{aL}-1}{a}\right),
\end{align*}
it follows from the choice of $x_k,y_k$ that 
$$
\lim_{k\rightarrow+\infty}e^{-c_k}{\diam}(Q_1,g_k)\leq e^{4\Lambda+a_0}\left((e^{aL}+1)\pi+\frac{e^{aL}-1}{a}\right).
$$
Since $e^{-c_k}g_k \to g_\infty$,  we have 
\begin{align*}
\lim_{k\rightarrow+\infty}\frac{\diam(Q_1,g_k)}{d_{g_k}(S_0,S_1)}&=\frac{\lim_{k\to+\infty}e^{-c_k}\diam(Q_1,g_k)}{d_{g_\infty}(S_0,S_1)}\leq e^{8\Lambda}\frac{(e^{aL}+1)\pi+\frac{e^{aL}-1}{a}}{\frac{e^{aL}-1}{a}}\\
&\leq e^{8\Lambda}\frac{2\pi|a|+1}{1-e^{a L}}\leq \frac{e^{8\Lambda}(1+8\Lambda)}{1-e^{-\kappa L}}\leq \frac{e^{8\Lambda}(1+8\Lambda)}{1-e^{16\Lambda}}.
\end{align*}
Together with \eqref{con.diam/SS}, we see a contradiction.
\eqref{decay.L} can be proved similarly.% we omit it.
\endproof

\vspace{.1cm}

Lemma \ref{3-circle} will be used in the blow-up analysis for
the convergence of distance when $\K(\{x\})<2\pi$. When $\K(\{x\})=2\pi$ and distance is finite, the lemma below can be used to deduce continuity of the distance at $x$.

\begin{lem}\label{3-circle'}
Let $g=e^{2u}(dt^2+d\theta^2)$ on $Q=S^1\times [-L,2L]$ and 
$
\|\nabla u\|_{L^1(S^1\times [-L,2L])}<\Lambda .
$
Then there is $\tau_0'=\tau_0'(\Lambda)<c_0$ so that if 
$$
|\K_g|(S^1\times[-L,2L])+\left\|\int_{S^1\times\{t\}}\frac{\partial u}{\partial t}\right\|_{L^\infty([-L,2L])} <\tau_0',
$$ 
then 
\begin{equation}\label{decay.SS/L}
e^{-8\Lambda-1}<\frac{d_{g,Q}(S_0,S_1)}{\ell_g(L_1)}<e^{8\Lambda+1}.
\end{equation} 
\end{lem}
%{\color{red} The same proof of before, omit?}

\proof 
Assume there is no such $\tau_0'$. Then we can find $g_k=e^{2u_k}(dt^2+d\theta^2)$ with
\begin{align*}
&\|\nabla u_k\|_{L^1(S^1\times[-L,2L])}<\Lambda,\\
& |\K_{g_k}|(S^1\times[-L,2L])\rightarrow 0,\\
&\left\|\int_{S^1\times\{t\}}\frac{\partial u_k}{\partial t}\right\|_{L^\infty([-L,2L])}\rightarrow 0,
\end{align*} 
but
\begin{equation}\label{either_or}
\mbox{either\s} \frac{d_{g_{k},Q}(S_0,S_1)}{\ell_{g_{ k}}(L_1)}\geq e^{8\Lambda+1}, \s
\mbox{or \s}
\frac{d_{g_{ k},Q}(S_0,S_1)}{\ell_{g_{ k}}(L_1)}\leq e^{-8\Lambda-1}.
\end{equation}

{ Choose $c_k$ as in the proof of Lemma \ref{3-circle}:
%Then, by the Poincar\'e inequality and the Sobolev embedding theorem, we may choose $c_k$, such that $u_k-c_k$ converges in $L^1(S^1\times[-1,2])$. Then it follows from Lemma \ref{W1q.convergence}, 
$u_k-c_k\to {\mbox {a harmonic $v$}}$ in $W^{1,q}_{\rm loc}(S^1\times[-L,2L])$ with 
%Since $\nabla v$ is also harmonic in $Q$ by the mean value theorem there is $t$ such that 
$$
\|\nabla v\|_{C^0(S^1\times[0,L])}\leq \frac{4}{\pi}\Lambda .%\frac{1}{|D_{{1}/{2}}|}\|\nabla  v\|_{L^1(S^1\times [t,t+1])}\leq\frac{4}{\pi}\lim_{k\to+\infty}\|\nabla (u_k-c_k)\|_{L^1(S^1\times[t,t+1])}\leq \frac{4}{\pi}\Lambda.
$$ 
%The harmonic function $v$ can be expanded as 
Expand 
$$
v=a_0+at+\sum_{k=1}^\infty (a_k(t)\cos(k\theta)+b_k(t)\sin(k\theta)):=a_0+at+v'
$$
%where $a_k(t)=c_ke^{kt}+\tilde{c}_ke^{-kt}$ and  $b_k(t)=c_k'e^{kt}+\tilde{c_k}'e^{-kt}$. 
%Then 
%$$
%a=\frac{1}{2\pi}\int_{S^1\times\{t\}}\frac{\partial v}{\partial t}\s \mbox{and}\s a_0+at=\frac{1}{2\pi}\int_{S^1\times\{t\}}v.
%$$
%It follows 
so 
$
\left|v'(t,\theta)\right|\leq 4\Lambda.%=\left| v(t,\theta) - \frac{1}{2\pi}\int_{S^1\times\{t\}}v \right|\leq \pi\max_\theta\left| \frac{\partial v}{\partial\theta}(t,\theta)\right|\leq 4\Lambda. 
$
%In particular,  this yields a lower bound on $a$. 
Note that
$$
a=\frac{1}{2\pi}\int_{S^1\times[0,L]}\frac{\partial v}{\partial t}\,d\theta dt=\lim_{k\rightarrow+\infty}
\frac{1}{2\pi}\int_{S^1\times[0,L]}\frac{\partial (u_k-c_k)}{\partial t}\,d\theta dt=0.
$$}
Let $g_\infty=e^{2v}(dt^2+d\theta^2)$ and $g_\infty'=
e^{2a_0}(dt^2+d\theta^2)$. Therefore $g_\infty=e^{2v'}g_\infty'$ and 
$$
e^{-8\Lambda}g_\infty'\leq g_\infty<e^{8\Lambda} g_\infty'.
$$
By a direct computation 
$$
e^{a_0}e^{-4\Lambda}L\leq\ell_{g_\infty}(L_1)= \int^{L}_0e^{a_0+v'(0,t)}dt\leq e^{a_0}e^{4\Lambda}L
$$ 
and (cf. \eqref{direct computation} with $a=0$)
$$
e^{a_0}e^{-4\Lambda}L=e^{-4\Lambda}d_{g_\infty'}(S_0,S_1)\leq d_{g_\infty}(S_0,S_1)\leq e^{4\Lambda}d_{g_\infty'}(S_0,S_1)=e^{a_0}e^{4\Lambda}L.
$$
Then 
$$
e^{-8\Lambda}\leq\frac{d_{g_\infty}(S_0,S_1)}{\ell_{g_\infty}(L_1)}
\leq  e^{8\Lambda}.
$$
In light of \eqref{either_or}, we then see 
$$
e^{8\Lambda+1}\leq e^{8\Lambda},\s \mbox{or}\s 
e^{-8\Lambda-1}\geq e^{-8\Lambda}.
$$
But this is clearly impossible. 
\endproof

\subsection{Finiteness of distance and curvature measure $\leq 2\pi$}

The first two items in the theorem below are observed in \cite[Theorem 3.1]{R1} (see \cite[Proposition 5.3]{Troyanov2}). %We will include a proof of it for convenience for proving the other two items where we need to apply the three-circle lemma. 

\begin{thm}\label{2 option} 
Let $(\Sigma,g_0)$ be a Riemannian surface and $g\in\M(\Sigma,g_0)$ with $|\K_g|(\Sigma)<+\infty$. 
Then 
\begin{enumerate}
\item $d_{g,\Sigma}$ is a distance function over ${\Sigma'}=\{x\in\Sigma:\K_g(\{x\})<2\pi\}.$
\item If $\K_g(\{x\})>2\pi$, then for any $\delta>0$ 
$$
\lim_{r\rightarrow 0}d_{g}(\partial B_{\delta}(x,g_0),\partial B_{r}(x,g_0))=+\infty.
$$
Consequently, if $d_{g,\Sigma'}<C$ on $U_x\backslash\{x\}\times U_x\backslash\{x\}$ for some constant $C>0$ and neighborhood $U_x$ of $x$ in $\Sigma$, then $\K_g(\{x\})\leq 2\pi$.
\item Let $\Sigma=\Sigma'\cup\{p_1,...,p_n\}$.  If there is $\delta$ so that $d_{g}(\partial B_\delta(p_i,g_0),\partial B_r(p_i,g_0))<C$ for any $r<\delta$, then $d_{g,\Sigma}$ is continuous on $\Sigma\times\Sigma$.   
Consequently, if  $d_{g,\Sigma'}<C$ on $U_{p_i}\backslash\{p_i\}\times U_{p_i}\backslash\{p_i\}$, then $d_{g,\Sigma'}$ continuously extends across $\{p_1,...,p_n\}$ to a distance function $d_{g,\Sigma}$ on $\Sigma\times\Sigma$.
\end{enumerate}
\end{thm}

\proof {By Corollary \ref{MTI} 2), $e^{u}\in W^{1,1}_{\rm loc}(\Sigma')$.
Then for any $x, y\in\Sigma'$, and a smooth curve $\gamma$ from
$x$ to $y$ in $\Sigma'$, it follows from the trace embedding theorem that
$
\int_\gamma e^{u}<+\infty.
$ }
Then $d_{g,\Sigma}(x,y)<+\infty$ for any $x,y\in\Sigma'$ 
%(similarly we can show $Area(D,g)<+\infty$). 
So (1) is established. 

\vspace{.1cm}

For (2),  without loss of generality, we work on $D\backslash\{0\}$.  If $\K(\{0\})>2\pi$,  take
$$
\kappa=\frac{\K_{g}(\{0\})-2\pi}{2}>0.
$$
As $|\K_g|(D\backslash\{0\})<+\infty$, we can take a small $r_0$ so that $|\K_g|(D_{r_0}\backslash\{0\})<\kappa$.  By \eqref{a.e.} and \eqref{GB}, for a.e. $\delta\in(0,r_0)$ we have 
$$
-\int_{\partial D_\delta}\frac{\partial u}{\partial r}=\K_g(D_\delta)=\K(D_\delta\backslash\{0\})+
\K_g(\{0\})>-\kappa+2\kappa+2\pi=2\pi+ \kappa.
$$ 
Change coordinates via
$re^{i\theta}\rightarrow (\theta,t)=(\theta,\log r)$
and view $g$ as a metric on $S^1\times(-\infty,0)$ by writing $g=e^{2v}(d\theta^2+dt^2)$ for $v(\theta,t)={u}(e^{t+i\theta})+t$. 
Then
$$
\int_{S^1\times\{t\}}\frac{\partial v}{\partial t}=\int_{\partial D_{e^{t}}}\frac{\partial u}{\partial r}+2\pi<-2\pi \kappa
$$ 
and 
$$
d_g(\partial D_r,\partial D_{r'})=d_g(S^1\times\{\log r\},S^1\times\{\log r'\}).
$$
By Lemma \ref{rq}, for any $t<0$ it holds 
\begin{align*}
\|\nabla v\|_{L^1(S^1\times [t,t+L])}&\leq e^{-(t+L)}\|\nabla u\|_{L^1(D_{e^{t+L}}\backslash D_{e^{t}})}+2\pi\\
&\leq C\left(e^{t+L}\|\nabla u\|_{L^1(D)}+|\K_g|(D)\right)+2\pi:=\Lambda',\\
|\K_g|(S^1\times(-\infty,0))&=|\K_g|(D\backslash\{0\}).
\end{align*}
By Lemma \ref{3-circle} (applied to $v$, $\Lambda'$) for large $L$, we get
\begin{align*}
d_g(S^1\times\{-2L\},S^1\times\{-L\})&\leq
e^{-\frac{\kappa}{2}L}d_g(S^1\times\{-3L\},S^1\times\{-2L\})\\
&\leq e^{-\frac{\kappa}{2}2L}d_g(S^1\times\{-4L\},S^1\times\{-3L\})\\
&\leq\cdots\leq e^{-\frac{\kappa}{2}(i-1)L}d_g(S^1\times\{-iL\},S^1\times\{-(i-1)L\}).
\end{align*}
Then
\begin{align*}
d_g(S^1\times\{-iL\},& \,S^1\times\{-L\})\geq  d_g(S^1\times\{-iL\},S^1\times\{-(i-1)L\})\\
&\geq e^{\frac{(i-1)\kappa}{2} L} d_g(S^1\times\{-2L\},S^1\times\{-L\})\to +\infty, \, \mbox{as $i\to+\infty$}. 
\end{align*}

Next, we  show (3). By (1) and (2), we may assume $\K_{g}(\{p_i\})=2\pi$.
We continue to use $D$ for an isothermal chart of a small disk around $p_i$ in $(\Sigma,g_0)$.
Since $|\K_{g}|(\Sigma)<+\infty$, 
$$
\lim_{r\rightarrow 0}|\K_{g}|(D_r\backslash\{0\})=
|\K_{g}|(\emptyset)=0,
$$
by rescaling, we may assume $|\K_{g}|(D\backslash\{0\})<\tau_0'$
in Lemma \ref{3-circle'}.  

Because 
\begin{align*}
\left|\int_{\partial S^1\times\{t\}} \frac{\partial u}{\partial t} d\theta\right|&=
\left|\int_{\partial S^1\times\{t\}} \frac{\partial (u+t)}{\partial t} d\theta-2\pi\right|=\left|-\int_{\partial D_{e^{-t}}}\frac{\partial u}{\partial r} rd\theta-2\pi\right|\\&=\left|\K_g(D_{e^{-t}})-\K_{g}(\{0\})\right|=\left|\K_{g}(D_{e^{-t}}\backslash\{0\})\right|\leq\tau_0',
\end{align*}
by Lemma \ref{3-circle'},
$$
\int_{e^{-k-1}}^{e^{-k}}e^{u(x,0)}dx\leq
Cd_{g,D}(\partial D_{e^{-k}},\partial D_{e^{-k-1}}).
$$
We define 
$$
a=\limsup_{r\rightarrow 0}d_{g,D}(\partial D_r,\partial D_\frac{1}{e})\in(0,+\infty).
$$ 
It is easy to check 
$$
\sum_{k=1}^\infty d_{g,D}(\partial D_{e^{-k}},\partial D_{e^{-k-1}})\leq a.
$$
Then
$$
\int_{0}^{e^{-1}}e^{u(x,0)}dx=\sum_{k=1}^\infty
\int_{e^{-k-1}}^{e^{-k}}e^{u(x,0)}dx\leq
C\sum_{k=1}^\infty d_{g,D}(\partial D_{e^{-k}},\partial D_{e^{-k-1}})<+\infty.
$$
So $d_{g,D}(0,(\frac{1}{e},0))<\infty$ and $d_{g,D}(x,0)<\infty$ for any $x$. By Corollary \ref{continuity.distance}, $d_{g,D}\in C^0(D\times D)$.% is continuous on $D\times D$.
\endproof

%\begin{cor}\label{finite}
%Let $g=e^{2u}g_{euc}\in\mathcal{M}(D)$ and {\color{blue}$\K_g(\{x\})<2\pi$ for any $x\in D\backslash\{0\}$} and $|\K_g|(D\backslash\{0\})<\tau_2'$.  If $d_{g,D\backslash\{0\}}<C$ on $D\backslash\{0\}\times D\backslash\{0\}$ for some $C>0$, then $d_{g,D}$ is well-defined and continuous on $D\times D$.
%\end{cor}

\vspace{.1cm}

For later discussion, we state a corollary of Lemma \ref{3-circle} in terms of $(r,\theta)$ instead of $(t,\theta)$. 

\begin{cor}\label{annulus}
Let $g=e^{2u}g_{\rm euc}\in\M(\C)$ with 
$\K_{g}^+(\C)<2\pi$ and $|\K_g|(\C)<+\infty$. Then  $$d_{g}(\partial D_{e^{mL}},\partial D_{e^{(m+1)L}})\rightarrow+\infty\s\mbox{ as $m\rightarrow+\infty$}.$$
\end{cor}

\proof
Change coordinates:
$
re^{i\theta}\rightarrow (\theta,t)=(\theta,-\log r),
$
and view $g$ as a metric on $S^1\times(-\infty,+\infty)$ and set $g=e^{2v}(d\theta^2+dt^2)$. Then $v(\theta,t)=u(e^{-t+\sqrt{-1}\theta})-t$. There is $t_0\in\R$ such that
$$
|\K_g|(S^1\times(-\infty,t_0])<\tau_0,\s \K_g^+(\C)<2\pi\kappa, \s\kappa<1.
$$

As in the proof of the above theorem, we have
$\int_{S^1\times\{t\}}\frac{\partial v}{\partial t}<-2\pi\kappa$. By \eqref{decay.SS} 
\[
d_{g}(S^1\times\{(i+1)L\},S^1\times\{(i+2)L\})
<e^{-\frac{\kappa}{2}L} d_{g}(S^1\times\{iL\},S^1\times\{(i+1)L\}),
\]
where $i\in\mathbb{Z}$ and $iL<t_0-L$.
Then for $m_0<m\in \mathbb{Z}^+$,
\begin{align*}
d_{g}&(S^1\times\{-m_0L\},S^1\times\{(1-m_0+1)L\})<e^{-\frac{\kappa}{2}L} d_{g}(S^1\times\{(-m_0-1)L\},S^1\times\{-m_0L\})\\
&<e^{-2\frac{\kappa}{2}L} d_{g}(S^1\times\{(-m_0-2)L\},S^1\times\{(-m_0-1)L\})\\
&<e^{-(m-m_0)\frac{\kappa}{2}L} d_{g}(S^1\times\{-mL\},S^1\times\{(-m+1)L\}).
\end{align*}
Then
%$$d_{g}(S^1\times\{-mL\},S^1\times\{(-m+1)L\})\geqe^{(m-m_0)\frac{\kappa}{2}L} d_{g}(S^1\times\{(-m_0)L\},S^1\times\{(-m_0+1)L\}).$$Since 
\begin{align*}
d_{g}(\partial D_{e^{mL}},\partial D_{e^{(m-1)L}})&=
d_{g}(S^1\times\{-mL\},S^1\times\{(-m+1)L\})\\
&\geq
e^{(m-m_0)\frac{\kappa}{2}L} d_{g}(S^1\times\{(-m_0)L\},S^1\times\{(-m_0+1)L\}),
\end{align*}
which goes to $+\infty$ as $m\to+\infty$.
\endproof

\subsection{Triviality of bubbles at a point where curvature $<2\pi$}
%Assume $\cap_{\epsilon>0}A_\epsilon =\{0\}$ or $\emptyset$ in $D$. 

We show that only ghost bubbles can develop in the sense of \eqref{diam0} below. We now describe how this will be used. Given a sequence of (singular) metrics of bounded integral curvature, their curvature measures may concentrate at a discrete set of points, namely, no matter how small the radius is $\K_{g_k}$  evaluated on the disk centered near these points  is not tending to 0 as $k\to+\infty$. Rescaling these disks develops the so-called bubbles in a limiting procedure while one may need to do this more than once (scale some part of the already scaled regions further) resulting in a bubble tree (cf. \cite{CL1},  \cite{CL2}).

%{\color{red}Given a measure sequence $\nu_k$ define on a surface $(\Sigma,g_0)$, we define
%$$
%\mathcal{C}(p,\nu_k)=\lim_{r\rightarrow 0}\varlimsup_{k\rightarrow+\infty}\nu_k(B_r^{g_0}(p)).
%$$}

\begin{pro}\label{ghost}
Assume that $g_k=e^{2u_k}g_{\rm euc}\in \M(D)$, $d_{g_k,D}$ is finite on $D$ and  
$\K_{g_k}^+$, $\K_{g_k}^-$  converges weakly to  Radon measures  $\mu^1$ and $\mu^2$ respectively, such that
\begin{enumerate}
\item $\mu^1(\{0\})<2\pi$, $\mu^2(\{0\})<\beta$,
\item $(\mu^1+\mu^2)(D\backslash\{0\})<\frac{\tau_0}{4}$,
\item $d_{g_k}(x,\partial D_{1/2}(x))\leq l_0$   for any $x\in \overline{D}_\frac{1}{2}$,
\item  $\|\nabla u_k\|_{L^1(D)}<A$
\end{enumerate}
where $\beta, l_0, A$ are positive constants and $\tau_0$ is as in Lemma \ref{3-circle}.  Then there is a subsequence $g_{k_i}$ such that
\begin{equation}\label{diam0}
\lim_{r\rightarrow 0}\lim_{i\rightarrow +\infty}\diam(D_r,g_{k_i})=0. 
\end{equation}
%\begin{equation}\label{area0}
%	\lim_{r\rightarrow 0}\lim_{i\rightarrow +\infty}Area(D_r,g_{k_i})=0. 
%\end{equation}
\end{pro}

\proof 
%It suffices to prove that {\color{red}after passing to a subsequence of $g_k$ } OR {\color{blue} 
%there is a subsequence for every subseqeunce of $g_k$} such that 
%\begin{equation}\label{diam0'}
%\lim_{r\rightarrow 0}\lim_{k\rightarrow +\infty}diam(D_r,g_k)=0.
%\end{equation}
 %(If  \eqref{diam0} is not true,
%then we can find $r_{i}\rightarrow 0$ such that 
%$\varlimsup_{k\rightarrow +\infty} %diam(D_{r_i},g_k)>\epsilon>0$.
%Then {\color{red}for each fixed $i$} we can find {\color{red} $g_{_{k^i_{j}}}$ such that $k^i_{j}\rightarrow+\infty$ and $diam(D_{r_i},g_{_{k^i_{j}}})>\epsilon$ as $j\to+\infty$. However, $g_{_{k^i_{j}}}$ satisfies (1)-(4), then it has a subsequence $g_{_{k^i_{{j_s}}}}$ satisfying 
%$$
%\lim_{r\rightarrow 0}\lim_{l\rightarrow %+\infty}diam(D_r,g_{_{k^i_{j_s}}})=0, 
%$$
%which contradicts  $diam(D_{r_i},g_{_{k^i_j}})>\epsilon$.})

For simplicity, we set $\mu=\mu^1-\mu^2$ and $\nu=\mu^1+\mu^2$. 
%Throughout the remainder of the proof, by ``a sequence converges", we will mean it  converges after passing to a subsequence.

\vspace{.1cm}

\noindent{\bf Step 1.} We prove \eqref{diam0} when $|\K_{g_k}|(D_{{1}/{2}})\leq \tau_0$ for all large $k$.  

Let $c_k$ be the mean value of $u_k$ on $D_{{1}/{2}}$. By the Poincar\'e and Sobolev inequalities, we may assume $\|u_k-c_k\|_{L^q(D_{1/2})}<C$ for any $q\in[1,2)$. By (4) and Lemma \ref{rq}, $\|\nabla u_k\|_{L^q(D_{1/2})}$ is bounded uniformly in $k$ for any $q\in[1,2)$. Then we may
find a subsequence of $u_k-c_k$, which  we still denote by $u_k-c_k$, converges to $u$ weakly in $W^{1,q}(D_{1/2})$,  and set $g=e^{2u}g_{\rm euc}$.

Since $\tau_0<c_0\leq\frac{4}{3}\pi$, by (3) and Corollary \ref{conv.length}, 
$$
e^{-c_k}l_0\geq e^{-c_k}d_{g_k,D}(0,\partial D_\frac{1}{2})\geq 
e^{-c_k}d_{g_k}(\partial D_\frac{1}{8},\partial D_\frac{1}{4})\rightarrow d_{g}(\partial D_\frac{1}{8},\partial D_\frac{1}{4})>0.
$$
Then $c_k<C$ for some constant $C$. In addition, for any $\varphi\in{C^\infty_0}(D_{1/2})$, we have
%\begin{eqnarray*}
\begin{align*}
\int_{D_\frac{1}{2}}\nabla\varphi\nabla u&=\lim_{k\rightarrow+\infty}\int_{D_\frac{1}{2}}\nabla\varphi\nabla (u_k-c_k)
=\lim_{k\rightarrow+\infty}\int_{D_\frac{1}{2}}\nabla\varphi\nabla u_k\\
&=\lim_{k\rightarrow+\infty}\int_{D_\frac{1}{2}}\varphi d\K_{g_k}
=\int_{D_\frac{1}{2}}\varphi d\mu.
%\end{eqnarray*}
\end{align*}
By the definition of $\K_g$, we have $\K_g=\mu$ on $D_{1/2}$. Corollary \ref{continuity.distance} asserts that $d_{g,D}$ is continuous on $D$, hence
$\sup_{x\in D_r}d_{g,D}(x,\partial D_r)\rightarrow 0$ as $r\rightarrow 0$.
%Then 
%\begin{align*}
%\lim_{r\rightarrow 0}\lim_{k\rightarrow +\infty}diam(D_r,d_{g_k,D})\leq
%C\lim_{r\rightarrow 0}\lim_{k\rightarrow +\infty}diam(D_r,e^{-2c_k}d_{g_k,D})
%=C\lim_{r\rightarrow 0}diam(D_r,d_{g,D})=0.
%\end{align*}

By Proposition \ref{conv0}, $d_{e^{-2c_k}g_k,D}\to d_{g,D}$ on $D_{r_0/4}$. When $r<\frac{r_0}{4}$, 
\begin{align*}
\sup_{x\in D_r}d_{g_k}(x,\partial D_r)&=e^{c_k}\sup_{x\in D_r}d_{e^{-2c_k}g_k,D}(x,\partial D_r)\\
&\leq C\sup_{x\in D_r}d_{e^{-2c_k}g_k,D}(x,\partial D_r)\rightarrow C\sup_{x\in D_r}d_{g,D}(x,\partial D_r)
\end{align*}
as $k\to\infty$. 
By \eqref{length.of.circle} and \eqref{diam.estimate}, we conclude
$$
\lim_{r\rightarrow 0}\lim_{k\rightarrow +\infty}\diam(D_r,g_k)=0.
$$
%{\color{red}This discussion holds for $|\K_{g_k}|(D_{\frac{1}{2}})<\tau_1$. }

\noindent{\bf Step 2.} We prove \eqref{diam0} when $|\K_{g_k}|(D_{{1}/{2}})>\tau_0$ holds for a subsequence (still write $g_k$). Set 
%infinitely many $k$.  Choose a subsequence (still write $g_k$) such that $|\K_{g_k}|(D_{1/2})>\tau_0$ and set
\begin{equation}\label{m}
m=\left\lceil \frac{4\nu(\{0\})}{\tau_0}\right\rceil
\end{equation}
 i.e. the least integer upper bound; since $\nu(\{0\})<\beta$, we see $m<+\infty$.
 
 We argue by induction on $m$. 
When $m\leq 2$, we have $\nu(\{0\})\leq m\frac{\tau_0}{4}\leq
\frac{3}{4}\tau_0$, so $\nu(\overline{D}_{{1}/{2}})<\tau_0$ by using (2). Then $|\K_{g_k}|(D_{{1}/{2}})<\tau_0$ for large $k$,  so \eqref{diam0} follows from Step 1  for this case.  

Now assume that $m>2$ and \eqref{diam0} holds for any $i<m$. To prove $\eqref{diam0}$ for $m$, we set 
\begin{equation}\label{r_k}
r_k=\inf\left\{r:\exists x\in \overline{D_\frac{1}{2}},\s s.t.\s|\K_{g_k}|(D_r(x))\geq\frac{\tau_0}{4}\right\}.
\end{equation}
Since $$(m-1)\frac{\tau_0}{4}<\nu(\{0\})\leq m\frac{\tau_0}{4}$$ for any  small $\epsilon<\frac{1}{8}$ and sufficiently large $k$,
% such that $|\mu|(\overline{D_\epsilon})\geq (m-1)\frac{\tau_2}{4}$. Then 
we have 
\begin{equation}\label{rk.lambdak.to.0}
|\K_{g_k}|(D_{2\epsilon})\geq |\K_{g_k}|(\overline{D}_\epsilon)\geq  \nu(\{0\})>(m-1)\frac{\tau_0}{4}.  
\end{equation}
Thus, the set in \eqref{r_k} is nonempty because $x=0$ is there for $r=2\epsilon$; so $r_k\leq 2\epsilon$  $\rightarrow 0$ when letting $\epsilon\to 0$. Let $r_i^k\rightarrow r_k$, $x_i^k\rightarrow x_k$ and $|\K_{g_k}|(D_{r_i^k}(x_i^k))\geq\frac{\tau_0}{4}$. 
(As $|\K_{g_k}|(D_t(x))$ may not be continuous in $(x,t)$, 
we cannot say $|\K_{g_k}|(D_{r_k}(x_k))\geq\frac{\tau_0}{4}$). Since 
$D_{r_i^k}(x_i^k)\subset D_{r_k+\epsilon}(x_k)$ for a fixed $\epsilon$ and large $i$, we get
$$
|\K_{g_k}|(D_{r_k+\epsilon}(x_k))\geq\frac{\tau_0}{4},
$$ 
which implies that
$$
|\K_{g_k}|(\cap_{\epsilon>0}D_{r_k+\epsilon}(x_k))\geq\frac{\tau_0}{4}.
$$
Consequently we see bubble developing near $x_k$: 
\begin{equation}\label{measure.bubble}
|\K_{g_k}|(\overline{D_{r_k}(x_k)})\geq\frac{\tau_0}{4}.
\end{equation}
Here we define $D_0(x_k)=\{x_k\}$ (we allow $r_k=0$).  
Denote
$$
\lambda_k=\inf\left\{r:|\K_{g_k}|(D_r(x_k))\geq (m-1)\frac{\tau_0}{4}\right\}.
$$ 
By \eqref{rk.lambdak.to.0}, $\lambda_k\rightarrow 0$.
If $\lambda_k>0$, we define
$t_k=\lambda_k$.  If $\lambda_k=0$, then for any $\epsilon>0$, $|\K_{g_k}|(D_\epsilon)(x_k)\geq (m-1)\frac{\tau_0}{4}$. 
Since $d_{g_k,D}$ is finite on $D$, it is continuous by Corollary \ref{continuity.distance}.  Therefore, by \eqref{length.of.circle}, we can
find $t_k$ such that $t_k<\frac{1}{k^2}$ and
$\diam(D_{kt_k}(x_k),g_k)<\frac{1}{k}$.
For both cases,  the choice of $t_k$ ensures 
$$
|\K_{g_k}|(D_\lambda{(x_k)})\geq (m-1)\frac{\tau_0}{4},\s\forall \lambda>t_k.
$$
By \eqref{m}, $\nu(\{0\})\leq m\frac{\tau_0}{4}$.
Then it follows from (2) that when $k$ is sufficiently large
\begin{equation}\label{neck.region.measure}
|\K_{g_k}|(D_\frac{1}{2}\backslash D_{2t_k}(x_k))<
\frac{\tau_0}{4}+m\frac{\tau_0}{4}-(m-1)\frac{\tau_0}{4}
=\frac{\tau_0}{2}.
\end{equation}

\vspace{.1cm}

We divide $D_r(x_k)$ into the ``neck" region
$D_r(x_k)\setminus D_{{t_k/r}}(x_k)$ and the 
``bubble" region $D_{\frac{t_k}{r}}(x_k)$.

\vspace{.1cm}

\noindent {\it Step A}. {We consider the ``neck" region.} Change coordinates
$
x_k+re^{i\theta}\rightarrow (\theta,t)=(\theta,-\log r)
$
and view $g_k$ as a metric on $S^1\times(\log 2,+\infty)$, and set $g_k=e^{2v_k}(d\theta^2+dt^2)$.  As in the proof of Theorem \ref{2 option} we have
$\int_{S^1\times\{t\}}\frac{\partial v_{k}}{\partial t}<-2\pi\kappa$ for all $t>\log 2$,  and
$$
\diam(D_r\backslash D_{{t_k/r}}(x_k),g_k)=\diam(S^1\times[-\log r,\log r-\log t_k], g_k).
$$
For the fixed $r$, there is an integer $m_0\geq 0$ so that $-\log r\in[m_0L,(m_0+1)L]$. 
Since $t_k\to 0$ we can choose $m_k\to\infty$ so that $\log r-\log t_k\in[m_kL,(m_k+1)L]$.  Using \eqref{decay.diam/SS} and \eqref{decay.SS}, on $S^1\times (L,(m_k+1)L)$ we have
\begin{align*}
&\diam(S^1\times[-\log r,  \log r-\log t_k],  g_k) \leq\sum_{i=m_0}^{m_k}\diam(S^1\times [iL,(i+1)L],g_k)\\
 &\leq C\sum_{i=m_0}^{m_k}d_{g_k}(S^1\times \{iL\},S^1\times\{(i+1)L\})\leq C\sum^{m_k}_{i=m_0}e^{-(i-1)\frac{\kappa}{2}}d_{g_k}(S^1\times\{L\},S^1\times\{2L\})\\
% &\leq C\sum_{i=m_0}^{m_k} e^{{-(i-1)L}\frac{\kappa}{2}}\,d_{g_k}(S^1\times\{{m_0}L\},S^1\times\{({m_0+1})L\})\\
 % &\leq{Cd_{g_k}(S^1\times\{m_0L\},S^1\times\{(m_0+1)L\})}\\
 &\leq Ce^{-m_0{L}\frac{\kappa}{2}}d_{g_k}(S^1\times\{L\},S^1\times\{2L\})\leq C e^{-m_0{L}\frac{\kappa}{2}}d_{g_k}(x_k,\partial D_{\frac{1}{2}}(x_k))\leq Ce^{-m_0{L}\frac{\kappa}{2}} l_0,
\end{align*}
where $l_0$ is given in (3).  
Noting $m_0\rightarrow+\infty$ as $r\to 0$,  we see 
\begin{equation}\label{neck}
\lim_{r\rightarrow 0}\lim_{k\rightarrow+\infty}\diam(D_r\backslash D_{\frac{t_k}{r}}(x_k),g_k)=0.
\end{equation}

\noindent {\it Step B}. {We consider the bubble region and show} 
\begin{equation}\label{bubble.region}
\lim_{r\rightarrow 0}\lim_{k\rightarrow+\infty}\diam(D_{\frac{t_k}{r}}(x_k),g_k)=0.
\end{equation}

If $\lambda_k=0$, we get
\begin{align*}
\lim_{r\rightarrow 0}\lim_{k\rightarrow\infty}&\diam\left(D_r ,g_k\right)\leq \lim_{r\rightarrow 0}\lim_{k\rightarrow\infty}\left(\diam(D_r\backslash D_{\frac{t_k}{r}}(x_k),g_k)+\diam(D_{ {\frac{t_k}{r}}}(x_k),g_k)\right)\\
&\leq \lim_{r\rightarrow 0}\lim_{k\rightarrow\infty}\left(\diam(D_r\backslash D_{\frac{t_k}{r}}(x_k),g_k)+\diam(D_{kt_k}(x_k),g_k)\right)\\
&\leq \lim_{r\rightarrow 0}\lim_{k\rightarrow\infty}\left(\diam(D_r\backslash D_{\frac{t_k}{r}}(x_k),g_k)+\frac{1}{k}\right)
=0.
\end{align*}
Now, we assume $\lambda_k>0$. Then $t_k=\lambda_k$.
We will use $t_k$ and $x_k$ to rescale $u_k$.
Let 
\begin{equation}\label{g'_k}
u_k'(x)=u_k(x_k+t_kx)-\log t_k \s \mbox{and} \s g_k'=e^{2u_k'}g_{\rm euc}.
\end{equation}
For any $\gamma$, we have
$$
\int_\gamma e^{u_k}=\int_{\frac{\gamma-x_k}{t_k}}e^{u_k'},
$$
then
\begin{equation}\label{bubble.metric}
d_{g'_k,(D_{t}-x_k)/t_k}(y,y')=d_{g_k,D_{t}}(t_ky+x_k,t_ky'+x_k),\s \forall t,
\end{equation}
hence \eqref{bubble.region} is equivalent to
\begin{equation}\label{bubble.region.rescaled}
\lim_{r\rightarrow 0}\lim_{k\rightarrow+\infty}\diam(D_{\frac{1}{r}},g_k')=0.
\end{equation}

Let $\mu'^1$ and $\mu'^2$ be the weak limits of  $\K_{g'_k}^+$ and $\K^-_{g'_k}$ respectively and $\nu'=\mu'^1+\mu'^2$. 
Define
$$
\S=\left\{x\in\R^2:\nu'(\{x\})>\frac{\tau_0}{4}\right\}.
$$	
To prove \eqref{bubble.region.rescaled},  it suffices to verify that 
\begin{align}\label{bubble.region.reg}
\lim_{r\to 0}	\lim_{k\to\infty}\diam( D_{\frac{1}{r}}\backslash\cup_{x\in\S} D_r(x),g_k')=0
\end{align}
and for any $x\in\S$
\begin{equation}\label{bubble.in.bubble}
\lim_{r\rightarrow 0}\lim_{k\rightarrow+\infty}\diam(D_r(x),g'_k)=0.
\end{equation}
%(The latter one will be proved by induction.)

 \vspace{.15cm}
 
We now prove \eqref{bubble.region.reg}. Let $c'_k$ be the mean value of $u_k'$ on $D$ and set
$$
\hat g_k=e^{2(u_k'-c'_k)}g_{\rm euc}.
$$ 
{\it Claim.} $c'_k\to -\infty$.  By Lemma \ref{rq}, 
$$
\|\nabla u'_k\|_{L^1(D_R)} = \frac{1}{r_k}\|\nabla u\|_{L^1(D_{r_kR})}\leq \frac{1}{r_k}C\left(r^2_kA+r_k\nu(D)\right)<C<+\infty. 
$$
By taking a subsequence  for the second time,
we may assume $u'_k-c'_k\rightarrow u'$  in $L^{1}_{\rm loc}(\mathbb R^2)$.  Set $g'=e^{2u'}g_{\rm euc}$. 
%It is easy to check that $\nu'=\K_{g'}$.
For any bounded Borel set $E\subset\R^2$ and
fixed $\delta>0$, $t_kE+x_k\subset \overline{D_{\delta/2}}$ for large $k$.  Then
$$
|\K_{\hat g_k}|(E)=|\K_{g_k'}|(E)=|\K_{g_k}|(x_k+t_kE)\leq |\K_{g_k}|(\overline{D_\frac{\delta}{2}})|
$$
when $k$ is sufficiently large. 
By \cite[Theorem 1.40 (ii)]{E-G}
$$
\varlimsup_{k\rightarrow+\infty}|\K_{\hat g_k}|(E)\leq\nu(\overline{D_\frac{\delta}{2}})< \nu(D_\delta).
$$
Since $\nu(D_\delta)\to \nu(\{0\})$ as $\delta\to 0$, we have
\begin{equation}\label{bubble.measure}
\varlimsup_{k\rightarrow+\infty}|\K_{\hat g_k}|(E)\leq \nu(\{0\})\leq m\frac{\tau_0}{4}.
\end{equation}
Similarly, we have
\begin{equation}\label{bubble.measure2}
\varlimsup_{k\rightarrow+\infty}\K_{\hat g_k}(E)\leq\varlimsup_{k\rightarrow+\infty}\K_{\hat g_k}^+(E)\leq \mu^1(\{0\})<2\pi,
\end{equation}
which implies as in proof of {Theorem \ref{2 option}}  that
$$
\int_{\partial D_r}\frac{\partial u'}{\partial r}<2\pi\kappa'
$$
for some $\kappa'<1$.  By Corollary \ref{annulus} and Corollary \ref{conv.length} (ii), we have 
\begin{equation}\label{c_k infinity}
\lim_{m\rightarrow+\infty}\lim_{k\rightarrow+\infty}d_{\hat g_k}(\partial D_{e^{mL}},\partial D_{e^{(m-1)L}})=\lim_{m\rightarrow+\infty}d_{g'}(\partial D_{e^{mL}},\partial D_{e^{(m-1)L}})=+\infty.
\end{equation}
Since 
$$
e^{c_k'}d_{\hat g_k}(\partial D_{e^{mL}},\partial D_{e^{(m-1)L}})=
d_{g_k}(\partial D_{e^{mL}t_k}(x_k),\partial D_{e^{(m-1)L}t_k}(x_k))
$$ 
is uniformly bounded by assumption (3), in light of \eqref{c_k infinity} $c_k'\rightarrow -\infty$ as claimed.

\vspace{.15cm}

If $\S=\emptyset$, then
 $e^{u_k'}$ is bounded in $W^{1,q}(D_{1/r})$ for any
$r$ for some $q>1$. By the
trace embedding theorem, we get $
\diam(D_{1/r},\hat g_k)<C(r)$, which implies that
$$
\lim_{r\rightarrow 0}{\lim_{k\to\infty}}\diam(D_{1/r},g_k')=0.
$$
We get \eqref{diam0}. If $\S\neq\emptyset$, by Lemma \ref{Sobolev.distance} we may assume $d_{\hat{g}_k,(D-x_k)/t_k}$ converges in $C^0((D_{1/r}\backslash\cup_{x\in\S} D_r(x))\times (D_{1/r}\backslash\cup_{x\in\S} D_r(x)))$. Since $
d_{g_k',(D-x_k)/t_k}=e^{c_k'}d_{\hat{g}_k,(D-x_k)/t_k}$ and $c_k'\rightarrow-\infty$, we have 
$d_{g_k',(D-x_k)/t_k}\to 0$ uniformly on $(D_{1/r}\backslash\cup_{x\in\S} D_r(x))\times (D_{1/r}\backslash\cup_{x\in\S} D_r(x))$. Then for any fixed $r$, 
$$
{\lim_{k\to\infty}\diam(D_\frac{1}{r}\backslash\cup_{x\in\S} D_r(x),g_k')=0.}
$$
In summary, we have established \eqref{bubble.region.reg}.

\vspace{.15cm}

Next, we show \eqref{bubble.in.bubble} under the assumption $\S\neq\emptyset$ (when  $\S=\emptyset$, we have proved \eqref{diam0} already).
By induction, it suffices to check $g_k'$ (scaled from $g_k$ defined in \eqref{g'_k})
satisfies (1)-(4) and 
\begin{equation}\label{bubble.con}
\left\lceil \frac{4\nu'(\{x\})}{\tau_0}\right\rceil
\leq m-1,\s \forall x\in\S.
\end{equation}
Indeed, (1), (2) follow from \eqref{bubble.measure},
\eqref{bubble.measure2} while (3) and (4) from \eqref{bubble.metric} and Lemma \ref{rq}, respectively. 

By the definition of $\lambda_k$, for any $i\in\mathbb{Z}^+$, 
$$
|\K_{g_k}|(D_{t_k+\frac{1}{i}}(x_k))\geq (m-1)\frac{\tau_0}{4},\s |\K_{g_k}|(D_{t_k-\frac{1}{i}}(x_k))< (m-1)\frac{\tau_0}{4}.
$$ 
Noting that $\overline{D_{t_k}(x_k)}=\cap_iD_{t_k+\frac{1}{i}}$ and $D_{t_k}(x_k)=\cup_iD_{t_k-\frac{1}{i}}$, we have 
\begin{equation}\label{second.bubble}
	|\K_{g_k}|(\overline{D_{t_k}(x_k)}) \geq (m-1)\frac{\tau_0}{4},\s
	|\K_{g_k}|(D_{t_k}(x_k)) \leq (m-1)\frac{\tau_0}{4}.
\end{equation}
(This will be used to get a second bubble, other than the one at $0$, in order to keep the induction proceed.)

\vspace{.15cm}

\noindent{\it Claim.} $\frac{r_k}{t_k}\rightarrow 0$.  Let $x_0\in\S$. Then 
$\nu'(\overline{D_\delta(x_0)})>\frac{\tau_0}{4}$ for any $\delta$, so $|\K_{\hat g_k}|(\overline{D_\delta(x_0)})>\frac{\tau_0}{4}$ for large 
$k$. Then
$$
|\K_{g_k}|(D_{2t_k\delta}(t_kx_0+x_k))\geq |\K_{g_k}|(\overline{D_{t_k\delta}(t_kx_0+x_k)})=|\K_{\hat g_k}|(\overline{D_{\delta}(x_0)})>\frac{\tau_0}{4}.
$$
Then $\frac{r_k}{t_k}\leq 2\delta$. Letting $\delta\rightarrow 0$ establishes the claim.

\vspace{.15cm}

A consequence of the claim is that if $\S\not=\emptyset$ then $0\in\S$. In fact, by the claim, for any $\delta>0$,  we have $\delta t_k>r_k$ when $k$ is sufficiently large. Then
$$
|\K_{g_k'}|(\overline{D_\delta})\geq|\K_{g_k'}|({D_{\delta}})=|\K_{g_k}|(D_{\delta t_k}(x_k))
\geq\frac{\tau_0}{4}.
$$ 
By \cite[Theorem 1.40 (ii)]{E-G} again, we have
$$
\nu'(\overline{D_\delta})\geq\frac{\tau_0}{4}.
$$
Letting $\delta\rightarrow 0$, we get
$\nu'(\{0\})\geq\frac{\tau_0}{4}$, we see $0\in\S$.  
By \eqref{second.bubble} and \cite[Theorem 1.40 (ii)]{E-G}, \begin{equation}\label{bubble.D}
	\nu'(D)\leq\frac{m-1}{4}\tau_0,
\end{equation}
and by \eqref{bubble.measure} and \cite[Theorem 1.40 (ii)]{E-G},
$\nu'(U)\leq \nu(\{0\})$ holds for any bounded open set $U$, then
$$
\nu'(\C)=\lim_{i\rightarrow+\infty}\nu'(D_i)\leq\nu(\{0\})\leq m\frac{\tau_0}{4},
$$
hence
\begin{equation}\label{bubble.C-0}
	\nu'(\C\setminus\{0\})=\nu'(\C)-\nu'(\{0\})\leq\frac{m-1}{4}\tau_0.
\end{equation}

Now, we are able to prove \eqref{bubble.con}. If $x=0$, by \eqref{bubble.D} we get
$$
\nu'(\{0\})\leq \nu'(D)\leq \frac{m-1}{4}\tau_0.
$$
If $x\in \S\setminus\{0\}$, then it follows from \eqref{bubble.C-0} that
$$
\nu'(\{x\})\leq\nu'(\C\setminus\{0\})\leq(m-1)\frac{\tau_0}{4}.
$$
Using the induction hypothesis on $g_k'$ and taking a subsequence for the third time,  we get 
\eqref{bubble.in.bubble}. 
Then 
\begin{align*}
\lim_{k\rightarrow\infty}\diam(D_{\frac{1}{r}}(x_k),g_k')
&\leq
\lim_{s\rightarrow 0}\lim_{k\rightarrow\infty}\diam(D_\frac{1}{r}\backslash\cup_{x\in\S}D_s(x),g_k')\\
&+
\sum_{x\in\S}\lim_{s\rightarrow 0}\lim_{k\rightarrow\infty}\diam(D_s(x),g_k')=0.
\end{align*}
Therefore \eqref{diam0} holds for $m$. The induction is complete.  
\endproof

\begin{lem}\label{main.lemma.1}
Let $(\Sigma,g_0)$ be a closed Riemannian surface and $\{d_k\}$ be a sequence of continuous distance functions on $\Sigma$. $S$ 
is a finite subset of $\Sigma$. Assume $d_k$ converges to some $d$ in $C^0_{\rm loc}((\Sigma\backslash S)\times (\Sigma\backslash S))$ and
for any $p\in\S$
$$
\lim_{r\rightarrow 0}\lim_{k\rightarrow+\infty}\diam(B_r(p,g_0),d_k)=0.
$$
Then $d$ extends continuously to a $C^0(\Sigma\times \Sigma)$ function and $d_k$ converges to $d$ uniformly on $\Sigma\times\Sigma$.
\end{lem}

\proof
We assume $S=\{p_0\}$.
By the assumption,  for any $\epsilon>0$ there are $\delta,K_0$ so that $d_k(x',x)<\epsilon$ for any $x,x'\in B_\delta(p_0,g_0)$ and $k>K_0$. As $d_k\rightarrow d$ in $C^0_{\rm loc}(\Sigma\backslash S)$, there is $K_1>K_0$ such that
$$
\left|d_{k}(x,x')-d_{{k'}}(x,x')\right|<\epsilon,
$$
for any $x,x'\notin B_\delta(p_0,g_0)$ and $k,k'>K_1$. When $x\notin B_\delta(p_0,g_0)$ and $x'\in B_\delta(p_0,g_0)$, for any  $y\in\partial B_\delta(p_0,g_0)$ it holds  
\begin{align*}
\left|d_k(x,x')\right.&\left.  - \,d_{k'}(x,x')\right|\leq \left|d_{k}(x,x')-d_{k}(x,y)\right|\\
&+
\left|d_{k}(x,y)-d_{k'}(x,y)\right|+\left|d_{k'}(x,y)-d_{k'}(x,x')\right|
<3\epsilon.
\end{align*}
Hence,  for any $\epsilon>0$ we can find $K_1$ such that 
$$
\left|d_{k}(x,x')-d_{k'}(x,x')\right|<3\epsilon,\s \forall x,x'\in\Sigma\s and\s k,k'>K_1.
$$
So $d_k$ converges uniformly on $\Sigma\times\Sigma$ to a continuous function that equals $d$ away from $S$.  \endproof

\subsection{Global Reshetnyak's theorem}

In this section, we prove Theorem \ref{main-ChenLi}.

\begin{lem}\label{main.lemma.2}
Let $(\Sigma,g_0)$ be a closed Riemannian surface and $d$ be a continuous semi-distance function. Let $g=e^{2u}g_0\in\mathcal{M}(\Sigma,g_0)$. Let $S$ be a finite subset of $\Sigma$.  Suppose that $|\K_g|(\Sigma)<\infty$ and $d_{g,\Sigma}$ is finite on $\Sigma$. Suppose that for any $x\notin S$ there exists $r$ such that $d=d_{g,\Sigma}$ on $B_r(x,g_0)\times B_r(x,g_0)$.
Then 
$
d_{g,\Sigma}\geq d.
$
\end{lem}

\proof
By Corollary \ref{continuity.distance}, $d_{g,\Sigma\backslash S}=d_{g,\Sigma}|_{\Sigma\backslash S}$. It suffices to show $d_{g,\Sigma\backslash S}\geq d$ 
on $\Sigma\backslash S$ for $S=\{p_0\}$.

Let $x,x'\in\Sigma\backslash B_\delta(p_0,g_0)$.  Cover the compact set $\Sigma\backslash B_\delta(p_0,g_0)$ by finitely many 
balls $B_{r_i}(x_i,g_0)\subset \Sigma\backslash \{p_0\}$ so that $r_i<\frac{\delta}{8}$ and $d=d_{g,\Sigma}$ on each $B_{4r_i}(x_i,g_0)$.  Let $r=\min\{r_i\}$. Then  $d=d_{g,\Sigma}$ on  $B_{r}(y,g_0)$ for any $y\notin B_\delta(p_0,g_0)$.
Let $\gamma:[0,1]\rightarrow\Sigma\backslash \{p_0\}$ with $\gamma(0)=x$, $\gamma(1)=x'$ and $$\ell_g(\gamma)\leq d_{g,\Sigma}(x,x')+\epsilon.$$
If $\gamma\cap B_\delta(p_0,g_0)=\emptyset$, we select 
$0=t_0<t_1<\cdots<t_j=1$, such that $d_{g_0}(\gamma(t_i),\gamma(t_{i+1}))<r$. Then 
$$
\ell_g(\gamma)\geq \sum_i d_{g,\Sigma}(\gamma(t_i),\gamma(t_{i+1}))=\sum_i d(\gamma(t_i),\gamma(t_{i+1})) \geq d(x,x').
$$
If $\gamma\cap B_\delta(S,g_0)\neq\emptyset$,
we let $t'$ and $t''$  be the smallest and the greatest $t$ with $\gamma(t)\in\partial B_\delta(S,g_0)$ respectively. Then $\gamma|_{[0,t']}\cap B_\delta(p_0,g_0)=\emptyset$ and $\gamma|_{[t'',1]}\cap B_\delta(p_0,g_0)=\emptyset$. Hence
\begin{eqnarray*}
\ell_g(\gamma)&\geq& \ell_g(\gamma|_{[0,t']})+\ell_g(\gamma|_{[t'',1]})\geq 
d(x,\gamma(t'))+d(\gamma(t''),x')\\
&\geq& d(x,x')-d(\gamma(t'),\gamma(t''))\geq d(x,x')-2\max_{y\in\partial B_{\delta}(p_0,g_0)} d(0,y). 
\end{eqnarray*}
Letting $\epsilon\rightarrow 0$, then $\delta\rightarrow 0$, we complete the proof.
\endproof

%For clairty of presentation, we first establish the result up to selection of subsequences:
\begin{thm}\label{main}
Let $(\Sigma,g_0)$ be a closed Riemannian surface and $g_k=e^{2u_k}g_0\in\M(\Sigma,g_0)$.  Assume  that $\K_{g_k}$ converges weakly to a signed Radon measure $\mu$,  and $\K_{g_k}^+$ converges weakly to a Radon measure $\mu'$, and one of the following holds:
\begin{enumerate}
\item[1)]  $\diam(\Sigma,g_k)=1$ and $\mu'(\{x\})<2\pi$ for any $x$ in $\Sigma$; or
\item[2)] $d_{g_k,\Sigma}$ converges to a continuous distance  function $d$ on $\Sigma$.
\end{enumerate}
\noindent Then $\{u_k\}$
converges weakly to a function $u$ in $W^{1,q}$ for any $1\leq q<2$, $\K_g=\mu$ and $d_{g_k,\Sigma}$ converges to
$d_{g,\Sigma}$ uniformly where $g=e^{2u}g_0$.
%Moreover, for any $\Omega\subset\Sigma$, 
%$Area(\Omega,g_k)\rightarrow Area(\Omega,g).$
\end{thm}

\proof  
{\bf Step 1.} We prove that a subsequence of $u_k$ converges weakly in $W^{1,q}$.

Let $c_k$ be the mean value of $u_k$ over $\Sigma$ in $g_0$. By Lemma \ref{global.gradient.estimate} and the Poincar\'e inequality, $u_k-c_k$ is bounded in $W^{1,q}$. Then a subsequence of $u_k-c_k$, which we still denote by $u_k-c_k$, converges weakly in $W^{1,q}$ to a function $u'$. Let $g'=e^{2u'}g_0$, $\mu'=\K_{g'}$ and
$$
\S=\left\{x:\nu(\{x\})>\frac{\tau_0}{2}\right\}, 
$$
where $\nu=2\mu'-\mu$ is the limit of $|\K_{g_k}|=\K_{g_k}^++\K_{g_k}^-=2\K^+_{g_k}-\K_{g_k}$.
Then $\S$ is finite. Without loss of generality, we assume $\S$ has only one point $p$.

Fix a ball $B_r(p_0,g_0)\subset\subset \Sigma\backslash \S$ and $p_1,p_2\in B_{r}(p_0,g_0)$. 
By Proposition \ref{conv0}, we can choose $r$ such that
$d_{g_k,\Sigma}(p_1,p_2)\to d_{g',\Sigma}(p_1,p_2)$ as $k\to\infty$.
Each of 1) and 2) implies 
$$
0<a<\diam(\Sigma,g_k)<b
$$
for some constants $a,b$ independent of $k$. 
Then 
$$
e^{-c_k}b\geq d_{e^{-2c_k}g_k,\Sigma}(p_1,p_2)\rightarrow d_{g',\Sigma\backslash\S}(p_1,p_2)>0,
$$
which implies that $c_k<C$ for some $C$.

Now we claim that the sequence $c_k$ is bounded below as well. Otherwise, after passing to a subsequence, we would have $c_k\rightarrow-\infty$.
By Lemma \ref{Sobolev.distance}, $\|d_{e^{-2c_k}g_k,\Sigma}\|_{W^{1,q}(\Omega\times\Omega)}<C(\Omega)$ for  $\Omega\subset\subset\Sigma\backslash \{p_0\}$ and some $q>2$, then  $d_{e^{2u_k-2c_k}g_0,\Sigma}$ converges  uniformly on any $\Omega\subset\subset \Sigma\backslash\S$, so $d_{g_k,\Sigma}  = e^{c_k}d_{e^{2u_k-2c_k}g_0,\Sigma}\to 0$ uniformly on any $\Omega\subset\subset \Sigma\backslash\S$. Therefore, when 1) holds, by Proposition \ref{ghost}, we know $d_{g_k,\Sigma}\rightarrow 0$ uniformly on $\Sigma$, but this is impossible for $\diam(\Sigma,g_k)=1$;
when 2) holds, $d=0$,  which is impossible for a distance function.

Now, $|c_k|$ is a bounded sequence, we may assume $u_k\to u$ weakly in $W^{1,q}$.

\vspace{.1cm}

{\bf Step 2.} We show $d_{g,\Sigma}$ is continuous. 

When 1) is satisfied,
the continuity of $d_{g,\Sigma}$ follows from Theorem \ref{2 option} (1).
When 2) holds, we choose $\delta$ such that $|\K_g|(B_\delta(p,g_0)\backslash\{p\})<\tau_0$, Corollary \ref{conv.length} (ii) yields 
\begin{align*}
d_g(\partial B_\delta(p,g_0),\partial B_r(p,g_0))
&=\lim_{k\rightarrow+\infty}d_{g_k}(\partial B_\delta(p,g_0),\partial B_r(p,g_0))\\
&=d(\partial B_\delta(p,g_0),\partial B_r(p,g_0))<C.
\end{align*}
%(Note that we only use the fact $d_g(\partial B_\delta^{g_0}(p),\partial B_r^{g_0}(p))\leq\lim_{k\rightarrow+\infty}d_{g_k}(\partial B_\delta^{g_0}(p),\partial B_r^{g_0}(p))$ here). 
By Theorem \ref{2 option} (3), we get the continuity of $d_{g,\Sigma}$.

\vspace{.1cm}

{\bf Step 3.} We show that a subsequence of $d_{g_k,\Sigma}$ converges in $C^0(\Sigma\times\Sigma)$ to a $C^0$ function.

We only need to prove the case when 1) is satisfied. By Proposition \ref{ghost},  
$$
\lim_{r\rightarrow 0}\lim_{k\rightarrow+\infty}\diam(B_r(p,g_0), d_{g_k,\Sigma})=0.
$$
By Lemma
\ref{Sobolev.distance}, $\|d_{g_k,\Sigma}\|_{W^{1,q}(\Omega\times\Omega)}<C(\Omega)$ for any $\Omega\subset\subset\Sigma\backslash \{p_0\}$ for some $q>2$, then
(a subsequence)
$d_{g_k,\Sigma}$ converges to a nonnegative function  $d$ in $C^{0,\alpha}_{\rm loc}(\Sigma\backslash \{p_0\}\times\Sigma\backslash \{p_0\})$. By Lemma \ref{main.lemma.1}, $d_{g_k,\Sigma}$
converges to a continuous function  $d$ uniformly on $\Sigma\times\Sigma$.

\vspace{.1cm}

{\bf Step 4.}  We claim $d = d_{g,\Sigma}$.  

First, by Proposition \ref{conv0}, Step 2 and Lemma \ref{main.lemma.2}, we have 
$
d\leq d_{g,\Sigma}.
$

Second,  we show 
$
d\geq d_{g,\Sigma}.
$
Let $\gamma_k$ be a curve from any $x_0$ and $y_0$ in $\Sigma$,
such that 
$$
d_{g_k,\Sigma}(x_0,y_0)\geq \ell_{g_k}(\gamma_k)-\epsilon.
$$

(i) If $x_0\not=p$ and $y_0\not=p$, we consider two cases: 

\vspace{.2cm}

\noindent Case 1: There exists $\delta$ such that $\gamma_k\subset\Sigma\backslash B_{\delta}(p,g_0)$.
For this case, it follows from Corollary \ref{conv.length} (i): 
$$
d(x_0,y_0)=\lim_{k\rightarrow+\infty}d_{g_k,\Sigma}(x_0,y_0)
\geq d_{g,\Sigma}(x_0,y_0).
$$

\noindent Case 2: For any $\delta$ there exist infinitely many $k$ such that 
$\gamma_k\cap B_\delta(p,g_0)\neq \emptyset$.  In this case,  let $t_1^k$ and $t_2^k$ be the first and the last $t$ such that $\gamma(t)\in\partial B_\delta(p,g_0)$ as $t$ increases, respectively.  Then 
$$
\lim_{k\rightarrow+\infty}\ell_{g_k}(\gamma_k)\geq \lim_{k\rightarrow+\infty}\ell_{g_k}(\gamma_k|_{[0,t_1]})
+\lim_{k\rightarrow+\infty}\ell_{g_k}(\gamma_k|_{[t_2,1]})
\geq d_{g,\Sigma}(x_0,x_1)+d_{g,\Sigma}(y_1,y_0),
$$
where $x_1, y_1\in\partial B_{\delta}(p,g_0)$ are the limits of $\gamma_k(t_1^k)$ and $\gamma_k(t_2^k)$ respectively.
Then 
$$
d(x_0,y_0)\geq d_{g,\Sigma}(x_0,y_0)-d_{g,\Sigma}(x_1,y_1)-\epsilon.
$$
Letting $\epsilon$ and $\delta\rightarrow 0$, we get the desired result.

(ii) We consider the case $x_0=p$.

Then we can choose $x_k\neq p$, $x_k\rightarrow x_0$, and get
$$
d(x_0,y_0)=\lim_{k\rightarrow+\infty}d(x_k,y_0)\geq\lim_{k\rightarrow+\infty}d_{g,\Sigma}(x_k,y_0)=d_{g,\Sigma}(x_0,y_0).
$$

In conclusion, we have established $d=d_{g,\Sigma}$. 

{\bf Step 5.} We prove $u,g$ are in fact independent of subsequences.  Assume $u,u'$ are two such limits and $g=e^{2u}g_0, g'=e^{2u'}g_0$. Since $\K_g=\K_{g'}$, 
$$
\int_\Sigma \left(K({g_0})\varphi + \nabla_{g_0}\varphi\nabla_{g_0}u\right) dV_{g_0} 
=\int_\Sigma \left(K({g_0})\varphi + \nabla_{g_0}\varphi\nabla_{g_0}u'\right) dV_{g_0}, \s \forall\varphi\in C^\infty_0(\Sigma). 
$$
Therefore $u-u'\in W^{1,q}$ is a weak solution of 
$$
\Delta_{g_0} (u-u') = 0
$$
and elliptic regularity implies $u-u'$ is smooth.  Since $\Sigma$ is closed,  $u-u' =c$ for some constant $c$.  The normalization $\diam(\Sigma,g_k)=1$ leads to $c=0$.  Therefore the given sequence $u_k$ converges weakly to a limit $u\in W^{1,q}$ and $d_{g_k,\Sigma}$ converges to $d_{g,\Sigma}$ uniformly, as claimed.  
\endproof

\medskip
%\vspace{.1cm}

\noindent{\em Proof of Theorem \ref{main-ChenLi}.} 
First, we show that for any sequence  $g_{k}$ satisfies the assumptions of Theorem \ref{main-ChenLi},  there exists $g=e^{2u}g_0\in \M(\Sigma,g_0)$, such that a subsequence of $d_{g_{k},\Sigma}$ converges to  $d_{g,\Sigma}$. 
We may assume $\K_{g_{k}}^+$
converges to $\mu'$ in the sense of distributions. By Hahn decomposition theorem (see \cite[Chapter 3]{Bogachev}), we may find a $\K_{g_k}$-measurable subset $A_k
\subset\Sigma$ such that 
$$
\K_{g_{k}}^+(E)=\K_{g_{k}}(E\cap A_k),\s \forall E.
$$
Then 
$$
\K_{g_{k}}^+(E)=\mu_{k}^1(E\cap A_k)-\mu_{k}^2(E\cap A_k)\leq \mu_{k}^1(E\cap A_k)\leq \mu_{k}^1(E).
$$
Then $\mu'\leq\mu^1$, so $\mu'(\{p\})<2\pi$ for any $p$. Applying Theorem \ref{main} completes the proof.
\endproof

\vspace{.1cm}

%In light of Theorem \ref{2 option}, Theorem \ref{main-ChenLi} has the following direct consequence:

%\begin{thm}\label{main2}
%Let $(\Sigma,g_0)$ be a closed Riemannian surface and $g_k=e^{2u_k}g_0\in\M(\Sigma)$.  Assume that $\K_{g_k}$ converges to a signed Radon measure $\mu$.  Assume $d_{g_k,\Sigma}$ converges to a continuous function $d$.Then $u_k$converges weakly to a function $u$ in $W^{1,q}$ for any $1\leq q<2$, $\K_g=\mu$ and $d_{g_k,\Sigma}$ converges to$d_{g,\Sigma}$ uniformly where $g=e^{2u}g_0$.
%\end{thm}

%\proof Similar to the proof of Theorem \ref{main}, we can obtain $\|u_k\|_{L^1}<C$. Let $\S$ be the finite set as in the proof of Theorem \ref{main}. Using the arguments in the Step 3 of the proof of Theorem \ref{main}, we have
%$$d|_{\Sigma\setminus S}\leq d_{g,\Sigma\setminus S}.$$
%By Theorem \ref{2 option}, $d_{g,\Sigma}$ is well-defined on $\Sigma$. By  Corollary \ref{continuity.distance}, $d_{g,\Sigma}$ is continuousand $d_{g,\Sigma}|_{\Sigma\setminus\S}=d_{g,\Sigma\setminus\S}$. Then we get $d\leq d_{g,\Sigma}$. The rest of the proof is the same with the proof of  Theorem \ref{main}. \endproof
%{\color{red} we need to demonstste all the requirements in these results are fufilled} {\color{red} Where is the $ diam = 1 $ assumption in Theorem 5.12 checked? } 

\subsection{Existence of approximation by smooth metrics} 
$\,$

\vspace{.2cm}

\begin{comment}
\begin{thm}\label{2pi}
Let $(\Sigma,g_0)$ be a closed Riemannian surface and
 $g=e^{2u}g_0\in \M(\Sigma,g_0)$ with $|\K_g|(\Sigma)<+\infty$. Assume $d_{g,\Sigma}$ is finite in $\Sigma\times\Sigma$. Then there exists a sequence of smooth metrics $g_k=e^{2u_k}g_0$, such that 
 \begin{enumerate}
 \item $u_k$ converges to  $u$ in $W^{1,q}(\Sigma,g_0)$ for any $q\in[1,2)$, 
 \item There exist nonnegative smooth functions $f_k^1$ and $f_k^2$, such that $K_{g_k}e^{2u_k}=f_k^1-f_k^2$ and $f_k^1dV_{g_0}$, $f_k^2dV_{g_0}$ converge to $\K_{g}^+,\K_{g}^-$ respectively in the sense of distributions,
 \item $d_{g_k,\Sigma}\to d_{g,\Sigma}$.  The metric $g$ is of bounded integral curvature in the sense of Alexandrov. 
\end{enumerate}
Moreover, for any $R$ and $x\in\Sigma$, we have
\begin{equation}\label{vol.com}
\frac{\vol(B_R(x,g))}{\pi R^2}\leq 1+\frac{1}{2\pi}\K_g^-(\Sigma).
\end{equation}
\end{thm}
\end{comment}

\noindent{\it Proof of Theorem \ref{2pi-main}}. {\bf Case 1.} $\K_{g}(\{x\})<2\pi$ for any $x\in \Sigma$.

By Proposition \ref{app.glo}, we can choose smooth metrics $g_k=e^{2u_k}g_0$ such that $u_k\to u$ in $W^{1,1}(\Sigma,g_0)$ and $\K_{g_k}\rightharpoonup\K_{g}$ as distributions. Since $\Sigma$ is closed, we can choose a scaling factor $\lambda_k\in\R$ so that 
$$
\diam(\Sigma,e^{2\lambda_k}g_k)=\diam(\Sigma,g).
$$  
Set $g_k':=e^{2\lambda_k}g_k=e^{2(u_k+\lambda_k)}g_0$.
By Theorem \ref{main} (applied to $g'_k$), $u_k+\lambda_k$
converges weakly in $W^{1,q}(\Sigma,g_0)$. So we can assume  $\lambda_k$ converges to a constant $\lambda$ as $u_k\rightharpoonup u$; so $u_k+\lambda_k\rightharpoonup u+\lambda$ in $W^{1,q}(\Sigma,g_0)$. By Theorem \ref{main} again,  $d_{g_k',\Sigma}\to d_{e^{2(u+\lambda)}g_0,\Sigma}=d_{e^{2\lambda}g,\Sigma}$ uniformly.
Then $\diam(\Sigma,e^{2\lambda}g)=\diam(\Sigma,g)$. Thus $\lambda=0$, so $d_{g_k,\Sigma}\to d_{g,\Sigma}$ in $C^0(\Sigma\times\Sigma)$. For $\epsilon>0$ and large $k$, 
$$
B_{R-\epsilon}(x,g_k)\subset B_R(x,g)\subset B_{R+\epsilon}(x,g_k).
$$

It is well-known that on the smooth Riemannian surface $(\Sigma,g_k)$ %(see Appendix), 
\begin{equation}\label{volume.comparison}
\frac{\vol(B_r(x),g_k)}{\pi r^2}\leq 1+\frac{1}{2\pi}\int_{B_r(x,g_k)}K_{g_k}^{-}dV_{g_k}.
\end{equation}
Let
$$
\S_0=\left\{x:|\K_g|(\{x\})\geq\frac{\pi}{2}\right\}\s \mbox{and}\s E_\epsilon=\bigcup_{y\in S_0}B_\epsilon(y,g_0).
$$ 
By 2) in Corollary \ref{MTI}, $e^{2u_k}$ is bounded in $L^3(\Sigma\backslash E_\epsilon)$, then $e^{2u_k}\to e^{2u}$ in $L^1(\Sigma\backslash E_\epsilon)$.
Hence
\begin{align}\label{vol.app}
\frac{\vol(B_R(x,g)\backslash E_\epsilon,g)}{\pi R^2}
&=\lim_{k\rightarrow\infty}\frac{\vol(B_R(x,g)\backslash E_\epsilon,g_k)}{\pi R^2}\nonumber\\
&\leq\frac{(R+\epsilon)^2}{R^2}\varlimsup_{k\rightarrow\infty} \frac{\vol(B_{R+\epsilon}(x,g_k),g_k)}{\pi (R+\epsilon)^2}\\&\leq \frac{(R+\epsilon)^2}{R^2}\left(1+\varlimsup_{k\rightarrow\infty}\frac{1}{2\pi}\int_\Sigma K_{g_k}^{-}dV_{g_k}\right)\nonumber\\
&\leq \frac{(R+\epsilon)^2}{R^2}\left(1+\varlimsup_{k\rightarrow\infty}\frac{1}{2\pi}\int_\Sigma f^2_k dV_{g_0}\right) \nonumber\\
&=\frac{(R+\epsilon)^2}{R^2}\left(1+\frac{1}{2\pi}\K_{g}^{-}(\Sigma)\right)\nonumber
\end{align}
where the last inequality follows from the Jordan decomposition theorem ($K^+_{g_k},K^-_{g_k}$ are mutually singular). 
Letting $\epsilon\rightarrow 0$, we have proved the theorem for Case 1.

\vspace{.1cm}

{\bf Case 2.} There is at least one $x\in\Sigma$ where $\K_g(\{x\})\geq 2\pi$.  

Since $d_{g,\Sigma}$ is finite,  in light of Theorem \ref{2 option}, we know $\K_{g}(\{x\})\leq 2\pi$ 
for any $x$ and (for simplicity) there is only
one point $p_0$ with $\K_{g}(\{p_0\})=2\pi$. In local isothermal coordinates with $p_0=0$ and $g_0=e^{2u_0}g_{\rm euc}$, we set $g=e^{2v}g_{\rm euc}$,  namely $u=v-u_0$.  %Put 
%$v=w+I_{\K_{g}}$, where $w$ is harmonic with $\|\nabla w\|_{L^q}<C(q)$ for any $q\in[1,2)$. Then we can write
%$$
%v=w-\log |x|, {\color{red} \mbox{why use $w_1$, $\log|x|$?}}
%%where $w_1$ is a harmonic function defined on $D$.
Let $\eta_\delta:\R\rightarrow[0,1]$ be a cut-off function which is 1 on $D_{2\delta}^c$ and 0 on $D_\delta$ with $|\nabla \eta_\delta|<\frac{2}{\delta}$ and $|\Delta\eta_\delta|<
\frac{C}{\delta^2}$. Set 
\begin{align*}
v_{\delta,k}
&=v+\frac{1-\eta_\delta}{k}\log |x|.
\end{align*}
Since
$$
|\nabla (v_{\delta,k}-v)|\leq\left|\nabla\left((1-\eta_\delta)\frac{\log|x|}{k}\right)\right|\leq\frac{|\log|x||+2}{k\delta}\chi_{D_{2\delta}},
$$
$$
|\Delta(v_{\delta,k}-v)|=\left|\frac{\log|x|}{k}\Delta\eta_\delta+2\nabla\eta_\delta\nabla\log|x|\right|\leq\frac{C(|\log|x||+1)}{k\delta^2},
$$
we can find $\delta_k\rightarrow 0$ such that $v_{\delta_k,k}\to v$ in $W^{1,q}$ and
$|\K_{e^{2v_k}g_{\rm euc}}(D)|\rightarrow |\K_{g}|(D)$.
Moreover, in $D_\delta$,
%$$
%v_{\delta,k}=w_1+I_{\K_g}+
%\frac{1}{k}\log|x|,
%$$
%then
$$
\K_{e^{2v_{\delta,k}}g_{\rm euc}}=-\Delta v_{\delta,k}=\K_{g}-\frac{2\pi}{k}\delta_0,
$$
hence
$$
\K_{e^{2v_{\delta,k}}g_{\rm euc}}(\{0\})=2\pi\left(1-\frac{1}{k}\right).
$$

Let $g_k'=g$ on $D_{2\delta}^c$ and $e^{2v_{\delta_k,k}}g_{\rm euc}$ on 
$D_{2\delta}$.  Let $g'_k=e^{2u'_k}g_0$. %Let $u_k'=\frac{1}{2}\log g_k'g_0^{-1}$, i.e. $g_k'=e^{2u_k'}g_0$.
Then $u_k'=u$ on $D_{2\delta}^c$, and over $D_{2\delta}$
$$
u_k'=v_{\delta_k,k}-u_0=v+(1-\eta_{\delta_k})\frac{1}{k}\log |x|-u_0=u+(1-\eta_{\delta_k})\frac{1}{k}\log |x|.
$$ 
Since $u_k'-u$ is smooth on $\Sigma\backslash\{p_0\}$, we have
$$
\K_{g_k'}(\{x\})=\K_{g}(\{x\})<2\pi,\s \forall x\neq p_0.
$$
Together with $\K_{g_k'}(\{p_0\})=2\pi(1-\frac{1}{k})$, we have $\K_{g_k'}(\{x\})<2\pi$ for any $x$. 
Clearly $u_k'\leq u$, then
$$
d_{g_k,\Sigma}\leq d_{g,\Sigma}.
$$ 

Next, we show that 
\begin{equation}\label{app.s}
d_{g_k,\Sigma}\geq d_{g,\Sigma}-\sup_{y,y'\in \partial D_{2\delta_k}}d_{g,\Sigma}(y,y').
\end{equation}
Let $x,x'\in\Sigma$. Assume $x\notin D_{2\delta_k}$ for simplicity. Let $\gamma$ be a curve from
$x$ to $x'$ with $\ell_{g_k}(\gamma)<d_{g_k,\Sigma}(x,x')+\epsilon$. If $\gamma\subset\Sigma\backslash D_{2\delta_k}$, then 
$\ell_{g_k}(\gamma)=\ell_{g}(\gamma)\geq d_{g,\Sigma}(x,x')$
which implies that 
$$
d_{g_k,\Sigma}(x,x')\geq d_{g,\Sigma}(x,x').
$$
If $\gamma\cap D_{2\delta_k}\neq\emptyset$ and $x'\in D_{2\delta_k}$, let $t_1$ be the first $t$ with
$\gamma(t_1)\in \partial D_{2\delta_k}$. Then
\begin{align*}
d_{g_k,\Sigma}(x,x')&+\epsilon\geq \ell_{g_k}(\gamma[0,t_1])= \ell_{g}(\gamma[0,t_1])\geq d_{g,\Sigma}(x,\gamma(t_1))\\
&\geq d_{g,\Sigma}(x,x')-d_{g,\Sigma}(\gamma(t_1),x')\geq d_{g,\Sigma}(x,x')-\sup_{y,y'\in \partial D_{2\delta_k}}d_{g,\Sigma}(y,y').
\end{align*}
If $\gamma\cap D_{2\delta_k}\neq\emptyset$ and $x'\notin D_{2\delta_k}$,  denote $t_1,t_2$ the first and the last $t$ with 
$\gamma(t)\in \partial D_{2\delta_k}$ respectively. Then
\begin{align*}
d_{g_k,\Sigma}(&x,x')+\epsilon\geq \ell_{g_k}(\gamma[0,t_1])+\ell_{g_k}(\gamma[t_2,1])= \ell_{g}(\gamma[0,t_1])+\ell_{g}(\gamma[t_2,1])\\
&\geq d_{g,\Sigma}(x,\gamma(t_1))+d_{g,\Sigma}(x',\gamma(t_2))\geq d_{g,\Sigma}(x,x')-\sup_{y,y'\in \partial D_{2\delta_k}}d_{g,\Sigma}(y,y').
\end{align*}
Letting $\epsilon\rightarrow 0$, we get \eqref{app.s}.
By Corollary \ref{continuity.distance},  $d_{g,\Sigma}$ is continuous, then 
$$
\lim_{k\rightarrow+\infty}
\sup_{y,y'\in \partial D_{2\delta_k}}d_{g,\Sigma}(y,y')= 0,
$$
hence 
$d_{g_k,\Sigma}$ converges to $d_{g,\Sigma}$. Also by \eqref{vol.app}, we get
\eqref{vol.com}.

\vspace{.1cm}

Since $\K_{g_k}(\{x\}) <2\pi$ for any $x$, we can find {smooth metric}
$g_k''=e^{2u_k''}g$ with $\|u_k'-u_k''\|_{W^{1,q}}<\frac{1}{k}$, $|\K_{g_k''}|(\Sigma)<C$ and $\|d_{g_k'',\Sigma}-d_{g_k,\Sigma}\|_{C^0}<\frac{1}{k}$. 
Moreover, it is easy to check \eqref{vol.app} still holds.

\vspace{.1cm}

Since $d_{g,\Sigma}(x,y)<+\infty$ for any $x,y\in\Sigma$ by the assumption, the metric $d_{g,\Sigma}$ is intrinsic by its definition ($\Sigma$ is a connected surface). We have just proved that we can approximate $g$ by smooth metrics, therefore $g$ has bounded integral curvature in the sense of Alexandrov.
%Then we complete the proof.
\endproof

\begin{cor}
Let $(\Sigma,g_0)$ be a closed Riemannian surface and $g_k=e^{2u_k}g_0\in\M(\Sigma,g_0)$.  If $u_k\to u$ in $L^1(\Sigma,g_0)$, $\K_{g_k}^-(\Sigma)<C$ and $d_{g_k,\Sigma}$ converges to a distance function $d$ uniformly, then 
$e^{2u_k}\to e^{2u}$ in $L^1$. 
\end{cor}

\proof
Let $\diam(\Sigma,d)=a$. Then $\diam(\Sigma,d_{g_k,\Sigma})\leq 2a$ for large $k$.  By \eqref{vol.com}, $\int_\Sigma e^{2u_k}dV_{g_0}<Ca^2$. Since $u_k\to u$ in $L^1$,  we can assume $u_k\to u$ almost everywhere.  By Fatou's Lemma,  
$$
\int_{\Sigma}e^{2u}dV_{g_0}\leq Ca^2<+\infty.
$$

Let
$$
\S_0=\left\{x:|\K_g|(\{x\})\geq \frac{\pi}{2}\right\}\s \mbox{and}\s E_\epsilon=\bigcup_{x\in\S_0}B_\epsilon(x,d).
$$ 
By 2) in Corollary \ref{MTI}, $e^{2u_k}$ is bounded in $L^3(\Sigma\backslash E_\epsilon)$, then $e^{2u_k}\to e^{2u}$ in $L^1(\Sigma\backslash E_\epsilon)$.

We may assume $B_\epsilon(x,d)\subset B_{2\epsilon}(x,d_{g_k,\Sigma})$ for each $x\in \S_0$. 
By \eqref{vol.com}, 
$$
\vol(B_{\epsilon}(x,d),g_k)\leq \vol(B_{2\epsilon}(x,d_{g_k,\Sigma}),g_k)\leq C\epsilon^2.
$$
Hence
\begin{align*}
\varlimsup_{k\rightarrow\infty}\int_{\Sigma}\left|e^{2u_k}-e^{2u}\right|dV_{g_0}
&\leq\varlimsup_{k\rightarrow\infty}\int_{\Sigma\setminus E_\epsilon}\left|e^{2u_k}-e^{2u}\right|dV_{g_0}\\
&+\varlimsup_{k\rightarrow\infty}\sum_{x\in\S_0}\left(\vol(B_\epsilon(x,d),g_k)+\int_{B_\epsilon(x,d)}e^{2u}dV_{g_0}\right)\\
&\leq\sum_{x\in\S_0}\left(C\epsilon^2+\int_{B_\epsilon(x,d)}e^{2u}dV_{g_0}\right).
\end{align*}
We complete the proof by letting $\epsilon\rightarrow 0$.
\endproof

\subsection{Complete noncompact orientable surfaces} 
We can truncate the ends of a compete noncompact surface then cap them off and verify the distance function on compactified surface is the restriction from the original one.  The previous results for compact surfaces then leads to existence of 
 smooth approximation for noncompact surfaces: 

\begin{cor}\label{complete}
Let $(\Sigma,g_0)$ be an open orientable surface and $g=e^{2u}g_0\in\mathcal{M}(\Sigma,g_0)$. If $d_{g,\Sigma}$ is finite for any two points and complete, then there exist complete smooth metrics $g_k=e^{2u_k}g_0$, such that 
\begin{enumerate} 
\item $u_k$ converges to  $u$ in $W^{1,q}_{\rm loc}(\Sigma,g_0)$ for any $q\in[1,2)$, 
\item There exist nonnegative smooth functions $f_k^1$ and $f_k^2$ such that $K_{g_k}e^{2u_k}=f_k^1-f_k^2$ and $f_k^1dV_{g_0}$, $f_k^2dV_{g_0}$ converge to $\K_{g}^+,\K_{g}^-$ as distributions, respectively,
\item $d_{g_k,\Sigma}\rightarrow d_{g,\Sigma}$ in $C^0_{\rm loc}(\Sigma\times\Sigma)$. 
\end{enumerate}
\end{cor}

\proof 
 
Since $r(x)=d_{g,\Sigma}(x,p)$ is continuous, we may choose a smooth function $f(x)$ on $\Sigma$ such that $\|r-f\|_{C^0(B_{10R}(p,g))}<\frac{1}{100}$.
Let $a,b\in (5R,6R)$ be regular values of $f$ and $a<b$, and define 
$$
\Omega_R=\{x:f<{b}\},\s\Omega'_R=\{x:a<f<b\}.
$$
Then $\partial\Omega_R$ consists of finitely many embedded closed curves $S^b_1,...,S^b_N$.  By Sard's theorem we can take $b$ close to $a$ so that $\partial\Omega'_R$ consists of embedded closed curves $S^a_1,...,S^a_N, S^b_1,...,S^b_N$ and each pair $\{S^a_i,S^b_i\}$ bounds a topological annulus ${\mathcal A}_i(R)$ in $\Omega_R$.  For simplicity, we assume $S_i^a\subset \{f=a\}$ and $S_i^b\subset \{f=b\}$. Clearly,  $\Omega_R'=\cup^N_{i=1}\mathcal{A}_i(R)$. %{\color{red} same as $\Omega'_R$?}

%Since $d_g$ is complete on the open surface $\Sigma$

Since $|\K_g|$ is locally finite on $\Sigma$, we may adjust $a$ so that $|\K_g|(\Omega'_R)<\frac{4}{3}\pi$, in turn, this implies
$$
|\K_{g}|(\{x\})<\frac{4}{3}\pi,\s\forall x\in\Omega'_R.
$$

By the uniformization theorem of Riemann surfaces,  each $({\mathcal A}_i(R),g_0)$ is conformal to an euclidean annulus $A_i(R):=D\backslash D_{\rho_i(R)}$  for some $\rho_i(R)\in(0,1)$, via a conformal diffeomorphism $\varphi_i:{\mathcal A}_i(R)\to A_i(R)$ mapping a neighborhood of $S^a_i$ in ${\mathcal A}_i(R)$ to a neighborhood of $\partial D$ in $A_i(R)$. %(outer boundary goes to {\color{red} outer inner?} boundary).  
Then we extend $\Omega_R$ to an  orientable surface by gluing a disk $D$ to $\Omega_R$ along each $S^b_i$: 
$$
\Sigma_R=\Omega_R\bigsqcup_{S_i^{b}}D.
$$
As $\Sigma$ is complete in $d_{g,\Sigma}$,  the surface $\Sigma_R$ is closed. 

Now we equip $\Sigma_R$ with a metric $g_R$ as follows.  First, in the standard complex coordinate $z\in\mathbb C$,  we write $g=e^{2u'}g_{\rm euc}$ and $g_0=e^{2u_0}g_{\rm euc}$ on $A_i(R)$ for each $i$, so $u=u'-u_0$ on $A_i(R)$. Select $b_i\in(\rho_i(R),1)$ and take a smooth cut-off function $\eta$ on $\C$ which is 1 on $D^c$ and 
0 on $D_{b_i}$ and $0\leq\eta\leq 1$. Define $g_{_{0,R}}=g_0$ on $\Omega_R\backslash \mathcal{A}_i(R)$ and $g_{_{0,R}}=e^{2\eta u_0}g_{\rm euc}$ on $D$ (here $D$ is a conformal parametrization of the union of $\mathcal{A}_i(R)$ and the glued disk). As $g_0$ is smooth $g_{_{0,R}}$ is a smooth metric on $\Sigma_R$ and it determines a conformal structure making $\Sigma_R$ a closed Riemann surface. 
Define $g_{_R}=g$ on $\Omega_R\backslash{\mathcal A}_i(R)$ and 
$g_{_R}=e^{2\eta u'}g_{\rm euc}$ on $D$. Then $g_{_R}\in\mathcal{M}(\Sigma_R,g_{_{0,R}})$ since
$$
-\Delta(\eta u') =\eta\K_g-2\nabla u'\nabla\eta-u'\Delta\eta.
$$
  Let $u_{_{R}}=u$ on $\Omega_R\backslash{\mathcal A}_i(R)$ and $u_{_R}=\eta u'-\eta u_0$ on $D$. Then
$g_{_R}=e^{2u_{_R}}g_{_{0,R}}$.
\vspace{.1cm}

Since $(\Sigma,d_{g,\Sigma})$ is complete,  there exists $R_n\to+\infty$ such that %{\color{red} not need:there corresponds well-defined $\Omega_{R_n}\not=\emptyset$ and 
%$d_{g,\Sigma}(\partial \Omega_{R_n},\partial\Omega_{R_{n+k}})\to+\infty$ as $k\to+\infty$. Hence, for each $n$ there is $k$ 
%such that 
%$$
%d_{g,\Sigma}(\partial\Omega_{R_n},\partial\Omega_{R_{n+k}})>diam(\Omega_{R_n},g)+10. 
%$$
%Then we can select a subsequence, which we still denote by $R_n$, such that 
\begin{equation}\label{+10}
d_{g,\Sigma}(\partial\Omega_{R_n},\partial\Omega_{R_{n+1}})>\diam(\Omega_{R_n},g)+10. 
\end{equation}
By Lemma \ref{suff.geodesic}, $\Omega_{R_{n+1}}$ is a quasi-geodesic convex neighborhood of $\Omega_{R_n}$.  %For simplicity, we may assume $\Omega_{R_n}$ is a quasi-geodesic convex neighborhood of $\Omega_{R_{n-1}}$ and $\Sigma=\cup_n^\infty\Omega_{R_n}$.  

%and $U_n=B_{\rho_n}^g(p)$. We can select $\rho_n$ and $R_n$ suitably, such that $\Omega_{n}$ is a geodesic neighborhood of $\Omega_{n-1}$.  Then
%$$diam(U_n,d_{g,\Omega_n})<d_{g}(\partial U_n,\partial\Omega_n)$$ which is equivalent to $$ diam(U_n,d_{g_{R_n},\Omega_n})<d_{g_{R_n}}(\partial U_{n},\partial\Omega_n).$$By Lemma \ref{quasiconvex}, we have$$diam(U_n,d_{g_{R_n},\Sigma_{R_n}})<d_{g_{R_n}}(\partial U_n,\partial\Omega_n).$$

It is easy to check that $|\K_{g_{_{R_n}}}|(\{x\})<\frac{4}{3}\pi$ for any $x\in\Sigma_{R_n}\backslash\Omega_{R_n}$ and $d_{g_{_{R_n}},\Sigma_{R_n}}$ is finite on $\Omega_{R_n}$. So $d_{g_{_{R_n}},\Sigma_{R_n}}$ is finite on $\Sigma_{R_n}$. %For each fixed $n$, a
Applying Theorem \ref{2pi-main} to $g_{_{R_n}}=e^{2u_{R_{n}}}g_{_{0,R_{n}}}\in{\mathcal M}(\Sigma_{R_{n}},g_{_{0,R_{n}}})$ we can find $u_{R_n}'\in C^\infty(\Sigma_{R_{n}})$ (then set $g'_{n}=e^{2u'_{{R_n}}}g_{_{0,R_n}}$) such that

\begin{enumerate}
\item[{1)}] $\|u_{_{R_n}}'-u_{_{R_n}}\|_{W^{1,q}(\Sigma_{R_n},g_{_{0,R_n}})}<\frac{1}{n}$, \hfill(existence of a converging sequence ) 
\item[{2)}] $|\K_{g_n'}|(\Omega_{R_i})\leq |\K_{g_{_{R_n}}}|(\Omega_{R_i})+\frac{1}{n}$, $i=1,\cdots,n,$ \hfill(convergence in measure)
\item[{3)}] $\|d_{g_n',\Sigma_{R_{n}}}-d_{g_{_{R_n}},\Sigma_{R_{n}}}\|_{C^0(\Sigma_{R_{n}}\times\Sigma_{R_{n}})}<\frac{1}{n}$. \hfill(convergence of distances)
\item[{4)}] $d_{g_n',\Sigma_{R_n}}(\partial\Omega_{R_n},\partial\Omega_{R_{n+1}})>\diam(\Omega_{R_n},g_n')$
\end{enumerate}
where 4) can be seen from 3) and the choice of $\Omega_{R_n}$ as follows 
\begin{align*}
d_{g_n',\Sigma_{R_n}}&(\partial\Omega_{R_n},\partial\Omega_{R_{n+1}})>d_{g,\Sigma}(\partial\Omega_{R_n},\partial\Omega_{R_{n+1}})-\frac{1}{n}\geq \diam(\Omega_{R_n},g)+10-\frac{1}{n}\\
&= \diam(\Omega_{R_n},g'_n)+10-\frac{1}{n}>\diam(\Omega_{R_n},g'_n).
\end{align*}
Consequently,  $\Omega_{R_{n+1}}$ is a quasi-geodesic convex domain of $\Omega_{R_n}$ in $(\Sigma_{R_n},g_n')$, by 4) and Lemma \ref{suff.geodesic}.

Select $\varphi_n\in\C^\infty(\Sigma)$ such that 
$$
e^{\varphi_n}|_{\Omega_{R_{n+2}}\backslash \Omega_{R_{n+1}}}\geq \frac{n}{d_{g_0}(\partial \Omega_{R_{n+2}},\partial \Omega_{R_{n+1}})}.
$$
On $\Sigma$, we define 
$$
u_n=\eta_nu_{R_n}'+(1-\eta_n)\varphi_n\s\mbox{and} \s g_n=e^{2u_n}g_0
$$ 
where we pick $\eta_n\in C^\infty(\Sigma)$ with $\eta_n=1$ on $\Omega_{R_n}$ and 
$\eta_n=0$ on $\Sigma\backslash\Omega_{R_{n+1}}$. Then
$d_{g_n,\Sigma}$ is complete on $\Sigma$ as any curve $\gamma$ connecting $\partial\Omega_{R_{n+1}}$ and $\partial\Omega_{R_{n+2}}$ %in $\Sigma\backslash\Omega_{R_n}$ 
enjoys $\ell_{g_n}(\gamma)\geq n$. 

Before move on, let us summarize the notations. On $\Omega_{R_n}\backslash\mathcal{A}(R_n)$: 
%$u_{_{R_n}}=u$ and 

\vspace{.1cm}

\noindent $\textcircled{1}\, g_{_{0,{R_n}}}\!= g_0,\textcircled{2}\, u_n \!= u'_{R_n},\textcircled{3}\, g_n\!=g_n' =e^{2u_n}g_0= e^{2u'_{R_n}} g_{_{0,R_n}}$, 
$\textcircled{4}\, u_{_{R_n}} \!= u$,  $\textcircled{5}\, g_{_{R_n}}\!= g. $%= e^{2u_{_{R_n}}}g_{_{0,R_n}}\!=e^{2u}g_{_{0}}$.

\vspace{.1cm}

%\begin{itemize}
%\item[$\textcircled{1}$]$g_{_{0,{R_n}}} = g_0$, \s $\textcircled{2}$ $u_n = u'_{R_n}$, \s $\textcircled{3}$ $g_n=g_n'\,(=e^{2u_n}g_0= e^{2u'_{R_n}} g_{_{0,R_n}}$), 
%\item[$\textcircled{4}$]$u_{_{R_n}} = u$, \s $\textcircled{5}$  $g_{_{R_n}}= g \,(= e^{2u_{_{R_n}}}g_{_{0,R_n}}=e^{2u}g_{_{0}})$.
%\end{itemize}

Substituting $\textcircled{1}$,  $\textcircled{2}$,  $\textcircled{4}$ into 1) yields 
$\|u_n-u\|_{W^{1,q}(\Omega_{_{R_n}},g_0)}\rightarrow 0$. Substituting $\textcircled{3},  \textcircled{5}$ into 2) gives
$$
|\K_{g_n}|(\Omega_{R_i})<C(i), \s \forall i.
$$
So we can assume $\K_{g_n}$ converges to a measure $\mu$
weakly, and $\mu=\K_g$ since $u_k\to u$ in $W^{1,q}_{\rm loc}(\Sigma, g_0)$.

By $\textcircled{3}$ and $\textcircled{5}$, $\Omega_{R_{n+1}}$ is a quasi-geodesic convex neighborhood of 
$\Omega_{R_n}$ in $(\Sigma_{R_n},g_{_{R_n}})$ and $(\Sigma,g_n)$. Hence
$$
d_{g_n,\Sigma}|_{\Omega_{R_n}\times\Omega_{R_n}}=d_{g_n,\Omega_{R_{n+1}}}|_{\Omega_{R_n}\times\Omega_{R_n}}\stackrel{\textcircled{3}}{=}d_{g_n',\Omega_{R_{n+1}}}|_{\Omega_{R_n}\times\Omega_{R_n}}=d_{g_n',\Sigma_{R_n}}|_{\Omega_{R_n}\times\Omega_{R_n}},
$$
%and
$$
d_{g_{_{R_{n}}},\Sigma_{R_n}}|_{\Omega_{R_n}\times\Omega_{R_n}}=d_{g_{_{R_{n}}},\Omega_{R_{n+1}}}|_{\Omega_{R_n}\times\Omega_{R_n}}\stackrel{\textcircled{5} }{=}d_{g,\Omega_{R_{n+1}}}|_{\Omega_{R_n}\times\Omega_{R_n}}=d_{g,\Sigma}|_{\Omega_{R_n}\times\Omega_{R_n}}.
$$
Substituting them into 3), we have
$$
\|d_{g_n,\Sigma}-d_{g,\Sigma}\|_{C^0(\Omega_{R_n}\times\Omega_{R_n})}\rightarrow 0.
$$

In conclusion, the sequence $(u_n,g_n)$ satisfies the requirements.% in the corollary. 
 \endproof

\vspace{.1cm}

Similar arguments lead to convergence of distance functions on complete orientable surface. 

\begin{cor}
Let $(\Sigma,g_0)$ be an open orientable surface and $g_k=e^{2u_k}g_0$, $g=e^{2u}g_0\in\mathcal{M}(\Sigma,g_0)$. Assume  $d_{g_k,\Sigma}$ and $d_{g,\Sigma}$ are finite and complete, and
$\K_{g}(\{x\})<2\pi$ for all $x\in\Sigma$. If $|\K_{g_k}|(E)<C(E)$
for any compact subset $E\subset \Sigma$ and $u_k$ converges to
$u$  in $L^{1}_{\rm loc}(\Sigma)$, then $d_{g_k,\Sigma}$ converges to
$d_{g,\Sigma}$ in $C^0_{\rm loc}(\Sigma\times\Sigma)$.
\end{cor}

\proof 

As in the proof of Corollary \ref{complete}, define $\Omega_{R}$, $g_{_R}$, $u_{_R}$, and $u_{_{k,R}}=u_k$ on $\Omega_R\backslash{\mathcal A}_i(R)$, $u_{_{k,R}}=\eta u_k'-\eta u_0$ on $D$ where $u_k'=u_k+u_0$. 
Define $g_{_{k,R}}=g_k$ on $\Omega_R\backslash{\mathcal A}_i$(R) and 
$g_{_{k,R}}=e^{2\eta u_k'}g_{\rm euc}$ on $D$. 
It is easy to check that $u_{_{k,R}}$ converges to $u_{_R}$ weakly in $W^{1,q}(\Sigma_{R})$, and $\K_{g_R}(\{x\})<2\pi$ for any $x\in\Sigma_R$. 
\vspace{.1cm}

As in the proof of Corollary \ref{complete}, we can select  $R_n\to+\infty$ such that  
$$
d_{g,\Sigma}(\partial\Omega_{R_n},\partial\Omega_{R_{n+1}})>\diam(\Omega_{R_n},g)+10. 
$$
By Lemma \ref{suff.geodesic}, $\Omega_{R_{n+1}}$ is a quasi-geodesic convex neighborhood of $\Omega_{R_n}$, then
$$
\diam(\Omega_{R_n},d_{g,\Sigma})=\diam(\Omega_{R_n},d_{g,\Omega_{R_{n+1}}})=\diam(\Omega_{R_n},d_{g_{_{R_n}},\Sigma_{R_n}}),
$$
hence
\begin{equation}\label{g.geodesic}
d_{g_{_{k,R_n}},\Sigma_{R_n}}(\partial\Omega_{R_n},\partial\Omega_{R_{n+1}})>\diam(\Omega_{R_n},g_{_{k,R_n}})+10. 
\end{equation}
We may assume $\Omega_{R_{n+1}}$ is a quasi-geodesic convex neighborhood of $\Omega_{R_{n}}$ and $\Sigma=\cup_n^\infty\Omega_{R_n}$.  

Select scalings $c_{_{k,R_n}}$ so that $\diam(\Sigma_{R_n},e^{2c_{_{k,R_n}}}g_{_{k,R_n}})=1$. 
Applying Theorem \ref{main} to $e^{2c_{_{k,R_n}}}g_{_{R_n}}$,  we know $u_{_{k,R_n}}+c_{_{k,R_n}}\rightharpoonup$ a function $v$ in $W^{1,q}(\Sigma_{R_n},g_{_{0,R_n}})$
and $d_{e^{2c_{{k,R_n}}}g_{_{k,R_n}},\Sigma_{R_{n}}}\to d_{e^{2v}g_{_{0,R_n}},\Sigma_{R_{n}}}$. As $u_{_{k,R_n}}\to u_{_{R_n}}$, we  assume $c_{_{k,R_n}}$ is convergent. Then $d_{g_{_{k,R_n}},\Sigma_{R_{n}}}\to d_{g_{_{R_n}},\Sigma_{R_{n}}}$ in $C^0$.  %(\Sigma_{R_n}\times\Sigma_{R_n})$. 
By \eqref{g.geodesic},
$$
d_{g_{_{k,R_n}}}(\partial\Omega_{R_n},\partial\Omega_{R_{n+1}})>\diam(\Omega_{R_n},g_{_{k,R_n}}),\s \mbox{for large $n$}.
$$
So $\Omega_{R_{n+1}}$
 is also a quasi-geodesic convex neighborhood of $\Omega_{R_{n}}$
in $(\Sigma_{R_n},g_{_{k,{R_n}}})$ and $(\Sigma,g_k)$. Then
$$
d_{g,\Sigma}|_{\Omega_{R_n}\times\Omega_{R_n}}=d_{g,\Omega_{R_{n+1}}}|_{\Omega_{R_n}\times\Omega_{R_n}}=d_{g_{_{R_n}},\Omega_{R_{n+1}}}|_{\Omega_{R_n}\times\Omega_{R_n}}=d_{g_{_{R_n}},\Sigma_{R_{n}}}|_{\Omega_{R_n}\times\Omega_{R_n}},
$$
and
$$ 
d_{g_{_{k,R_{n}}},\Sigma_{R_n}}|_{\Omega_{R_n}\times\Omega_{R_n}}=d_{g_{_{k,R_{n}}},\Omega_{R_{n+1}}}|_{\Omega_{R_n}\times\Omega_{R_n}}=d_{g_k,\Omega_{R_{n+1}}}|_{\Omega_{R_n}\times\Omega_{R_n}}=d_{g_k,\Sigma}|_{\Omega_{R_n}\times\Omega_{R_n}}.
$$
It follows 
$$
\|{d_{g_k,\Sigma}}-d_{g,\Sigma}\|_{C^0(\Omega_{R_n}\times\Omega_{R_n})}\rightarrow 0.
$$
\endproof

\section{Convergence of distance functions in varying conformal classes}

%In this section, we will treat closed orientable surfaces via results from Riemann surfaces.  The closed nonorientable surfaces will be dealt in the next section by going to the orientation double cover. The term {\it moduli space} is for closed Riemann surfaces of fixed genus. 

%\subsection{A Mumford type lemma for length of loops in nonsmooth metrics}

%Suppose that $\Sigma$ is a closed orientable surface of genus $\geq 1$.  Assume that $h_k$ and $h_0$ are smooth metrics on $\Sigma$ with $h_k\rightarrow h_0$ in the $C^2$-topology {\color{red} and the Gauss curvature $K_{h_k}=-1$ or $0$.

\subsection{A Mumford type lemma for nonsmooth metrics}
A sequence of metrics conformal to metrics of constant curvature with nonsmooth conformal factors can be convergent  under suitable assumption in the spirit of Mumford's compactness theorem. This is important when we consider distance convergence for varying background conformal classes. 
\begin{lem}\label{Mumford}
Let $\Sigma$ be a closed (orientable or nonorientbale) surface, and let $h_k$ be a smooth metric with Gauss curvature $K_{h_k}=-1,0$. 
%When $K_{h_k}=0$, we assume $Area(h_k)=1$. 
Let $g_k=e^{2u_k}h_k$ with $u_k\in W^{1,1}(\Sigma)$ and $\vol(g_k)+|\K_{g_k}| (\Sigma)<C$.  If  the conformal class of $h_k$ converges to the boundary of the moduli space,
 then there exists $\gamma_k:S^1\rightarrow\Sigma$,
such that $[\gamma_k]\neq 1$ in $\pi_1(\Sigma)$
and $\ell_{g_k}(\gamma_k)\rightarrow 0$.   
\end{lem}

\proof 

We first consider the case $\Sigma$ is orientable.  Assume $K_{h_k}=-1$. Let $\gamma_k$ be the shortest geodesic loop, and set
$$
w_k= \operatorname{arcsinh}\frac{1}{\sinh(\frac{1}{2}\ell(\gamma_k))}.
$$
By Mumford's compactness theorem (cf. \cite{M}), $\ell(\gamma_k)\rightarrow 0$.
Applying \cite[Theorem 4.1.1]{Buser}, there exists a domain $U_k\subset\Sigma$, which is isometric to
$ S^1\times(-w_k,w_k)$ with the metric 
$$
g=d\rho^2+ \ell^2(\gamma_k)\cosh^2\rho ds^2=\left(\frac{\ell(\gamma_k)\cosh\rho}{2\pi}\right)^2\left(\left(\frac{2\pi d\rho}{\ell(\gamma_k)\cosh\rho}\right)^2+d\theta^2\right),
$$ 
where $s=\frac{\theta}{2\pi}$.
Let 
$$
(t,\theta)=\phi_k(\rho,\theta)=\left(\frac{4\pi\arctan e^\rho}{\ell(\gamma_k)},\theta\right).
$$
Then $\phi_k$ is a diffeomorphism from $S^1\times (-w_k,w_k)$
to 
$$
Q_k=\left(\frac{4\pi\arctan e^{-w_k}}{\ell(\gamma_k)},\frac{4\pi\arctan e^{w_k}}{\ell(\gamma_k)}\right)\times S^1
$$
with
$$
\phi_k^*(g)=\ell^2(\gamma_k)\cosh^2\rho(dt^2+d\theta^2)=\left(\frac{\ell(\gamma_k)}{2\pi\sin\frac{\ell(\gamma_k)t}{2\pi}}\right)^2(dt^2+d\theta^2).
$$
Hence $\Omega_k$ is conformal to $Q_k$. Note that 
$$\frac{\arctan e^{w_k}-\arctan e^{-w_k}}{\ell(\gamma_k)}\rightarrow+\infty.$$ After a translation, we may assume $\Omega_k$ is conformal to $S^1\times(-T_k,T_k)$ with $T_k\rightarrow+\infty$.

When $K_{h_k}=0$,
 $(\Sigma,h_k)$ is induced by lattice
$\{1,a_k+b_k\sqrt{-1}\}$ in $\C$, where $-\frac{1}{2}<a_k\leq\frac{1}{2}$, $b_k>0$, $a_k^2+b_k^2\geq 1$, and
$a_k\geq 0$ whenever $a_k^2+b_k^2=1$.
$(\Sigma,h_k)$ is conformal to $(S^1\times\R)/G_k$,
where $G_k\cong \mathbb{Z}$ is the transformation group of
$S^1\times\R$ generated by
$
(\theta,t)\rightarrow (\theta+2\pi a_k,t+2\pi b_k).
$
In the moduli space $\mathcal{M}_1$ of tori, $(\Sigma,h_k)$ diverges if and only if  $b_k\rightarrow+\infty$.
Then $S^1\times (1,b_k-1)$ is a domain of $\Sigma$ which is conformal to $S^1\times(-T_k,T_k)$ with $T_k\rightarrow+\infty$.

%If the conformal class of $h_k$ converges to the boundary of the moduli space,  

In conclusion, there exists a domain $U_k\subset\Sigma$,
such that $(U_k,h_k)$ is conformal to $S^1\times (-T_k,T_k)$
with $T_k\rightarrow+\infty$, and $S^1\times\{t\}$ is
nontrivial in $\pi_1(\Sigma)$.
Then we can consider $g_k|_U$,  as a conformal 
metric with nonsmooth coefficients on $S^1\times(-T_k,T_k)$ since we can write $g_k=e^{2v_k}(dt^2+d\theta^2)$. 

Since 
$$
\mbox{Area}(U_k,g_k)=\int_{-T_k}^{T_k}\int_{0}^{2\pi}e^{2v_k}dtd\theta\leq C, 
$$
we can find $a_k\in[-T_k+1,T_k-2]$, such that
$$
\int_{a_k}^{a_k+1}\int_{0}^{2\pi}e^{2v_k}dtd\theta\rightarrow 0,
$$
which yields that 
$$
\int_{a_k}^{a_k+1}\int_{0}^{2\pi}e^{v_k}dtd\theta\rightarrow 0.
$$
Then we can find $b_k\in[a_k,a_k+1]$, such that
$$
\int_{0}^{2\pi}e^{v_k(b_k,\theta)}d\theta\rightarrow 0.
$$
Let $\gamma_k=S^1\times\{b_k\}$. Then $\ell_{g_k}(\gamma_k)\rightarrow 0$.

\vspace{.1cm}

Next, when $\Sigma$ is non-orientable, let $\pi:\widetilde\Sigma\rightarrow\Sigma$ be the orientation covering map and $\widetilde h_k=\pi^*(h_k)$. Then $\K_{h_k}=-1$ or $0$. Moreover, we have $\sigma^*(\widetilde h_k)=\widetilde h_k$, where $\sigma$ is the nontrivial covering transformation.  Let $\widetilde g_k=e^{2\tilde{u}_k}\widetilde h_k$ where $\widetilde u_k(x)= u_k(\pi(x))$. We have
$$
\mbox{Area}(\widetilde\Sigma,\widetilde g_k)=2\mbox{Area}(\Sigma,g_k)<C.
$$
Suppose $\{h_k\}$ is  not compact in $C^\infty$. Then $\{\widetilde h_k\}$
 is also not compact in $C^\infty(\widetilde{\Sigma}$).  By Lemma \ref{Mumford}, there exists a noncontractible circle
$\widetilde\gamma_k:S^1\rightarrow\widetilde\Sigma$, such that $\ell_{\widetilde g_k}(\widetilde\gamma_k)\rightarrow 0$. Let $\gamma_k=\pi(\widetilde\gamma_k)$. Then $\ell_{g_k}(\gamma_k)=\ell_{\widetilde g_k}(\widetilde\gamma_k)\rightarrow 0$.  By \cite[Theorem 4.1]{Massey},  the loop $\gamma_k$ is nontrivial in $\pi_1(\Sigma)$. %Hence $\co(\Sigma,g_k)\rightarrow 0$. But this contradicts the assumption that there is a positive lower bound for $\co(\Sigma,g_k)$. 
\endproof

\subsection{Metric convergence with nondegenerating conformal classes}

\begin{lem}\label{conv.conf}
Let $\Sigma$ be a closed surface of genus $\geq 1$. 
Assume that $h_k$ and $h_0$ are smooth metrics on $\Sigma$ with  $h_k\rightarrow h_0$ in the $C^2$-topology and $K_{h_k}=-1$ or $0$. Let 
$g_k=e^{2u_k}h_k\in\mathcal{M}(\Sigma,h_k)$, $g=e^{2u}h_0\in\mathcal{M}(\Sigma,h_0)$.  Suppose $|\K_{g_k}|$ and
$\K_{g_k}^+$ converges to measure $\nu$ and $\mu$ respectively. Then
\begin{itemize}
\item[0)] If $\|u_k\|_{L^1(B_\delta^{h_0}(p),h_k)}<A$, and $|\K_{g_k}|(B_\delta^{h_0}(p))<\tau$, then there exists $r$ so that $\int_{B_r^{h_0}(p)}e^{\frac{4\pi}{\tau}u_k}dV_{h_k}<C(\tau)$.
Moreover, for any $p'$ with $\K_{g_k}(\{p'\})<2\pi$,  $\|\nabla^{h_0} d_{g_k,\Sigma}(p',x)\|_{L^\frac{4\pi}{\tau}(B_r^{h_0}(p),h_0)}<C(\tau)$.
\item[1)] If $u_k\rightarrow u$
in $L^1(\Sigma,h_0)$, and $|\K_{g_k}|(B_\delta^{h_0}(p))<\tau_0$, then there exists $r$, such that 
$d_{g_k,\Sigma}$ converges to $d_{g,\Sigma}$ in 
$C^0(B_{r}^{h_0}(p)\times B_{r}^{h_0}(p))$;
\item[2)] Let  $K$ be a compact subset. If $u_k\rightarrow u$
in $L^1(\Sigma,h_0)$ and  $\mu(\{x\})<\tau_0$ for any $x\in K$, then for any $\gamma_k\subset K$ with $\gamma_k(0)\rightarrow x$,
$\gamma_k(1)\rightarrow y$, 
there holds
$$
\varliminf_{k\rightarrow+\infty}\ell_{g_k}(\gamma_k)\geq d_{g,\Sigma}(x,y).
$$
\item[3)]Let $U,V$ be compact domains in $\Sigma$ and $U\subset V$.
If $u_k\rightarrow u$
in $L^1(\Sigma,h_0)$ and $\mu(\{x\})<c_0$ in  $\overline{V\backslash U}$
then 
$$
d_{g}(\partial U,\partial V)= \lim_{k\rightarrow+\infty}d_{g_k,\Sigma}(\partial U,\partial V).
$$
\item[4)] If $\diam(\Sigma,d_{g_k})<\ell_0$, then for any $p$ with $\mu(\{p\})<2\pi$, we have
$$
\lim_{r\rightarrow 0}\lim_{k\rightarrow+\infty}\diam(B_r^{h_0}(p),d_{g_k,\Sigma})=0.
$$

\end{itemize}

\end{lem}

\proof Assume first that $\Sigma$ is orientable. Let $(\widetilde\Sigma,\tilde h)$ be the hyperbolic plane ${\mathbb H}^2$ or ${\mathbb R}^2$. Let $\pi_k:\widetilde{\Sigma}\to\Sigma$ be the Riemannian covering of $(\Sigma,h_k)$ such that the deck transformation on $\widetilde\Sigma$ by the fundamental group of $\Sigma$ acts by isometries on $\widetilde\Sigma$.  Fix any $p\in\Sigma$,  and then fix $\tilde p\in\widetilde\Sigma$ such that $\pi_0(\tilde p)=p$. Take $\rho< \min\left\{\frac{1}{2}\mbox{inj}(h_0),\delta\right\}$.  As $h_k\to h_0$ in $C^2$,  the injectivity radius of $h_k$ converges to $\mbox{inj}(h_0)$ and $\pi_k\to \pi_0$ in $C^0$. 
For each $k$, let $\widetilde{U}_k$ be the connected component of $\pi^{-1}_k(B_\rho^{h_0}(p))$ that contains $\tilde p$.  The restriction of $\pi_k$ on $\widetilde{U}_k$ to $B^{h_0}_\rho(p)$ is diffeomorphic. The intersection of all $\widetilde{U}_k$'s is a non-empty open set and denote it by $\widetilde{U}_p$.  Then $\pi_k$ is isometric from $\widetilde{U}_p$ to an open set $U_p\subset B_\rho^{h_0}(p)$ that contains $p$.  In particular, $\pi_k$ is a harmonic map from $(\widetilde U_p,\tilde h)$ to $(U_p,h_k)$ and the energy density of $\pi_k$ is uniformly bounded as $h_k\to h_0$, and then together with the $C^0$ convergence $\pi_k\to\pi_0$ we conclude $\pi_k\to\pi_0$ in $C^\infty_{\rm loc}$.
So $\pi_k^{-1}\to \pi^{-1}_0$ in $C^\infty_{\rm loc}(U_p)$.  

Let $\phi_{\tilde p}:D\rightarrow \widetilde{U}_p$ define an isothermal coordinate system.  Then 
$\vartheta_{k}:=\pi_k\circ\phi_{\tilde p}$ define an isothermal chart of $(U_p,h_k)$ and $\vartheta_{k}\to\vartheta_{0}$ in $C^0$. Pulling back to $D$, there are $\varphi_k,\varphi_0\in C^\infty(D)$ so that
$$
\vartheta_{k}^*(h_k)=e^{2\varphi_k}g_{\rm euc}\s {\rm and} \s \vartheta_{0}^*(h_0)=e^{2\varphi_0}g_{\rm euc}.
$$

Put
$$
g'_k=e^{2u_k\circ\vartheta_{k}+2\varphi_{k}}g_{\rm euc},\s g'=e^{2u\circ\vartheta_{0}+2\varphi_{0}}g_{\rm euc}.
$$
By Corollary \ref{MTI}, 
$$
\int_{D_\rho}e^{\frac{4\pi}{\tau}(u_k\circ\vartheta_{k}+\varphi_{k})}dx<C.
$$

By Lemma \ref{Sobolev.distance},  $|\nabla^{h_k}d(p,x)|\leq e^{u(x)}$, where $p$ is a fixed point. Then $$\|\nabla^{h_k}d_{g_k}(p,\cdot)\|_{L^\frac{4\pi}{\tau}( B_{r_x}^{h_0}(x),h_k)}<C$$ where $C$ is independent of $k$. 
However, 
$$
\|\nabla^{h_0}d_{g_k}(p,\cdot)\|_{L^\frac{4\pi}{\tau}( B_{r_x}^{h_0}(x),h_0)}\leq C\|\nabla^{h_k}d_{g_k}(p,\cdot)\|_{L^\frac{4\pi}{\tau}( B_{r_x}^{h_0}(x),h_k)}.
$$
Thus 0) holds.

Now, we prove 1):  
By Proposition \ref{conv0},
$d_{g_k',D}\to d_{g',D}$ on  $D_{r_0}$,  it follows $d_{g_k,\Sigma}(\vartheta_k(x),\vartheta_k(x'))\to d_{g,\Sigma}(\vartheta_0(x),\vartheta_0(x'))$ on $D_{r_0}$.
Choose $r<\rho$ such that $B_r^{h_0}(p)\subset\subset U_p$ and $\phi^{-1}_{\tilde p}(\pi_k^{-1}B_r^{h_0}(p))\subset D_{r_0}$ for large $k$. Note that
$$
\vartheta_0\circ\vartheta_k^{-1}=\pi\circ \pi_k^{-1}.
$$
Then we may assume $\vartheta_0\circ\vartheta_k^{-1}$ converges to identity map in $C^1(U_p)$. Let $y,y'\in B_r^{h_0}(p)$ and $x_k=\vartheta^{-1}_k(y)$ and $x_k'=\vartheta^{-1}_k(y')$. Then
\begin{align*}
&\left|d_{g_k,\Sigma}(y,y')\right.-\left.d_{g,\Sigma}(y,y')\right|\leq \left|d_{g_k',D}(x_k,x_k')-d_{g',D}(x_k,x_k')\right|+\left|d_{g',D}(x_k,x_k')-d_{g,\Sigma}(y,y')\right|\\
&\leq \left|d_{g_k',D}(x_k,x_k')-d_{g',D}(x_k,x_k')\right|+\left|d_{g,\Sigma}(\vartheta_0 \circ \vartheta_k^{-1}(y),\vartheta_0 \circ \vartheta_k^{-1}(y'))-d_{g,\Sigma}(y,y')\right|\\
&\leq \|d_{g_k',D}-d_{g',D}\|_{C^0(D_{r_0}\times D_{r_0})}+\left|d_{g,\Sigma}(\vartheta_0 \circ \vartheta_k^{-1}(y),\vartheta_0 \circ \vartheta_k^{-1}(y'))-d_{g,\Sigma}(y,y')\right|.
\end{align*}
Together with the continuity of $d_{g,\Sigma}$ on $B_{2r}^{h_0}(p)$, we get 1). 

\vspace{.1cm}

The proof of 2) is almost the same as the proof of Corollary \ref{conv.length} (i), we omit it.

\vspace{.1cm}

Next, we prove 3): For any $\epsilon>0$, let $\gamma$ be a curve in $\overline{V\backslash U}$ with $\gamma(0)\in\partial U,\gamma(1)\in\partial V$ and
$$
\ell_g(\gamma)\leq d_{g,\Sigma}(\partial U,\partial V)+\epsilon.
$$ 
By 1), we can select $t_0=0<t_1<\cdots<t_m=1$ such that
$$
d_{g,\Sigma}(\gamma(t_i),\gamma(t_{i+1}))=
\lim_{k\rightarrow+\infty}d_{g_k,\Sigma}(\gamma(t_i),\gamma(t_{i+1})).
$$
Then 
\begin{align*}
\ell_g(\gamma)&\geq \sum d_{g,\Sigma}(\gamma(t_i),\gamma(t_{i+1}))=
\sum\lim_{k}d_{g_k,\Sigma}(\gamma(t_i),\gamma(t_{i+1}))\\
&\geq \varlimsup_{k}d_{g_k,\Sigma}(\gamma(t_0),\gamma(1))\geq \varlimsup_k d_{g_k}(\partial U,\partial V).
\end{align*}
Hence
$$
\varlimsup_k d_{g_k}(\partial U,\partial V)\leq  d_{g}(\partial U,\partial V)+\epsilon.
$$
On the other hand, we select a curve $\gamma_k$ in $\overline{V\backslash U}$ with $x_k=\gamma_k(0)\in\partial U,y_k=\gamma_k(1)\in\partial V$ and
$$
\ell_{g_k}(\gamma_k)\leq d_{g_k,\Sigma}(\partial U,\partial V)+\epsilon.
$$ 
Assume $x_k\rightarrow x_\infty$ and $y_k\rightarrow y_\infty$. By 2), we get
$$
\varliminf_{k\rightarrow+\infty}d_{g_k,\Sigma}(x_k,y_k)\geq d_{g,\Sigma}(x_\infty,y_\infty)\geq d_{g}(\partial U,\partial V),
$$
which implies that
$$
\varliminf d_{g_k,\Sigma}(\partial U,\partial V)+\epsilon\geq  d_{g}(\partial U,\partial V)+\epsilon.
$$
Letting $\epsilon\rightarrow 0$, we get  3).

Lastly, we prove 4): under the assumptions of 3), $g_k'$ satisfies 1)-4) of Proposition \ref{ghost}, then 
$$
\lim_{r\rightarrow 0}\lim_{k\rightarrow+\infty} {\rm diam} (D_r,g_k')=0,
$$
which implies that
$$
\lim_{r\rightarrow 0}\lim_{k\rightarrow+\infty}{\rm diam} (B_r^{h_0}(p),d_{g_k,\Sigma})=0.
$$

Finally,  we can treat the case for non-orientable $\Sigma$ by going to the orientation double covering of $\Sigma$.
\endproof

\vspace{.2cm}

\noindent {\it Proof of Theorem \ref{main3}}. 
The proof follows from that of Theorem \ref{main} by replacing Proposition \ref{conv0},  Corollary \ref{conv.length} (i) (ii), Proposition \ref{ghost} therein with Lemma \ref{conv.conf} 1), 2), 3), 4) respectively. So we omit the repetition.
\endproof
\begin{rem}Suppose another subsequence converge to $u',h'$.  Then $\K_{e^{2u}h}=\K_{e^{2u'}h'}$, but now, $h,h'$ may have different conformal classes so the argument in Step 5 in the proof of Theorem \ref{main-ChenLi} is no longer valid.  For varying conformal class case we only conclude convergence of subsequences. 
\end{rem}
Theorem \ref{main3} has immediate consequences.  
Before stating our results, we comment on that for convergence of metrics we may need pulling back the metrics via diffeomorphisms of the surface, and they pull back other geometric quantities such as the conformal factor $e^{2u}$,  the curvature measure $\K_g$ (even nonsmooth) and the distance function as well.  We will abbreviate by using the term ``up to diffeomorphisms'' in statements involving pullbacks, especially for the reason that we will only deal with nondegenerating sequence of conformal classes.   

\begin{cor}\label{Riemann surface}
Let $\Sigma$ be a closed surface, and $h_k$ be a smooth metric with $K_{h_k}=-1$ or $0$ for each $k\in\mathbb N$. When $K_{h_k}=0$ we 
assume $\vol(h_k)=1$. Let $g_k=e^{2u_k}h_k\in\mathcal{M}(\Sigma,h_k)$ and $|\K_{g_k}|(\Sigma)<C$. Assume $d_{g_k,\Sigma}$ converges uniformly to a distance function $d$.
Then, after passing to a subsequence and up to composing diffeomorphisms of $\Sigma$, $h_k$ converges smoothly to a metric $h$ and  
$d_{e^{2u}h,\Sigma}=d$. Further, $u_k$
converges weakly to a function $u$ in $W^{1,q}$ for any $1\leq q<2$.
\end{cor}
\proof
It suffices to show that the conformal classes $\{c_k\}$ of $\{h_k\}$
is sequentially compact in the moduli space.  
We fix a smooth metric $g$ and define $r_0=\inf_{d_{g,\Sigma}(x,y)=a}d(x,y)$, where 
$a$ is smaller than the injectivity  radius of $(\Sigma,g)$.
Then for any $x$, there exists a simply connected domain
$U_x$, such that $B_{r_0}^d(x)\subset U_x$. 
If there was a subsequence of $\{c_k\}$ converges to the boundary
of the moduli space.  by \eqref{vol.com} and Lemma \ref{Mumford}, we can select
 $\gamma_k$ which is nontrivial in $\pi_1(\Sigma)$ with $\ell_{g_k}(\gamma_k)\rightarrow 0$.  Take $x_k\in\gamma_k$
and assume $x_k\rightarrow x_0$.
Since $d_{g_k,\Sigma}\to d$,  $B_{r_0/2}^{g_k}(x_k)\subset
B_{r_0}^d(x_0)\subset U_{x_0}$. Then $\gamma_k\subset U_{x_0}$. A contradiction. \endproof

%\subsection{Contractbility radius}

In \cite{Debin}, the contractibility radius at $x$ is defined as
$$
\co (\Sigma,g,x)=\sup\left\{r > 0 \, \big| \,\overline{B(x,s)} \mbox{ is homeomorphic to a closed disc for every } s < r\right\}
$$
and
$$
\co(\Sigma,g)=\inf_x \co (\Sigma,g,x).
$$
It is shown in \cite{Debin} that when $\Sigma$ is closed, if $\mbox{Area}(\Sigma,g_k)<C$, $\co(\Sigma,g_k)>\delta$, $\K^{+}_{g_k}(B_r(x))<2\pi-\epsilon,$ for some $C,\delta,\epsilon,r>0$ and all $x,k$, then $d_{g_k,\Sigma}$ converges subsequentially, up to pullbacks by diffeomorphisms of $\Sigma$, to a metric $d_{g,\Sigma}$ in $C^0$.  The lower bound of $\co(\Sigma,g)$ guarantees that for a fixed small
$r$ there are converging conformal parametrizations on $B_r^{g_k}(x)$ and this can be used to construct a limiting metric.

We will prove that $\co(\Sigma,g_k)>\delta$ implies convergence of the conformal classes (hence up to composing diffeomorphisms of $\Sigma$, a sequence of the constant curvature metrics converges to a limiting metric $h_\infty$) and it also yields a uniform lower bound of the average of $u_k$ on $(\Sigma,h_\infty)$. Consequently, the assumption $\K_{g_k}(B_r(x))<2\pi-\epsilon$ can be removed from the assumptions in \cite[Main theorem]{Debin}.  

\begin{cor}\label{Debin-improve}
Let $\Sigma$ be a closed surface, and $h_k$ be a smooth metric with Gauss curvature $\K_{h_k}\in \{-1,0,1\}$. When $K_{h_k}=0$, we 
assume $\vol(h_k)=1$. Let $g_k=e^{2u_k}h_k\in \mathcal{M}(\Sigma,h_k)$ and $|\K_{g_k}| (\Sigma)<C$. 
We assume ${\rm diam} (\Sigma,g_k)=\pi$ and $\co(\Sigma, g_k)>\delta>0$.
Then, after passing to a subsequence and up to pulling back by diffeomorphisms $\phi_k:\Sigma\to\Sigma$,  we have 
\begin{enumerate}
\item[1)] $h_k$ converges to a metric $h_\infty$ in $C^\infty$ and $u_k$ converges to some $u_\infty$
weakly in $W^{1,q}(\Sigma,h_\infty), \forall q\in[1,2)$. 
\item[2)] $\K_{g_k}$ converges to $\K_{g_\infty}$ in the sense of distributions where $g_\infty=e^{2u_\infty}h_\infty$.
\item[3)] $d_{g_k,\Sigma}$ converges to $d_{g_\infty,\Sigma}$ in $C^0(\Sigma\times\Sigma)$.
\end{enumerate}
\end{cor}

\proof 
Since ${\rm diam}(\Sigma,g_k)=\pi$ and $|\K_{g_k}|(\Sigma)<C$, by 
\eqref{volume.comparison}
$
\mbox{Area}(\Sigma,g_k)<C'.
$

First, we consider the case that $\mathbb{K}_{h_k}\in\{-1,0\}$.
By Lemma \ref{Mumford}, if the conformal classes converge to the boundary in the moduli space, then we can find a nontrivial circle $\gamma_k$ with $\ell_{g_k}(\gamma_k)\rightarrow 0$. Then $\co(\Sigma,g_k)\rightarrow 0$, which contradicts our assumption. Hence we may
assume $h_k\to h_\infty$ smoothly in a bounded domain of the moduli space.
Then by Lemma \ref{global.gradient.estimate} $$r^{2-q}\int_{B_r^{h_\infty}(x)}|\nabla u_k|^qdV_{h_\infty}<C(q).$$

Let $\K_{g_k}$ converge to a signed Radon measure $\mu$ weakly. 
Since 
$$
\int_{\{u_k>0\}}e^{2u_k^+}dV_{h_\infty}\leq C\int_{\Sigma}e^{2u_k}dV_{h_k}=C\mbox{Area}(\Sigma,g_k)<C,
$$ 
then by Jensen's inequality,  
$$
\int_\Sigma u_k^+ dV_{h_\infty}=\int_{\{u_k\geq 0\}} u_k^+ dV_{h_\infty}<C.
$$

Let $c_k$ be the mean value of $u_k$  on $(\Sigma,h_\infty)$.  Since Area$(\Sigma,h_\infty)=-2\pi\chi(\Sigma)$ when $K_{h_k}=-1$ by the Gauss-Bonnet theorem or equals $1$ when $K_{h_k}=0$ by assumption,  the sequence $c_k$ is bounded from above. By the Poincar\'e inequality,  $\|u_k-c_k\|_{W^{1,q}(\Sigma,h_\infty)}<C$.  We can extract a subsequence (still use the same notation) so that  $u_k-c_k\rightharpoonup u'$ in $W^{1,q}(\Sigma,h_\infty)$.   Setting $g_k'=e^{2(u_k-c_k)}h_k$ and $g'=e^{2u'}h_\infty$, then $\mu=\K_{g'}$.

Since 
$\K_{g_k}=\K_{g_k'}$, we can find a curve $\gamma$ which is
nontrivial in $\pi_1(\Sigma)$, such that $|\K_{g_k'}|(U)<\tau_0$ in a domain $U\supset\gamma$.  By Lemma \ref{Sobolev.distance} and Lemma \ref{conv.conf} 0), 
$$
\|d_{g_k',\Sigma}\|_{W^{1,q}(U\times U,h_\infty)}\leq C\|d_{g_k',\Sigma}\|_{W^{1,q}(U\times U,h_k)}<C. 
$$
Then $a_k=\sup_{x,y\in\gamma}d_{g'_k,\Sigma}(x,y)\rightarrow a$. Fix an $x_0\in\gamma$. Since $\gamma\subset B^g_{a_k}(x_0)$, 
$$
\co(\Sigma,g_k)\leq e^{c_k} a_k
$$
in turn,  $c_k$ is bounded below as well. Hence we may assume $u_k$ is weakly convergent in $W^{1,q}(\Sigma,h_\infty)$.  So, we finish step 1 of 
the proofs of Theorems \ref{main} and \ref{main3}.
%$\ell_{g_k}(\gamma)=O(e^{c_k})$. Since $\ell_{g_k}(\gamma)\geq \co(g_k)$,  $c_k$ is bounded from below, which yields the claim.

Assume $|\K_{g_k}|$ converges to $\nu$ weakly, and
let 
$$
\S=\{x:\nu(\{x\})\geq\tau_0\}=\{p_1,\cdots,p_m\}.
$$
Let $u$ be the weak limit of $u_k$ in $W^{1,q}(\Sigma,h_\infty)$ and set $g=e^{2u}h_\infty$. By Lemma \ref{conv.conf} 0), we may assume $d_{g_k,\Sigma}\to d$ in $C^0_{\rm loc}((\Sigma\setminus\S)\times(\Sigma\setminus\S))$. 
Using the arguments in step 2 of proofs of Theorems \ref{main} and \ref{main3}, we conclude the continuity of $d_{g,\Sigma}$ on $\Sigma$. Now, to complete the step 3 and step 4 there, i.e. to show $d_{g_k,\Sigma}\to d$ uniformly on $\Sigma$ and  $d=d_{g,\Sigma}$, we only need to check that for any $i$ it holds 
$$
\lim_{r\rightarrow 0}\varlimsup_{k\rightarrow+\infty}{\rm diam}(B_r^{h_\infty}(p_i),d_{g_k,\Sigma})=0.
$$
Assume this was not true. Then in an isothermal coordinates of $(\Sigma,h_\infty)$ around $p_i$, 
$$
\lim_{r\rightarrow 0}\varlimsup_{k\rightarrow+\infty}{\rm diam}(D_r,d_{g_k,\Sigma})=b_0>0.
$$
Then for any sufficiently small $r$, after passing to a subsequence, we may assume
${\rm diam}(D_r,d_{g_k,\Sigma})>b_0/2$. 
By \eqref{length.of.circle}, we may 
choose $r$, such that Lemma \ref{main.lemma.2} holds and $\ell_{g}(\partial D_r)<\epsilon<\min\{d(\partial D_{1/4},\partial D_{1/2}),b_0\}/{100}$.  Then
$$
\sup_{x,y\in\partial D_r}d(x,y)\leq\sup_{x,y\in\partial D_r}d_{g,\Sigma}(x,y)\leq \frac{1}{2}\ell_{g}(\partial D_r)<\frac{1}{2}\epsilon.
$$
Then $b_k:=\sup_{x,y\in\partial D_r}d_{g_k,\Sigma}(x,y)<\epsilon$ when $k$ is large.  Hence, there must be a point $x_k\in D_r$ with $d_{g_k}(x_k,\partial D_r)>b_0/8$. 
Take a point $y\in\partial D_r$. Then 
$\partial D_r\subset B^{g_{k}}_{2b_k}(y)\subset D_{1/2}\backslash \{x_k\}$ is not trivial in $B_{2b_k}^{g_k}(y)$ since it is not contractible in $D_{1/2}\backslash \{x_k\}$, so
$
\co(\Sigma,g_k)\leq \epsilon;
$
this is a contradiction if we select $\epsilon<\delta$.

Next, when $K_{h_k}=1$ we can assume $h_k=h= \mbox{the round metric on } \mathbb S^2$.  By composing a M\"obius transformation (specify three points), we assume 
$d_{g_k,\mathbb{S}^2}(N,S)=\pi$ where $N,S$ are the north and south poles,  and there exists $x_k$ so that 
$$
d_{g_k,\mathbb{S}^2}(x_k,N)=d_{g_k,\mathbb{S}^2}(x_k,S)=\frac{\pi}{2}.
$$
We assume $x_k\to x_0$. 

Let $c_k$ be the mean value of $u_k$ on $({\mathbb S}^2,h)$.  We can extract a subsequence so that  $u_k-c_k\rightharpoonup u'$ in $W^{1,q}(\mathbb S^2,h)$.  

We will show $\inf_kc_k>-\infty$ by contradiction.  Suppose $c_k\to -\infty$. We set $g_k'=e^{2(u_k-c_k)}h$, $g'=e^{2u'}h$.  Since $|\K_{g_k}(\mathbb S^2)|<C$,  we assume $|\K_{g_k'}|=|\K_{g_k}|$ converges to $\nu$ weakly (up to a subsequence).  For the finite set 
$$
\S=\{x:\nu(\{x\})\geq\tau_0\}
$$
we can select a small $r$ such that $\Gamma:=\partial B^h_r(x_0)\cap\S=\emptyset$. Then,  by Corollary \ref{MTI} and the trace embedding theorem,  (up to a subsequence) $\delta_k=\ell_{g_k}(\Gamma)=e^{c_k}\ell_{g_k'}(\Gamma)\rightarrow 0$.  Fix a point $p\in\Gamma$. Since $d_{g_k,\mathbb{S}^2}(p,S)+d_{g_k,\mathbb{S}^2}(p,N)\geq d_{g_k,\mathbb{S}^2}(S,N)=\pi$, without loss of generality, we may assume $d_{g_k,\mathbb{S}^2}(p,S)\geq\pi/4$. 

We have two cases. 

Case 1: $d_{g_k,\mathbb{S}^2}(x_k,p)\geq a>0$ for all $k$. For this case $\Gamma\subset B_{2\delta_k}^{g_k}(p)\subset\mathbb{S}^2\backslash\{S,x_k\}$. Since $\Gamma$ is not contractible in $\mathbb{S}^2\backslash\{S,x_0\}$,  we would have $\co(\mathbb S^2,g_k,p)\leq 2\delta_k$,  which contradicts the assumption $\co(\mathbb S^2,g_k)>\delta$. 

Case 2: by passing to a subsequence, $d_{g_k,\mathbb{S}^2}(x_k,p)\rightarrow 0$.  For this case we have 
$$
d_{g_k,\mathbb{S}^2}(p,S),\s d_{g_k,\mathbb{S}^2}(p,N)\geq \frac{\pi}{4}.
$$
If the loop $\Gamma$ is not contractible in $\mathbb{S}^2\backslash \{S,N\}$, then 
$\co(\mathbb{S}^2,g_k,p)\leq 2\delta_k$, which leads to a contradiction. If $\Gamma$ is contractible in $\mathbb S^2\backslash\{S,N\}$, we select a loop $\Gamma'$ passing through $p$ in $\mathbb S^2\backslash(\S\cup\{S,N\})$ which is not contractible in $\mathbb S^2\backslash\{S,N\}$. Then $\delta'_k=\ell_{g_k}(\Gamma')=e^{c_k}\ell_{g'_k}(\Gamma')\to 0$, in turn, $\co(\mathbb S^2,g_k,p)\leq 2\delta'_k$, a contradiction. 

\vspace{.1cm}

The rest of the proof is the same as that for the case $\K_{h_k}\in\{-1,0\}$. 
\endproof

\section{Appendix}

%\begin{lem}
%Then for almost every $s,t\in(R_1,R_2)$,  $s<t$,  we have 
%$$
%t\frac{du^*}{dr}(t)- s\frac{du^*}{dr}(s)=-\frac{1}{2\pi}\,\mu(D_t\backslash {D}_s).
%$$
%There exists $E\subset[0,R]$ with ${\mathscr L}^1(E)=0$, such that for any $\{t_k\}\subset[0,R]\backslash E$ with $t_k\rightarrow 0$ it holds
%$$
%\mu(\{0\}) =-2\pi\lim_{t_k\to 0}t_k\frac{du^*}{dr}(t_k).
%$$
%\end{lem} 

\subsection{Sharpness of the curvature threshold}  

We construct a sequence of metrics on $\mathbb S^2$ that satisfies the assumptions in Theorem \ref{main-ChenLi} except $\mu^1(\{x\})<2\pi$ for all $x\in \mathbb S^2$ and show the distances do not converge to a distance function.

First, we recall that Hulin-Troyanov constructed \cite{Hulin-Troyanov} a metric on ${D_{1/2}}\subset\C$ by
$$
g=\frac{|dz|^2}{|z|^2|\log|z||^{2a}}.
$$
Its Gauss curvature measure is 
$$
\K_g=2\pi\delta_0-a|\log|z||^{2a-2}dV_g.
$$
If $a>0$, the singular point $0$ is a cusp and $K<0$ for $z\not=0$.  If $a>1$ the cusp is of finite distance from any $z\not=0$ but complete; if $0<a\leq 1$ the cusp is of infinity distance so $g$ is complete; the area is finite if $a> 1/2$ and infinite if $a\leq 1/2$.
The special case $a = 1$ is the Beltrami pseudosphere so $\K_g=2\pi\delta_0-dV_g$ and the limit case $a=0$ is a half cylinder so $\K_g=2\pi\delta_0$. Corollary \ref{complete} applies to all of these cases for $(D\backslash\{0\},g)$. %Below, we demonstrate that  Theorem \ref{main-ChenLi} is not true if the key assumption $\mu^1(\{x\})<2\pi$ is dropped. 

Fix an $a\in(\frac{1}{2},1)$ and extend the above 
$g$ to a smooth metric $\tilde{g}$ on $\C\backslash\{0\}$ by taking $
\frac{|dz|^2}{(1+\frac{1}{4}|z|^2)^2}$ on $D^c$ and smoothly joining it with $g$ on $D\backslash D_{1/2}$.
As $(\C,\frac{|dz|^2}{(1+\frac{1}{4}|z|^2)^2})$ is isometric to ${\mathbb S}^2\backslash \{\mbox{a point}\}$, $\tilde{g}$ can be viewed as a (singular) metric on $\mathbb S^2$.
Let $\eta\geq 0$ be a $C^\infty$ function on $\R$ which is 1 in $(-\infty,\frac{1}{8})$ and $0$ in $(\frac{1}{4},+\infty)$. Define
$
g_\epsilon=e^{2\epsilon\eta(|z|) \log|z|}\tilde{g}.
$
Then 
$$
\K_{g_\epsilon}=\K_g-2\pi\epsilon\delta_0+\epsilon F|dz|^2,
$$
for some $F\in C^\infty_0(D_\frac{1}{4}\backslash D_\frac{1}{8})$.
Set $g_\epsilon=e^{2u_\epsilon}g_{\mathbb S^2}$ and $g_\epsilon'=e^{2(u_\epsilon-c_\epsilon)}g_{\mathbb S^2}$, where $c_\epsilon$ is a chosen constant so that $\mbox{diam}\,(\mathbb S^2,d_{g_\epsilon'})=1$. Since $d_g(0,x)=+\infty$ for any $x\neq 0$, we see $c_\epsilon\rightarrow+\infty$.  As $\epsilon\rightarrow 0$, $\K_{g_\epsilon'}=\K_{g_\epsilon}$ converges to $\K_g$ weakly, $u_\epsilon-c_\epsilon$ converges
to $-\infty$ almost everywhere  and $d_{g_\epsilon'}$ converges to 0 on any compact subset of $\mathbb C\backslash\{0\}$. But the zero function cannot be a distance. Near $0\in D$,
 we can write $\K_{g_\epsilon'}$
as 
\begin{align*}
\K_{g_\epsilon'}=\left(2\pi(1-\epsilon)\delta_0+\epsilon F^+|dz|^2\right)-\left(a|\log|z||^{2a-2}dV_g+\epsilon F^-|dz|^2\right)
:=\mu_\epsilon^1-\mu_\epsilon^2.
\end{align*}
If we let $\mu_\epsilon^1\rightharpoonup\mu^1$ and $\mu_\epsilon^2\rightharpoonup\mu^2$ as $\epsilon\rightarrow 0$, then $\mu^1(\{0\})=2\pi$. 
%Hence, the assumption `$\mu^1(\{x\})<2\pi$ for any $x$' in Theorem \ref{main-ChenLi} cannot be dropped.

\vspace{.15cm}

\begin{comment}
\noindent ({\bf D}) $(D, g)$ with $g = e^{2u(z)}|dz|^2$ where $u$ is the real part of a holomorphic function with isolated singularity of order $k$ at 0. The real part of $ z^{-k}$ with $k>0$ is the harmonic function $\frac{\cos k\theta}{r^k}$ for $r>0$. The corresponding $u\not\in W^{1,1}_{\rm loc}(D)$. Nevertheless, $u\in W^{1,1}_{\rm loc}(D\backslash\{0\})$ and 
$\K_g (A) = 0$ for any $A\subseteq D\backslash\{0\}$. Note that $0$ is an essential singularity of $e^{1/z^{k}}$.  The metric $g$ is incomplete as seen in the direction $k\theta=\pi$. Corollary \ref{complete} does not apply to $(D\backslash\{0\},g)$.
\end{comment}

%If we can get a gradient estimate of $u^{har}$, then a $L^q$-gradient estimate of $u$ follows. 

\subsection{Proof of Lemma \ref{Gauss-Bonnet-Green}}
We divide our arguments into four steps. 

\vspace{.15cm}

\noindent {\bf Step 1.} We show that $u^*$ is absolutely continuous. 
Choose a sequence $u_k\in C^\infty(\R^2)$ converging to $u$ in $W^{1,1}(D_{R_2}\backslash D_{R_1})$. By the trace embedding theorem, 
$$
\|\hbox{tr}(u_k-u)\|_{L^1(\partial B_{R_2}\cup \partial{B_t})}\leq C\|u_k-u\|_{W^{1,1}(B_{R_2}\backslash B_t)}\rightarrow 0.
$$
In particular, $u^*_k(t)\rightarrow u^*(t)$ for any $t$.  Since $u_k$ is smooth, for any $R_1\leq s\leq t\leq R_2$, 
%$$u_k^*(t)-u_k^*(s)=\int^t_s\frac{du_k^*}{dr}dr=\frac{1}{2\pi}\int_{D_{t}\backslash D_s}\frac{\partial u_k}{\partial r}\frac{1}{r}dx.$$
\begin{align*}
u^*(t)-u^*(s)&=\lim_{k\to\infty}\left(u^*_k(t)-u^*_k(s)\right) =\lim_{k\to\infty} \frac{1}{2\pi}\int_{D_t\backslash D_s}\frac{\partial u_k}{\partial r}\frac{1}{r}dx\\
&=\frac{1}{2\pi}\int_{D_t\backslash D_s}\frac{\partial u}{\partial r}\frac{1}{r}dx=\frac{1}{2\pi}\int_{s}^t\int_{0}^{2\pi}\frac{\partial u}{\partial r}d\theta dr.
\end{align*}
By Fubini's Theorem, $\int_0^{2\pi}\frac{\partial u}{\partial r}d\theta\in L^1([R_1,R_2])$. So 
$u^*$ is absolutely continuous on $[R_1,R_2]$ with 
\begin{equation}\label{in proof}
\frac{du^*}{dr}(r)=\frac{1}{2\pi}\int_0^{2\pi}\frac{\partial u}{\partial r}d\theta,\s  \s \mbox{a.e.}\s r\in[R_1,R_2].
\end{equation}
As $u_k\rightarrow u$ in $W^{1,1}(D_{R_2}\backslash D_{R_1})$, Fubini's Theorem asserts 
$$
\int_{R_1}^{R_2}\left|\int_0^{2\pi}\left(\frac{\partial u_k}{\partial r}-\frac{\partial u}{\partial r}\right)d\theta\right|dr\rightarrow 0.
$$
We may therefore assume $$\frac{du_k^*}{dr}(r)\to\frac{du^*}{dr}(r)\s\s\s\mbox{a.e. $r\in[R_1,R_2]$.}
$$

\vspace{.15cm}

\noindent{\bf Step 2.} Denote  
$
\mathcal K(R_1,R_2,r)=\left\{\varphi\in C^\infty_0(D_{R_2}\backslash D_{R_1}): {\hbox{ $\varphi=1$ on $\partial D_r$}}\right\}
$  a set of test functions for $R_1<r<R_2$. 
For $\varphi\in\mathcal K(R_1,R_2,r)$, define 
\begin{equation}\label{Gauss-Bonnet1}
\lambda(r,u,\varphi,R_1)=\frac{1}{2\pi}\int_{D_{r}\backslash D_{R_1}}\nabla\varphi\nabla u\,dx-\frac{1}{2\pi}\int_{D_{r}\backslash D_{R_1}}\varphi \,d\mu.
\end{equation}

Since $u$ solves $-\Delta u=\mu$ in $D_{R_2}\backslash D_{R_1}$ we know $u\in W^{1,q}$, we claim that $\lambda(r,u,\varphi,R_1)$ is independent of $\varphi\in\mathcal K(R_1,R_2,r)$. To verify this, take another $\varphi'$ there. Then 
$\varphi-\varphi'$ is 0 on the boundary of $D_r\backslash D_{R_1}$ and the support of $\varphi-\varphi'$ does not intersect $\partial D_{R_1}$. Now,
$$
\int_{D_r\backslash D_{R_1}}\nabla(\varphi-\varphi')\nabla u \,dx=\int_{D_r\backslash D_{R_1}}(\varphi-\varphi')\,d\mu,
$$
in turn
$
\lambda(r,u,\varphi,R_1)=\lambda(r,u,\varphi',R_1).
$
This leads to
$
\lambda(r,u,\varphi,R_1)=\lambda(r,u,\varphi,R_1')
$
for any $R_1'\in (R_1,r)$ because if $\varphi\in\mathcal K(R_1',R_2,r)$ then $\varphi\in\mathcal K(R_1,R_2,r)$. So we will simply write $\lambda(r,u)$. 

\vspace{.1cm}

For $u\in L^1(D)$ with $\Delta u$ as a signed Radon measure, as in Lemma \ref{local} we can select smooth $u_k,f^1_k,f^2_k$ with $-\Delta u_k=f_k^1-f_k^2$,
$f_k^idx\rightharpoonup\mu^i$, $\mu_k=(f_k^1-f_k^2)dx\rightharpoonup\mu$. Moreover,  $\frac{d u_k^*}{dr}(r)\to\frac{d u^*}{dr}(r)$ a.e. $r$, $\mu=\mu^1-\mu^2$ and $|\mu|\leq \mu^1+\mu^2$.

We claim that there is a countable set $A\subset(R_1,R_2)$ such that for any $r\notin A$ and $\varphi\in \mathcal{K}(R_1,R_2,r)$
$$
\int_{D_r\setminus D_{R_1}}\varphi d\mu_k\rightarrow \int_{D_r\setminus D_{R_1}}\varphi d\mu.
$$
Define
$
A=\left\{r\in(R_1,R_2):(\mu^1+\mu^2)(\partial D_r)>0\right\}.
$
Since $(\mu^1+\mu^2)(D_r)<\infty$, $A$ is a countable set.  Fix $r\in (R_1,R_2)\backslash A$. Let $t\in (R_1,r)\backslash A$ and $\eta_{t}:\R\rightarrow [0,1]$ be a cut-off function which is 1 on $[-\infty,t]$ and 0 on $[r,+\infty)$. By 
\cite[Theorem 1.40 (iii)]{E-G}, 
$$
\lim_{k\rightarrow\infty}\int_{D_r\backslash D_t}(f_k^1+f_k^2)dx=(\mu^1+\mu^2)(D_r\backslash D_t).
$$
Then
\begin{align*}
\varlimsup_{k\rightarrow\infty}&\left|\int_{D_r\backslash D_{R_1}}
\varphi d\mu_k-\int_{D_r\backslash D_{R_1}}\varphi d\mu\right|\leq
\varlimsup_{k\rightarrow\infty}\left(\left|\int_{D_{R_2}\backslash D_{R_1}}
\eta_t\varphi d\mu_k-\int_{D_{R_2}\backslash D_{R_1}}\eta_t\varphi d\mu\right|\right.\\
&+\left|\int_{D_r\backslash D_{R_1}}\left.
(1-\eta_t)\varphi d\mu_k-\int_{D_r\backslash D_{R_1}}(1-\eta_t)\varphi d\mu\right|\right)\\
&\leq
C\|\varphi\|_{C^0}\varlimsup_{k\rightarrow\infty}\left(\int_{D_r\backslash D_t}(f_k^1+f_k^2)dx
+|\mu|(D_r\backslash D_t)\right)\\
&=C\|\varphi\|_{C^0}(\mu^1+\mu^2+|\mu|)(D_r\backslash D_t).
\end{align*}
Letting $t\rightarrow r$, we get
$$
\varlimsup_{k\rightarrow+\infty}\left|\int_{D_r\backslash D_{R_1}}
\varphi d\mu_k-\int_{D_r\backslash D_{R_1}}\varphi d\mu\right|=0.
$$
For smooth $u_k$, $\lambda(r,u_k)=r\frac{du^*_k}{dr}$. It follows from \eqref{Gauss-Bonnet1} that
$\lambda(r,u_k)\rightarrow\lambda(r,u)$ a.e.  $r$. Then
\begin{equation}\label{a.e.}
\lambda(r,u)=r\frac{d u^*(r)}{dr},\s
\hbox{ a.e.\, $r$}.
\end{equation}

\noindent{\bf Step 3.}
Let $\varphi\in\mathcal K(R_1,R_2,t)$ such that $\varphi|_{D_t\setminus D_s}=1$. Then $\varphi$ is also in
$\mathcal K(R_1,R_2,s)$. Hence
\begin{align*}
 &t\frac{du^*}{dr}(t)\, - s\frac{du^*}{dr}(s)= \lambda(t,u)-\lambda(s,u) \\
&=\frac{1}{2\pi}\int_{D_{t}\setminus D_{R_1}}\nabla\varphi\nabla u\,dx-\frac{1}{2\pi}\int_{D_{t}\setminus D_{R_1}}\varphi \,d\mu-\frac{1}{2\pi}\int_{D_{s}\setminus D_{R_1}}\nabla\varphi\nabla u\,dx+\frac{1}{2\pi}\int_{D_{s}\setminus D_{R_1}}\varphi \,d\mu\\
&=-\frac{1}{2\pi}\int_{D_{t}\setminus D_{s}}\varphi \,d\mu= -\frac{1}{2\pi}\mu(D_t\backslash D_s).
\end{align*}
Then for almost every $s,t\in[R_1,R_2]$ with $s<t$ 
\begin{eqnarray*}
t\frac{du^*}{dr}(t)-s\frac{du^*}{dr}(s)&=& \lambda(t,u)-\lambda(s,u) = -\frac{1}{2\pi}\mu(D_t\backslash D_s).
\end{eqnarray*}
This proves the first statement in the lemma. 

\vspace{.1cm}

\noindent{\bf Step 4.}
When $-\Delta u = \mu$ on a disk $D_R$, set 
$\mathcal K'(R,r)=\left\{\varphi\in C^\infty(D_{R}): {\hbox{ $\varphi=1$ on $\partial D_r$}}\right\}$.
For $\varphi\in\mathcal K'(R,r)$, define 
\begin{equation*}%\label{Gauss-Bonnet2}
		\lambda'(r,u,\varphi)=\frac{1}{2\pi}\int_{D_{r}}\nabla\varphi\nabla u\,dx-\frac{1}{2\pi}\int_{D_{r}}\varphi \,d\mu.
\end{equation*}
As in Step 2, $\lambda'(r,u,\varphi)$ is independent of the choice of $\varphi$, we will write $\lambda'(r,u)$. 
The same arguments justify  \eqref{a.e.}  for $\lambda'(r,u)$. 
%\begin{eqnarray}\label{a.e.'}
%\lambda'(r,u)=r\frac{d u^*(r)}{dr},\s \mbox{a.e.}\s r.
%\end{eqnarray}
Letting $\varphi=1$, we get
\begin{equation*}\label{GB}
\lambda'(r,u)=-2\pi\mu(D_r).
\end{equation*}

Let $r_k$ go to 0 with \eqref{a.e.} holds. Since
$$
\lim_{k\rightarrow+\infty}
|\mu(D_{r_k})-\mu(\{0\})|\leq\lim_{k\rightarrow+\infty}
|\mu|(D_{r_k}\backslash\{0\})=
|\mu|(\cap_kD_{r_k}\backslash\{0\})=|\mu|(\emptyset)=0,
$$
we conclude the proof with 
$$
\mu(\{0\}) = \lim_{r_k\to 0 }\mu(D_{r_k})=-2\pi\lim_{r_k\to 0}r_k\frac{du^*}{dr}(r_k). 
$$
%\endproof
%$$\lim_{r\rightarrow 0}\lambda(r,u)=-2\pi\mu(\{0\}).$$

\subsection{Hausdorff measure of sets defined by Rayleigh quotient} The goal of this section is to prove Lemma \ref{diff.of.d.1}. The special case $n=2,s=1$ is used to establish the distance comparison theorem  in section 3.2.

Denote $B_r(x)$ the ball in $\R^n$ and write $B_r$ for $B_r(0)$. For $u\in W^{1,p}(B_2), 1\leq p<n$, let $u_{x,r}$ be the average of $u$ over $B_r(x)$
%we set $$ u_{x,r} =\frac{1}{\omega_nr^n}\int_{B_r(x)} u(y)\, dy$$ where $\omega_n$ is the volume of $B_1$ and $B_r(x)\subset B_2$, 
and 
$$
%A(u)=\{x\in B_1:\lim_{\tau\rightarrow 0}\underset{r\in(0,\tau]}{\hbox{osc}} u_{x,r}>0\}
A(u)= \left\{x\in B_2: \lim_{r\to 0^+}u_{x,r}\, \hbox{does not exist or ${\limsup}_{r\to 0^+}|u_{x,r}|=\infty$}\right\}.
$$
According to a theorem of Federer and Ziemer (\cite{Fed-Zie},  also see \cite[Theorem 2.1.2]{Lin-Yang}), 
the Hausdorff dimension of $A(u)$ satisfies
$
\dim_{\mathcal H} A(u)\leq n-p.
$
For any $x\notin A(u)$, we define 
$$
\hat{u}(x)=\lim_{r\rightarrow 0} u_{x,r}.
$$
Note that $\hat{u}$ is well-defined for $\mathcal{H}^s$-a.e. $x\in B_1$, for $s> n-p$. The Sobolev function $u$ can be altered over an $\mathcal H^s$ measure zero set without changing its $W^{1,p}$ norm.  We always assume the alternation is done, namely, $\hat{u}(x) = u(x)$, $\mathcal H^s$-a.e.

%In order to get the further properties about $\hat{u}$, we first consider decay of $u_{x,r}$ as $r\to 0$. 
\begin{lem}\label{decay}
Suppose $u\in W^{1,p}(B_{r_0})$ with
$$   
\frac{1}{r^s}\int_{B_r}|\nabla u|^p<M,\s \forall r<r_0.
$$
Then for any $r_1<r_0, s\in(n-p,n]$, we have
$$
|u_{0,r_0}-u_{0,r_1}|\leq \Lambda M^\frac{1}{p}r_0^\theta,
$$
where $\Lambda=\Lambda(n,s,p)$ and $\theta=\frac{p-n+s}{p}>0$.
\end{lem} 

\proof
Recall the Poincar\'e inequality
$$
\frac{1}{|B_r|}\int_{B_r}|u-u_{0,r}|^p\leq
\Lambda_1r^{p-n}\int_{B_r}|\nabla u|^p,
$$
where $\Lambda_1$ only depends on $n$. This gives
\begin{align}\label{repeat}
|u_{0,\frac{r}{2}}-u_{0,r}|&=\frac{1}{|B_{\frac{r}{2}}|}\left|\int_{B_\frac{r}{2}}(u-u_{0,r})\right|\nonumber\leq \frac{1}{|B_{\frac{r}{2}}|}\int_{B_r}\left|u-u_{0,r}\right|\nonumber\\
&\leq\frac{1}{|B_{\frac{r}{2}}|}\left(\int_{B_r}|u-u_{0,r}|^p\right)^\frac{1}{p}|B_r|^{1-\frac{1}{p}}=2^n\left(\frac{1}{|B_r|}\int_{B_r}|u-u_{0,r}|^p\right)^\frac{1}{p}\nonumber\\
& \leq 2^n\left(\Lambda_1r^{p-n}\int_{B_r}|\nabla u|^p\right)^\frac{1}{p}  \leq\Lambda_2\,r^\theta M^\frac{1}{p}, 
\end{align}
where $\theta=\frac{p-n+s}{p}$ and $\Lambda_2=2^n\Lambda_1^{{1}/{p}}$.

Assume $r_1\in[2^{-k}r_0,2^{-k+1}r_0)$ where $k\in\mathbb{N}$.  We have
\begin{equation}\label{triangle1}
|u_{0,r_0}-u_{0,2^{-k}r_0}|\leq \Lambda_2 M^{\frac{1}{p}}\sum_{i=0}^{k-1} (2^{-i}r_0)^\theta 
\leq \Lambda_3 M^\frac{1}{p}r_0^\theta.
\end{equation}
Repeating the argument in \eqref{repeat} leads to 
\begin{align}\label{triangle2}
|u_{0,2^{-k}r_0}-u_{0,r_1}|&=\frac{1}{|B_{2^{-k}r_0}|}\left|\int_{B_{2^{-k}r_0}}(u-u_{r_1})\right| \leq \frac{1}{|B_{2^{-k}r_0}|}\int_{B_{r_1}}\left|u-u_{r_1}\right| \nonumber\\
&\leq \frac{|B_{r_1}|}{|B_{2^{-k}r_0}|} \left(\frac{1}{|B_{r_1}|}\int_{B_{r_1}} |u-u_{0,r_1}| \right)^{\frac{1}{p}} \nonumber\\
&\leq \Lambda_2 r_1^\theta \ M^\frac{1}{p}\leq \Lambda_2 2^{(1-k)\theta} r_0^\theta \ M^{\frac{1}{p}} \leq \Lambda_2  r_0^\theta \ M^{\frac{1}{p}}. 
\end{align}
The desired result follows from the triangle inequality, \eqref{triangle1} and \eqref{triangle2}.
\endproof

In the proof of Lemma \ref{measure.estimate} below, we will cover $E(u,\lambda)$ with countable balls $\overline{B_{r_i}(x_i)}$; however, we can only do this with $r_i<1$ not  with $r_i<\delta$ for any fixed $\delta$. Thus, we do not have an estimate of $\mathcal{H}^1$. Instead, we use \cite{E-G}: for $A\subset\R^n$
$$
\mathcal{H}^s_\infty(A):=\inf\big\{\sum_{j=1}^\infty\alpha(s)\left(\frac{{\rm diam} C_j}{2}\right)^s| A\subset\bigcup_{j=1}^\infty C_j\big\}.
$$

\begin{lem}\label{measure.estimate}
Suppose $u\in W^{1,p}(B_2)$ with $\|u\|_{L^1(B_2)}\leq \frac{\omega_n}{4}\lambda$ and $ p<n$. Let
$
E(u,\lambda)=\{x\in B_1\backslash A(u): |{u}(x)|>\lambda\}.
$
Then for any $s\in(n-p,n]$ we have
$$
\mathcal{H}^s_\infty(E(u,\lambda))\leq\frac{\Lambda'}{\lambda^p}\int_{B_2}
|\nabla u|^p,
$$
where $\Lambda'=\Lambda'(n,s,p)$. There is 
a cover $\left\{\overline{B_{r_i}(x_i)}\right\}$ of $E(u,\lambda)$ such that 
$x_i\in E(u,\lambda)$ and for $\omega_s= \pi^{\frac{s}{2}}/\Gamma(\frac{s}{2}+1)$ it holds 
$$
\omega_s\sum_i r_i^s\leq\frac{\Lambda'}{\lambda^p}\int_{B_2}|\nabla u|^p.
$$
%where $\omega_s= \pi^{\frac{s}{2}}/\Gamma(\frac{s}{2}+1)$. 
\end{lem}
\proof Let $x\in E(u,\lambda)$. 
Set $\Lambda M^\frac{1}{p}=\lambda/4$ for $r_0=1$ in Lemma \ref{decay}. 
Then, if 
$$   
\frac{1}{r^s}\int_{B_r(x)}|\nabla u|^p<M,\s \forall r<1,
$$
we would have
$$
|u_{x,r}-u_{x,1}|<\Lambda M^\frac{1}{p}=\frac{\lambda}{4}
$$
where $\theta=\frac{p-n+s}{p}\leq 1$. Letting $r\rightarrow 0$, we would get
$$
|{u}(x)-u_{x,1}|\leq\frac{\lambda}{4}.
$$
Then
$$
|{u}(x)|\leq \frac{\lambda}{4}+|u_{x,1}|\leq
\frac{\lambda}{2}+\frac{1}{|B_1(x)|}\|u\|_{L^1(B_2)}<\lambda.
$$
This contradicts $x\in E(u,\lambda)$. Thus, for any $x\in E(u,\lambda)$ 
there exists $r<1$ such that 
$$
\frac{1}{r^s}\int_{B_r(x)}|\nabla u|^p\geq M =\left(\frac{\lambda}{4\Lambda}\right)^p. 
$$
By the Vitali Covering Lemma,  there exists pairwise disjoint $\overline{B_{r_i}(x_i)}$ such that 
$$
\frac{1}{r^s}\int_{B_{r_i}(x_i)}|\nabla u|^p\geq {M},\s
E(u,\lambda)\subset\bigcup_i\overline{B_{5r_i}(x_i)}.
$$
Then,
\begin{equation*}
\mathcal{H}_\infty^s(E(u,\lambda))\leq \sum_i\omega_s(5r_i)^s
\leq
\frac{5^s\omega_s}{M}\int_{\cup B_{r_i}(x_i)}|\nabla u|^p
\leq \frac{5^s\omega_s}{M}\int_{B_2}|\nabla u|^p.
\end{equation*}
\endproof

%We will use the following result for $n=2$ to prove a distance comparison theorem in subsection 3.3.
\begin{lem}\label{diff.of.d.1}
Let $u\in W^{1,p}(B_2)$ with $p\in (1,2)$. Then
for any $\epsilon>0$ there exists $\lambda=\lambda(\epsilon)>0$ such that for any $s>n-p$ it holds 
$$
\mathcal{H}^s_\infty\left(\left\{x\in B_1:|u(x)-u_{0,1}|>\lambda\,{\|\nabla u\|_{L^p(B_2)}}\right\}\right)\leq\epsilon.
$$
\end{lem}
\proof Recall $\dim_{\mathcal H}A(u)\leq n-p<s$ and $\hat u(x)$ exists for $\mathcal H^s$-a.e. $x$. We assume $\|\nabla u\|_{L^p(D_2)} \not=0$ as the lemma is trivially true otherwise. 
By the Poincar\'e inequality, 
$$
\frac{\|u-u_{0,1}\|_{L^1(D_2)}}{\|\nabla u\|_{L^p(D_2)}}<\frac{c\omega_n}{4},
$$
where $c$ is a uniform constant and $\omega_n$ is the volume of $B_1$.  
For any $\epsilon>0$, applying Lemma \ref{measure.estimate} to $$v:=\frac{u-u_{0,1}}{\|\nabla u\|_{L^p(D_2)}}, \ \ \ \ \lambda\geq\max\left\{c,{\left(\frac{\Lambda'}{\epsilon}\right)^{1/p}}\right\},$$ 
we see 
$$
\mathcal H^s_\infty(\{x\in D:|\hat v(x)|>\lambda\}\cup A(u))<\epsilon.
$$
We finish the proof by noting that $v(x)=\hat v(x)$ for $\mathcal H^s_\infty$-a.e. $x$. 
\endproof

%\vspace{-.1cm}

\end{document}